\documentclass[pdflatex,sn-mathphys-num]{sn-jnl}

\usepackage{graphicx}%
\usepackage{multirow}%
\usepackage{amsmath,amssymb,amsfonts}%
\usepackage{amsthm}%
\usepackage{mathrsfs}%
\usepackage[title]{appendix}%
\usepackage{xcolor}%
\usepackage{textcomp}%
\usepackage{manyfoot}%
\usepackage{booktabs}%
\usepackage{algorithm}%
\usepackage{algorithmicx}%
\usepackage{algpseudocode}%
\usepackage{listings}%
\usepackage{lineno}
\usepackage[symbol]{footmisc}
\usepackage{moresize}

\usepackage{dirtytalk}

\newtheorem{theorem}{Theorem}
\newtheorem{lemma}{Lemma}
\newtheorem{prop}{Proposition}
\theoremstyle{definition}

\newtheorem{remark}{Remark}

\newcommand{\rev}[1]{{\color{black} #1}}

\begin{document}

\title[An insight on some properties of high order nonstandard linear multistep methods]{An insight on some properties of high order nonstandard linear multistep methods}


\author[1,2,*]{\fnm{Balint M.} \sur{Takacs}\footnote[0]{*Corresponding author. E-mail: takacs.balint.mate@ttk.bme.hu}}



\affil[1]{\orgdiv{Department of Analysis and Operations Research}, \orgname{Budapest University of Technology and Economics}, \orgaddress{\street{Egry József u. 1.}, \city{Budapest}, \postcode{H-1111},  \country{Hungary}}}

\affil[2]{\orgname{HUN-REN-ELTE Numerical Analysis and Large Networks Research Group}, \orgaddress{\street{Pázmány Péter sétány 1/C}, \city{Budapest}, \postcode{H-1117}, \country{Budapest}}}



\abstract{In this paper, nonstandard multistep methods are \rev{considered}. It is shown that under some (sufficient and necessary) conditions, these methods attain the same order as their standard counterparts - to prove this statement, a nonstandard version of Taylor's series is constructed. The preservation of some qualitative properties (boundedness, the linear combination of the components\rev{, and a property similar to monotonicity}) is also proven for all step sizes. The methods are applied to a one-dimensional equation and a system of equations, in which the numerical experiments confirm the theoretical results.}

\keywords{nonstandard finite difference, multistep method, positivity preservation, SSP linear multistep methods}


\pacs[MSC Classification]{65L06, 92D25, 92D30}

\maketitle


\section{Introduction}\label{sec:intro}
Two key aspects of the numerical modeling of differential equations have always been the high accuracy of the methods along with the preservation of important features of the initial continuous problem. It has been known that the application of high-order methods does not necessarily imply the latter: see e.g. \cite{mickens}, along with the specific cases of multistep methods \cite{hunsd2, hundsdorfer}.

During the previous decades, numerous different approaches have been taken to construct methods with such properties. One \rev{of such approaches}, which has had an increasing prevalence in recent years, has been the approach of nonstandard methods: only last year, multiple articles applied these ideas to epidemic problems \cite{aljohani, buhader, butt, hoang243, maa}, to equations of population dynamics \cite{appadu, hoang242, hoang244, mokni}, to models in microbiology and chemistry \cite{alal24, pez}, to problems in finance \cite{dube, moaddy} and also to equations with derivatives of fractional order \cite{abou, alkat, cai, namjoo}, just to name a few. \rev{For a comprehensive list of earlier applications, the author recommends the remarkably insightful survey articles of Patidar \cite{patidar1, patidar2} and the references therein.} The main advantage of the use of such methods is that they are customarily constructed in a way that some desired properties of the original, continuous model are preserved for any choice of step size.

One thing that most former applications of the nonstandard ideas lack is a discussion regarding the order of the methods: the majority of them achieved only the order of one. Nonetheless, in recent years, there has been an emerging trend for the construction of higher-order nonstandard methods. 

Earliest examples include the application of these ideas to the Yee algorithm \cite{cole97, cole02} and the \rev{analysis} of the second-order nonstandard Runge-Kutta method \cite{dimitrov, dimitrov2}. Nevertheless, nowadays one leading school of thought is the use of theta methods, which can achieve an order of two \cite{ gupta2, kojouharov}, can be combined with second order Runge-Kutta methods \cite{dimitrov3}, and can also be applied to reaction-diffusion \cite{anguelov03, anguelov2} and singularly perturbed problems \cite{lubuma}.

Another possibility is the application of Richardson extrapolation \cite{bassenne}, which can also be applied to population models \cite{gonzalez}. There are also approaches involving the use of a more complicated form of the denominator function \cite{alal, alal23, gupta, hoang24, hoang24new, conte25}, while other ideas incorporate innovative modifications of the right-hand side of the equation \cite{anguelov22, hoang22}. Other works include the general order conditions for one-step nonstandard methods applied to some special types of equations \cite{chen, kojo04} or methods with variable step sizes \cite{martin}. These ideas have also been applied to several real-life applications \cite{alal24, hoang23, hoang232, hoang242}.

It should also be mentioned that \rev{the core idea in the preservation of the properties in the case of nonstandard multistep methods, namely,} the application of strongly stable methods in the context of nonstandard methods, has appeared in some previous works: some have combined the nonstandard methods with total variation diminishing (TVD) schemes \cite{khal14} and Runge-Kutta methods written in a strongly stability preserving (SSP) form \cite{khal17} - the authors in \cite{anguelov3} even got the order condition for second order methods akin to the one proven in Section \ref{sec:order} (other similar instances include \cite{alal} and \cite{hoang24}). 

One article to highlight is the work of Dang and Hoang \cite{dang} involving nonstandard Runge-Kutta methods. In this, it is shown that the nonstandard versions of these schemes have the same order as their standard counterparts, and (under some conditions) they also preserve the properties of the initial differential equation. It should also be pointed out that the condition of convergence in that work corresponds strongly to the one that is proven in this present article in Section \ref{sec:order}. \rev{However, their condition is only sufficient while the one in this paper is also necessary.}

In this paper, we focus on nonstandard multistep methods. Such schemes were introduced in \cite{anguelov} by Anguelov and Lubuma. However, due to the issues mentioned in Section \ref{sec:prelim}, these ideas were abandoned and not considered by other authors. This work aims to establish that such issues can be resolved by imposing some additional conditions on the denominator function. In Section \ref{sec:prelim}, we state the most important definitions and core concepts. The main results are presented in Section \ref{sec:main} with their proofs detailed in Section \ref{sec:proof}. Two possible applications are mentioned in Section \ref{sec:Num} with several numerical simulations illustrating the effectiveness of the methods. The article concludes with some possible extensions in Section \ref{sec:conc}.

\section{Preliminaries}\label{sec:prelim}
In this paper, we \rev{analyze} an initial problem of an autonomous ordinary differential equation given in the form
\begin{equation}\label{eq:ODE}
    \begin{cases}
        y'(t) & = f(y(t)), \\
        y(t_0) & = \tilde{y},
    \end{cases}  
\end{equation}
where $y: [t_0,\;T] \rightarrow \mathbb{R}^m$ ($[t_0,\;T] \subset \mathbb{R}$, $m \in \mathbb{Z}^+$), $y(t) = \left( y_1(t),\; y_2(t), \; \dots, \; y_m(t) \right)^T$, $f: \mathbb{R}^m \rightarrow \mathbb{R}^m$ and $\tilde{y} \in \mathbb{R}^m$. We assume that \eqref{eq:ODE} has a unique \rev{global} solution\rev{:} this can be ensured if $f$ has the \rev{global} Lipschitz property\rev{, but other, weaker conditions can also be assumed}. 

To get an approximate solution of \eqref{eq:ODE}, we split our time-interval $[t_0, \;T]$ with equidistant points \linebreak $\{ t_n := t_0 + n \Delta t \}_{n = 0,\; 1, \; \dots, \; \mathcal{N}}$, where $t_{\mathcal{N}}=T$ and $\Delta t > 0$ is the step-size. For the approximations of $y(t)$ at point $t=t_n$, we use the notation $u^n \approx y(t_n)$. In line with the notations of \cite{anguelov}, the resulting numerical \rev{scheme} can be written in the form
\begin{equation}\label{eq:gennum}
    D_{\Delta t} u^n = F_{\Delta t} (f; u^n),
\end{equation}
where $D_{\Delta t} u^n$ is the approximation of the derivative of the exact solution $y(t)$ \rev{at $t_n$}, and $F_{\Delta t} (f; u^n)$ approximates the right-hand side of the differential equation. Later, the notation $F_{\Delta t} (u^n)$ is also used.

The main goal of this paper is to extend the well-known class of explicit linear multi-step methods with $k$-many steps (see e.g. \cite{hairer}) given in the form
\begin{equation}\label{eq:LMM}
   \sum_{j=0}^k \alpha_j u^{n+j} = \Delta t \sum_{j=0}^{k-1} \beta_j f(u^{n+j}). 
\end{equation}
Subsequently, we will also consider the following form used in the theory of strongly stability preserving methods (see e.g. \cite{sspbook})
\begin{equation}\label{eq:LMM2}
   u^{n+1} = \sum_{j=1}^s \left( \tilde{\alpha}_j u^{n+1-j} + \Delta t \tilde{\beta}_j f(u^{n+1-j})\right),
\end{equation}
which is beneficial when \rev{analyzing} the preservation of certain qualitative properties of the method. We can also include implicit methods if the summation starts at $j=0$.

In \cite{anguelov}, \rev{in line with the nonstandard modeling rules of Mickens, see e.g. \cite{mickens}}, Anguelov and Lubuma proposed the nonstandard version of the previous (\rev{\say{standard}}) method, which takes the form
\begin{equation}\label{eq:NSLMM}
   \sum_{j=0}^k \alpha_j u^{n+j} = \varphi(\Delta t) \sum_{j=0}^k \beta_j f(u^{n+j}),
\end{equation}
where $\varphi$ is a nonnegative real-valued one-variable function such that
\begin{equation}\label{eq:angelovcond}
    \varphi(x) = x + O(x^2) \qquad \text{as } x \rightarrow 0.
\end{equation}
\rev{It should be noted that although method \eqref{eq:NSLMM} is implicit, in this paper we only consider the explicit case.}

\rev{Moreover, another important remark is that from the two usual guiding principles of nonstandard methods, only one of them (the use of a denominator function) is taken into account, while the application of nonlocal approximations in not considered. The case when both of them are implemented could be a topic of further research.} 

In \cite{anguelov}, the consistency and the zero-stability of method \eqref{eq:NSLMM} were proven. Namely,
\begin{theorem}[\cite{anguelov}, Theorem 6]
    If \eqref{eq:LMM} represents a standard finite difference scheme that is consistent and zero-stable, then any corresponding nonstandard finite difference scheme \eqref{eq:NSLMM} (that is, with the same $\alpha_j$ and $\beta_j$ constants) is necessarily consistent. 
    Furthermore, scheme \eqref{eq:NSLMM} is zero-stable provided that the operator $F_{\Delta t}$ \rev{(as defined before)} satisfies, for some $M>0$ independent of $h$ and for any bounded sequences $(u_n)$ and $(\tilde{u}_n)$ in $\mathbb{R}^n$, the Lipschitz condition
\begin{equation}\label{eq:num_lip_cond}
    \sup_k \Vert F_{\Delta t}(u_n) - F_{\Delta t}(\tilde{u}_n) \Vert \leq M \sup_k \Vert u_n - \tilde{u}_n \Vert.
\end{equation}
\end{theorem}
Note that condition \eqref{eq:num_lip_cond} holds if scheme \eqref{eq:LMM} is zero-stable (see e.g. \cite{lambert}).

In \cite{anguelov}, Anguelov and Lubuma \rev{analyzed} the elementary stability properties of the linear multistep methods. A numerical method is called elementary stable if the numerical scheme has the same equilibrium points as the continuous model, and their stability is also the same. It was shown that under some mild conditions (with the absolute stability region of the numerical method being bounded, and the original differential equation having at least one stable equilibrium), the \rev{\say{classical}} linear multistep methods are not elementary stable when their order is greater than one. Because of this, they focused on one-step methods.

In the subsequent section, we demonstrate that under some conditions, the nonstandard linear multistep method \eqref{eq:NSLMM} attains the same order as the \rev{\say{standard}} counterpart \eqref{eq:LMM}, and it also preserves some key qualitative properties: namely, \rev{boundedness, a property similar to monotonicity}, and the linear combination of the elements of $u^n$.

\section{Main results} \label{sec:main}

In this section, we present the main results of the paper. First, the order of the nonstandard multistep methods is \rev{analyzed}, and then the preservation of qualitative properties is shown. The proofs of these statements can be found in Section \ref{sec:proof}.

\subsection{The order of nonstandard multistep methods}\label{sec:order}

Before presenting the theorem about the order of nonstandard linear multistep methods, the lemma used in the theorem's proof is highlighted. \rev{In the proof of the convergence of the nonstandard multistep methods, this nonstandard form of the usual Taylor expansion will be used.}

\begin{lemma} [Nonstandard form of Taylor's theorem] \label{lem:taylor}
    Let $n \geq 0$ be an integer and the functions $f, \varphi: \mathbb{R} \rightarrow \mathbb{R}$ be $n$-many times continuously differentiable on $[a,x]$ and $(n+1)$-many times differentiable on $(a,x)$, where $a,x \in \mathbb{R}$ ($a<x$). Moreover, let us also assume that neither $\varphi$ nor its derivatives vanish on $(a,x)$. Then, the following two claims are equivalent:
    \begin{itemize}
        \item[(C1)] If $n\geq 1$, then $\varphi(0)=0$, $\varphi'(0)=1$, and if $n\geq 2$, then $\varphi^{(k)}(0)=0$ holds for $2 \leq k \leq n$.
        \item[(C2)] There exists a function $e_n(x): \mathbb{R}\rightarrow\mathbb{R}$ such that $\displaystyle \lim_{x \rightarrow a} e_n(x) = 0$ and
        \begin{equation}\label{eq:Taylor}
            f(x) = \sum_{k=0}^n \dfrac{f^{(k)}(a)}{k!} \left(\varphi(x-a)\right)^k + e_n(x) \left(\varphi(x-a)\right)^n.
        \end{equation}
        Moreover, $e_n(x) = O(x-a)$.
    \end{itemize}
\end{lemma}

\begin{remark}
Note that if $\varphi(x)=x$, then we get the \rev{\say{usual}} form of Taylor's theorem. Moreover, one could rephrase Lemma \ref{lem:taylor} in the following way: formula \eqref{eq:Taylor} holds if and only if the first $(n+1)$-many terms of the Taylor expansions (around zero) of the functions $\varphi(x)$ and $f(x)=x$ are the same. In Section \ref{sec:phi} we show several possible choices of $\varphi$ for which condition (C1) holds for different values of $n$, and these choices are compared in Section \ref{sec:Num}.
\end{remark}

\begin{remark}
    It is worth mentioning here that one can generalize the previous lemma in the following way: let us assume that\rev{,} other than the smoothness assumption in the statement of Lemma \ref{lem:taylor}\rev{,} $\varphi(0)=0$, $\varphi(x)>0$ as $x>0$, and $\varphi'(x)\neq 0$ as $x>0$. Then, \rev{by using the intermediate value theorem,} a similar form can be proven , namely
    \begin{equation}\label{eq:Taylor2}
     \begin{aligned}
        f(x) &= \sum_{k=0}^n \dfrac{f^{(k)}(a)}{k!} \left(\varphi(x-a)\right)^k + \dfrac{f^{(n+1)}(\xi)}{(n+1)!} \dfrac{\left(\varphi(x-a)\right)^{n+1}}{\varphi'(x-\xi)} +\\
        &+  \sum_{k=0}^{n-1} \dfrac{f^{(k+1)}(\xi)}{k!} \left(\varphi(x-\xi)\right)^k \left( 1 - \varphi'(x-\xi) \right),
    \end{aligned} 
    \end{equation}
    where $\xi \in (a,x)$. Note that for $\varphi(x)=x$, this is Taylor's theorem with the Lagrange form of the remainder. Nevertheless, because of the last \rev{sum} of \eqref{eq:Taylor2}, this form is not beneficial in our case.
\end{remark}

Now, we state the main result about the order of nonstandard multistep methods.

\begin{theorem}\label{th:conv_order}
Let us assume that for equation \eqref{eq:ODE} with a smooth solution $y \in C^{p+1} ([t_0, T])$, the \rev{\say{standard}} linear multistep method \eqref{eq:LMM} is convergent of order $p$. Let us assume that for function $\varphi(x)$, it is $(p+1)$-times continuously differentiable, and neither $\varphi$ nor its derivatives vanish for $x>0$. 
Then, the nonstandard linear multistep method \eqref{eq:NSLMM} is convergent of order $p$ if and only if condition (C1) of Lemma \ref{lem:taylor} holds with $n=p$.
\end{theorem}

\rev{\begin{remark}
    Note that Theorem \ref{th:conv_order} can also be proven without the use of Lemma 1, and directly by the fact that condition (C1) is equivalent to 
    \begin{enumerate}
        \item[(C3)]  $\varphi(\Delta t) = \Delta t + O((\Delta t)^{n+1})$
    \end{enumerate}
    (see the proof of Lemma 1). However, the author believes that Lemma \ref{lem:taylor} gives a further insight into the necessity of condition $\varphi(\Delta t) = \Delta t + O((\Delta t)^{n+1})$: namely, the nonstandard Taylor expansion works if and only if (C1) (and consequently, (C3)) is true, and this nonstandard expansion ensures the convergence of the nonstandard methods (similarly as the \say{standard} Taylor expansion in the case of \say{standard} methods). In other words, by using Lemma \ref{lem:taylor}, one can get a sufficient and necessary condition for the order of the method, and the reason for the necessity of these conditions is illustrated by Lemma \ref{lem:taylor}. 
\end{remark}}

\rev{\begin{remark}
    It should be mentioned here that the absolute stability regions (see e.g. \cite{hairer2}, Ch. V.1 Def. 1.1.) of nonstandard multistep methods coincide with the same regions of the standard ones - the main reason for this is that the corresponding characteristic equation
    \begin{equation*}
        \rho(\zeta) - \mu \sigma(\zeta)=0
    \end{equation*}
    is the same as in the nonstandard case, the only difference being that the variable $\mu$ is now defined as $\mu=\varphi(\Delta t) \lambda$ instead of $\mu=\Delta t \lambda$, but because there are no restrictions on the range of $\lambda \in \mathbb{C}$, variable $\mu$ can also attain any complex value, even if $\varphi$ is bounded (see the next Section). Because of this, the second Dahlquist barrier (see \cite{dahl}, along with \cite{hairer2} and \cite{leveque}) of standard methods is also true in this, nonstandard setting.
\end{remark}}

\subsection{Preservation of qualitative properties}\label{sec:pres}
In this section, we \rev{consider} the preservation of some qualitative properties of the original differential equation \eqref{eq:ODE}. In all three cases, we consider the Forward Euler method applied to equation \eqref{eq:ODE} with timestep $\Delta t_{FE}$, and denote the resulting numerical \rev{solution} by $v^n = (v_1^n, \; v_2^n, \; \dots, \; v_m^n)^T$. 

Then, we assume that this method preserves some qualitative property under some conditions for the timestep, and then \rev{analyze} the same property for the explicit nonstandard multistep method with $s$-many steps written in the form
\begin{equation}\label{eq:NSLMM2}
   u^{n+1} = \sum_{j=1}^s \left( \tilde{\alpha}_j u^{n+1-j} + \varphi(\Delta t) \tilde{\beta}_j f(u^{n+1-j})\right),
\end{equation}
where we assume that $\tilde{\alpha}_j, \tilde{\beta}_j \geq 0$ and if $\tilde{\alpha}_j=0$, then $\tilde{\beta}_j=0$. Both the Euler method and the multistep methods are started from $v^0 = u^0 = \tilde{y} = y(0)$.

\rev{\begin{remark}
    It should be noted here that the derivation of the bound $\Delta t_{FE}$ can be challenging sometimes, and a general algorithm for its calculation might be complicated. However, for simple equations it can be determined: for example, if the solution of the equation is monotone. Two calculations can be found among the numerical experiments in Section \ref{sec:Num}.
\end{remark}}

\rev{\begin{remark} \label{rem:initialvalues}
    A topic that has been omitted until this point is the choice of the initial values of the numerical method. In \cite{dang}, it was shown that nonstandard versions of Runge-Kutta methods (i.e. ones, where the usual $\Delta t$ term is interchanged with a $\varphi(\Delta t)$ term), have (at least) the same order $p$ as their \rev{\say{standard}} counterparts and also preserve the linear stability of the system if the following conditions are met:
\begin{itemize}
    \item Let us assume that for the \rev{\say{standard}} numerical method, the linear stability and positivity are preserved if $\Delta t < \varphi^*$. Then, the nonstandard method preserves these properties if $\varphi(\Delta t)<\varphi^*$.
    \item $\varphi(\Delta t) = \Delta t + O((\Delta t)^{p+1})$.
\end{itemize}
Note that both of these conditions are met if we choose a function $\varphi$ by the ideas discussed in Section \ref{sec:phi} with $\mathcal{B} = \varphi^*$. 

Moreover, in the next Theorems we are always assuming that the one-step methods (that are used to calculate the initial values) preserve a given qualitative property. By using the strongly stability preserving form of \say{standard} Runge-Kutta methods, conditions under which these requirements are fulfilled can be obtained (see for example \cite{takacs}). Then, if we replace the usual stepsize $\Delta t$ by $\varphi(\Delta t)$ in these one-step methods, the required properties can be guaranteed for the nonstandard versions of these methods for any stepsize $\Delta t > 0$.
\end{remark}}

First, we focus on the boundedness of solutions. We say that a component $y_j(t)$ of the solution  $y(t) = (y_1(t),\; y_2(t),\;\dots, \; y_m(t))^T$ is bounded from above, if $y_j(t)\leq M$ holds for every $t \in [t_0,T]$ for some $M \in \mathbb{R}$. The boundedness from below is also similar. 

We also consider the discrete version of the same property - namely, we say that the numerical solution  $u^n = (u_1^n, \; u_2^n, \; \dots, \; u_m^n )^T$ preserves the boundedness (from above) for element $u_j^n$, if in the case of $u_j^0\leq M$ we have $u_j^n \leq M$ for every $n= 1, \; 2,\; \dots, \; \mathcal{N}$. Again, a similar definition can also be made for boundedness from below.

Let us assume that for component $y_k(t)$ of the solution of equation \eqref{eq:ODE} we have $y_k(t)\leq M$ for $t \in [t_0, T]$. The next theorem states that in this case (under some conditions) the numerical method \eqref{eq:NSLMM2} also preserves this property for any choice of the step size.

\begin{theorem}\label{th:bounded}
     Let us assume that the following conditions are met:
     \begin{itemize}
         \item If $\Delta t_{FE}\leq B_{FE}$ for some constant $B_{FE}\in \mathbb{R}^+$, then the Forward Euler method applied to \eqref{eq:ODE} preserves the boundedness property, namely $v_k^1 \leq M$ holds.
         \item For the starting values of the multistep method \eqref{eq:NSLMM2}, we have $u_k^0, u_k^1, \dots, u_k^{s-1}\leq M$.
         \item $\varphi(x) \leq \mathcal{C} B_{FE} $, where $\mathcal{C}\in \mathbb{R}^+$ is the SSP coefficient defined as (see \cite{sspbook})
    $$ \mathcal{C} = \min_j \dfrac{\tilde{\alpha}_j}{\tilde{\beta}_j}. $$
     \end{itemize}
Then, for the nonstandard multistep method \eqref{eq:NSLMM2}, $u_k^n \leq M$ also holds for $n=s,\dots \mathcal{N}$ for any stepsize $\Delta t\geq 0$.
\end{theorem}

The same theorem also holds if we consider boundedness from below.

Another theorem can be stated if we consider the monotone property of the solution. Let us assume that the component $y_k(t)$ of the solution of equation \eqref{eq:ODE} is monotone increasing for $t \in [t_0, T]$ if $\tilde{y} \in (A,B)$ element-wise. \rev{Here} $A,B \in \mathbb{R}, \; A<B$, where we also allow that $A=-\infty$ or $B=\infty$, \rev{meaning that we also allow cases like $\tilde{y} \in [0, \infty)$ and not only $\tilde{y} \in [0, C]$ ($C >0$)}.

\begin{theorem}\label{th:monotone}
    Let us assume that the following conditions are met:
    \begin{itemize}
        \item If $\Delta t_{FE}\leq B_{FE}$ for some constant $B_{FE}\in \mathbb{R}^+$, then the Forward Euler method applied to \eqref{eq:ODE} preserves the monotone increase property, namely $v_k^1 \geq v_k^0$, and \linebreak $v_k^1 \in (A,B)$ holds if $v_k^0 \in (A,B)$. 
        \item For the starting values of the multistep method \eqref{eq:NSLMM2}, we have 
        $$u_k^0, \; u_k^1, \; \dots, \; u_k^{s-1} \in (A,B).$$
        \item $\varphi(x) \leq \mathcal{C} B_{FE} $, where $\mathcal{C}\in \mathbb{R}^+$ is the SSP coefficient defined as in Theorem \ref{th:bounded}. 
    \end{itemize}
    Then, for the nonstandard multistep method \eqref{eq:NSLMM2}, 
    \begin{equation}\label{eq:weakmon}
        \min_{\ell=n-s+1, \; n-s+2, \; \dots, \; n}\{u_k^{\ell}\} \leq u_k^{n+1}
    \end{equation}
    also holds for $n=s,\dots \mathcal{N}$ for any stepsize $\Delta t\geq 0$.
\end{theorem}

An analogous theorem can be proven if we consider monotone decrease: namely,  the final inequality in this case takes the form
$$\max_{\ell=n-s+1, \; n-s+2, \; \dots, \; n}\{u_k^{\ell}\} \geq u_k^{n+1}.$$

\begin{remark}
    Note that property \eqref{eq:weakmon} is only a necessary, but not sufficient condition for the monotone increase of numerical solution $u_k^{n}$. Numerical experiments indicate that the numerical scheme can indeed preserve the \rev{\say{classical}} monotone property \rev{(i.e. $u_k^n \leq u_k^{n+1}$)} if the bound for $\varphi(x)$ is even smaller than $\mathcal{C} B_{FE}$, but if we only assume that $\varphi(x)\leq \mathcal{C} B_{FE}$, then only inequality \eqref{eq:weakmon} holds (see the tests in Section \ref{sec:Num}). The determination of a bound under which the \rev{\say{classical monotone property}} is also preserved could be a topic of further research.
\end{remark}

The next theorem is about a property concerning the preservation of some linear combination of the elements of $u^n$. Let us assume that if $y(t_0) \in S$ for some $S \subset \mathbb{R}^m$ hyperrectangle, then for some linear function $P:\mathbb{R}^m \rightarrow \mathbb{R}$ we have 
    $$ P(y_1(t), \; y_2(t),\; \dots , y_m(t)) = \mathcal{M} \qquad \forall t \in [t_0, \; T]$$
    for some $\mathcal{M}\in \mathbb{R}$. 
\begin{theorem}\label{th:conserv}
 Let us assume that the following conditions are met:
 \begin{itemize}
     \item Consider the Forward Euler method starting from $v^0 \in S$. If $\Delta t_{FE}\leq B_{FE}$ for some constant $B_{FE}\in \mathbb{R}^+$, then $P(v_1^n, \; v_2^n,\; \dots , v_m^n) = \mathcal{M}$ for $n=1, 2, \dots$ assuming that $P(v_1^0, \; v_2^0,\; \dots , v_m^0) = \mathcal{M}$. 
     \item For the starting values of the multistep method \eqref{eq:NSLMM2}, we have 
     $$P(u_1^n, \; u_2^n,\; \dots , u_m^n) = \mathcal{M}$$
      and $u^n \in S$ holds for $n=0,\; 1,\; \dots, \; s-1$.
     \item $\varphi(x) \leq \mathcal{C} B_{FE} $, where $\mathcal{C}\in \mathbb{R}^+$ is the SSP coefficient defined as in Theorem \ref{th:bounded}. 
 \end{itemize}
  Then,
    $$P(v_1^n, \; v_2^n,\; \dots , v_m^n) = \mathcal{M}$$ 
    also holds for $n=s,\dots \mathcal{N}$.  
\end{theorem}

\begin{remark}
    If $P(y_1(t), \; y_2(t), \; \dots, \; y_m(t) ) = \sum_{j=1}^m y_j(t)$, then the previous theorem is about the preservation of the sum of all components (see Section \ref{sec:seir} for such an example).
\end{remark}

\subsection{\texorpdfstring{The choice of function $\varphi$}{The choice of function phi}}\label{sec:phi}
Here, we briefly address the issue of the choice of function $\varphi$. As it was mentioned in the previous sections, for a nonstandard multistep method to be of order $p$ and to preserve some important qualitative properties, the following two conditions should be met:
\begin{itemize}
    \item[(A)] $\varphi$ should be $(p+1)$-times continuously differentiable (for $x>0$) and $\varphi(x)>0$ for $x>0$. Moreover, if $p\geq 1$, then $\varphi(0)=0$, $\varphi'(0)=1$ and if $p\geq 2$ then $\varphi^{(k)}(0)=0$ holds for $2 \leq k \leq p$.
    \item[(B)] $\varphi(x)\leq \mathcal{B} := \mathcal{C} B_{FE}$ for every $x>0$. 
\end{itemize}
In the following, we outline different possible choices for $\varphi$. These are chosen with condition (B)  fulfilled, and they are arranged according to the value $p$ for which condition (A) is satisfied.  

\begin{itemize}
    \item With $p=1$:
    \begin{itemize}
        \item $\varphi_1(x)=\mathcal{B}(1-e^{-\frac{x}{\mathcal{B}}})$ (here $\varphi''(0)=-\frac{1}{\mathcal{B}} \neq 0$),
        \item $\varphi_2(x)=xe^{-\frac{x}{\mathcal{B} e}}$ (here $\varphi''(0)=-\frac{2}{\mathcal{B} e} \neq 0$),
        \item $\varphi_3(x) = \dfrac{\mathcal{B} x}{\mathcal{B}+x}$ (here $\varphi''(0) = -\frac{2}{\mathcal{B}}\neq 0$).
    \end{itemize}
    \item With $p=2$:
    \begin{itemize}
        \item $\varphi_4(x)=\dfrac{\mathcal{B} 2}{\pi}\mathrm{arctan}\left(\dfrac{x \pi}{2 \mathcal{B}}\right)$ (here $\varphi'''(0)=-\dfrac{\pi^2}{2 \mathcal{B}^2}\neq 0$),
        \item $\varphi_5(x)=\mathcal{B}\mathrm{tanh}\left(\dfrac{x}{\mathcal{B}}\right)$ (here $\varphi'''(0)=-\dfrac{2}{\mathcal{B}^2}\neq 0$).
        \item $\varphi_6(x) = \dfrac{\mathcal{B} x}{\left( \mathcal{B}^2 + x^2 \right)^{1/2}}$ (here $\varphi'''(0) = - \frac{3 }{\mathcal{B}^2}\neq 0$).
    \end{itemize}
    \item With $p=3$: $\varphi_7(x) = \dfrac{\mathcal{B} x}{\left( \mathcal{B}^3 + x^3\right)^{1/3}}$ (here $\varphi^{(4)} (0) =-\dfrac{8}{\mathcal{B}^3}\neq 0$).
    \item With $p=4$: $\varphi_8(x) = \dfrac{\mathcal{B} x}{\left( \mathcal{B}^4 + x^4\right)^{1/4}}$ (here $\varphi^{(5)} (0) =-\dfrac{30}{\mathcal{B}^4}\neq 0$).
    \item For the general $p \geq 5$ case: $\varphi(x) = \dfrac{\mathcal{B} x}{\left( \mathcal{B}^p + x^p\right)^{1/p}}$ (here $\varphi^{(p+1)} (0) =-\dfrac{(p-1)! + p!}{\mathcal{B}^p}\neq 0$).
\end{itemize}
We should also note here that in \cite{dang}, Dang and Hoang demonstrated a different construction of functions fulfilling these conditions up to an arbitrary order.

\rev{\begin{remark}
    It should be noted that although the condition $\varphi(x)\leq \mathcal{B}$ ensures the preservation of the aforementioned qualitative properties, such choices where the function $\varphi(x)$ is close to the constant function $\varphi_0(x)\equiv \mathcal{B}$ might still result in considerably large errors when using large timesteps, similarly as in the case of \say{standard} methods. The benefit of these nonstandard methods is that even for larger timesteps, the numerical solution behaves in a qualitatively reasonable way, although the errors might still be relatively big. For smaller timesteps, the methods behave as expected resulting in small errors, like other, high order standard methods. 
\end{remark}}

In Section \ref{sec:Num} we perform several numerical experiments comparing these choices. \rev{There, it turns out that it is usually beneficial to apply such $\varphi$ functions which allow a higher order of convergence (e.g. $\varphi_8$), since their application generally results in smaller errors, and sometimes even higher rates of convergence than the one which is guaranteed by the multistep method itself, see the third order methods in Section \ref{sec:seir} that attain an order of four. Numerical tests (see Table \ref{tab:phiruntimes}) indicate that there is no significant difference between the runtime of methods using e.g. $\varphi_3$ or $\varphi_8$, so it is advised to use the latter: the functions that require more time to compute are $\varphi_1$ and $\varphi_2$, but their application is not recommended since they decrease the order of the methods significantly.}

\begin{table}[!h]
    \centering
    \begin{tabular}{c|c|c|c|c|c|c|c|c}
         & $\varphi_1$ & $\varphi_2$ & $\varphi_3$ & $\varphi_4$ & $\varphi_5$ & $\varphi_6$ & $\varphi_7$ & $\varphi_8$  \\ \hline
        runtime (s) & $7.2136$ & $7.0121$ & $0.4383$ & $0.4377$ & $0.4358$ & $0.4447$ & $0.4382$ & $0.4410$
    \end{tabular}
    \caption{The time (seconds) it takes to calculate the value of a given function $10^9$-many times with $\mathcal{B}=0.1$ and $x=0.1$. The tests were made by using Matlab R2024b.}
    \label{tab:phiruntimes}
\end{table}

\section{Proofs of the statements of Section \ref{sec:main}} \label{sec:proof}
In this section, we present the proofs of the statements of the previous section.

\begin{proof}[Proof of Lemma \ref{lem:taylor}]
\rev{The main idea of the proof is that condition (C1) is equivalent to the following condition: \small
\begin{enumerate}
    \item[(C3)]  \small the (standard) Taylor series of function $\varphi$ around $0$ at point $x-a$ is in the form
    \begin{equation*}
    \varphi(x-a)=(x-a)+O((x-a)^{n+1}).
\end{equation*}
\end{enumerate} 
It is easy to see that if (C1) holds, then by the standard Taylor expansion of $\varphi$ around $0$ we get (C3). Moreover, (C3) can only hold if (C1) is true (again, by the form of the standard Taylor expansion of $\varphi$).

We prove the two directions of the lemma separately.
\begin{enumerate}
    \item  \small Let us start with the direction (C1)$\Rightarrow$(C2). By the initial remark, it is enough to prove (C3)$\Rightarrow$(C2). Let us consider the standard Taylor expansion of $f$ around point $a$:
    \begin{equation}\label{eq:fTaylor}
        f(x) = \sum_{k=0}^n \dfrac{f^{(k)}(a)}{k!} (x-a)^k + O((x-a)^{n+1}).
    \end{equation}
    Next, we substitute the statement of (C3) into the right-hand side of (C2):
    $$ \sum_{k=0}^n \dfrac{f^{(k)}(a)}{k!} \left(\varphi(x-a)\right)^k + O(x-a) \left(\varphi(x-a)\right)^{n} = $$
     $$ = \sum_{k=0}^n \dfrac{f^{(k)}(a)}{k!} \left((x-a)+O((x-a)^{n+1})\right)^k + O(x-a) \left((x-a)+O((x-a)^{n+1})\right)^{n} = $$
     $$ = \sum_{k=0}^n \dfrac{f^{(k)}(a)}{k!} \left(x-a)\right)^k + O\left((x-a)^{n+1}\right) = f(x), $$
     where we used \eqref{eq:fTaylor}.

     \item Now consider the direction (C2)$\Rightarrow$(C1). Let us substitute the standard Taylor expansion $\varphi$ around $0$ into (C2):
     $$
     \begin{aligned}
         f(x) = \sum_{k=0}^n \dfrac{f^{(k)}(a)}{k!} \left( \varphi(0) + \varphi'(0) (x-a) + \sum_{\ell=2}^n \dfrac{\varphi^{(\ell)}(0)}{\ell!} (x-a)^{\ell} + O((x-a)^{n+1}) \right)^k + \\
         + O(x-a) \left( \varphi(0) + \varphi'(0) (x-a) + \sum_{\ell=2}^n \dfrac{\varphi^{(\ell)}(0)}{\ell!} (x-a)^{\ell} + O((x-a)^{n+1}) \right)^n.
     \end{aligned}   $$
     By comparing this expression with \eqref{eq:fTaylor}, they can only hold if
     $$ \left( \varphi(0) + \varphi'(0) (x-a) + \sum_{\ell=2}^n \dfrac{\varphi^{(\ell)}(0)}{\ell!} (x-a)^{\ell}  \right)^k = (x-a)^k, \qquad k=1, \; 2, \; \dots \; n, $$
     which can only hold if (C1) is true.\qedhere
\end{enumerate}
}
\end{proof}


\begin{proof}[Proof of Theorem \ref{th:conv_order}]
We prove the two directions separately.
\begin{enumerate}
    \item  \small First we show that from condition (C1), we get a convergence of order $p$.

    Let us \rev{consider} the following linear difference operator $L$ (see e.g. \cite{hairer}):
    $$ L(y,t,\Delta t) = \sum_{j=0}^k \left( \alpha_j y(t+j \Delta t) - \varphi(\Delta t) \beta_j y'(t+j \Delta t) \right), $$
    which we got by substituting the exact solution $y(t)$ into scheme \eqref{eq:NSLMM}. Now, let us apply the nonstandard Taylor expansion (around $a=0$) defined in Lemma \ref{lem:taylor} for functions $y$ and $y'$ at points $x=t+j \Delta t$:
    $$ L(y,t,\Delta t) = \sum_{j=0}^k \alpha_j \sum_{q=0}^p \dfrac{\left(\varphi(j \Delta t)\right)^{q}}{q!} y^{(q)}(t)  - \varphi(\Delta t) \sum_{j=0}^k \beta_j \sum_{r=0}^p \dfrac{\left(\varphi(j \Delta t)\right)^{r}}{r!} y^{(r+1)}(t)  + O((\Delta t)^{p+1}) = $$
    $$ = y(t) \sum_{j=0}^k \alpha_j +  \sum_{q=1}^p \dfrac{y^{(q)}(t)}{q!} \sum_{j=0}^k \alpha_j \left(\varphi(j \Delta t)\right)^{q} - \varphi(\Delta t)  \sum_{r=0}^p \dfrac{y^{(r+1)}(t)}{r!} \sum_{j=0}^k \beta_j \left(\varphi(j \Delta t)\right)^{r} + O((\Delta t)^{p+1})= $$
    $$ = \sum_{q=1}^p \dfrac{y^{(q)}(t)}{q!} \sum_{j=0}^k \alpha_j \left(\varphi(j \Delta t)\right)^{q} - \varphi(\Delta t)  \sum_{r=1}^{p+1} r \dfrac{y^{(r)}(t)}{r!} \sum_{j=0}^k \beta_j \left(\varphi(j \Delta t)\right)^{r-1} + O((\Delta t)^{p+1}), $$
    where we used that the scheme \eqref{eq:LMM} is consistent, meaning that $\sum_{j=0}^k \alpha_j=0$ holds. Moreover, if condition (C1) is true for $n=p$, then it means that $\varphi(h)=h+O(h^{p+1})$, which, substituted into the previous formula, gives
    $$ \sum_{q=1}^p \dfrac{y^{(q)}(t)}{q!} \left(\Delta t\right)^{q} \sum_{j=0}^k \alpha_j j^q  -  \sum_{r=1}^{p+1} r \dfrac{y^{(r)}(t)}{r!} \left(\Delta t\right)^{r} \sum_{j=0}^k \beta_j j^{r-1}  + O((\Delta t)^{p+1}) = $$
    $$ = \sum_{q=1}^p \dfrac{y^{(q)}(t)}{q!} \left(\Delta t\right)^{q} \left(\sum_{j=0}^k \alpha_j j^q  -  q\sum_{j=0}^k \beta_j j^{q-1}\right)  + O((\Delta t)^{p+1}). $$
    If scheme \eqref{eq:LMM} has an order of $p$, then $\displaystyle \sum_{j=0}^k \alpha_j j^q  -  q\sum_{j=0}^k \beta_j j^{q-1}=0$ for $q=1 ,\dots p$ (see e.g. \cite{hairer}, Theorem 2.4.), meaning that the nonstandard method also has an order of $p$.

\item Now we show that an order of $p$ can only hold if (C1) is satisfied with $n=p$.

We prove the statement by induction. Let us consider the case $n=1$, i.e., the method's order is at least one. It is obvious that consistency requires \linebreak $\varphi(0)=0$. Moreover, if we consider the Taylor expansion of $\varphi$ around zero, then (since  $\varphi(0)=0$) we have $\varphi(\Delta t) = \varphi'(0) \Delta t + O((\Delta t)^2)$. By substituting this into operator $L$ (defined in the first part of the proof), we get
$$ L(y,t,\Delta t) = \sum_{j=0}^k \left( \alpha_j y(t+j \Delta t) - \varphi'(0) \Delta t \beta_j y'(t+j \Delta t) \right) + O((\Delta t)^2).$$
Let us apply the \rev{\say{standard}} Taylor expansion for the functions $y(t+j\Delta t)$ and $y'(t+j\Delta t)$ at points $x=t+j \Delta t$:
$$ L(y,t,\Delta t) = \sum_{j=0}^k \alpha_j \left( y(t) + j \Delta t y'(t) \right)  - \varphi'(0) \Delta t \sum_{j=0}^k \beta_j \left( y'(t) + j\Delta t y''(t) \right)  + O((\Delta t)^{2}) = $$
$$ = \sum_{j=0}^k \alpha_j y(t) +  y'(t) \Delta t \sum_{j=0}^k\left(j \alpha_j - \varphi'(0) \beta_j  \right) + O((\Delta t)^2)= $$ 
$$ = y'(t) \Delta t \sum_{j=0}^k\left(j \alpha_j - \varphi'(0) \beta_j  \right) + O((\Delta t)^2), $$
where we used that $\sum_{j=0}^k \alpha_j=0$ holds. Moreover, we also know that \linebreak $\displaystyle \sum_{j=0}^k \alpha_j j  -  \sum_{j=0}^k \beta_j =0$, meaning that the above expression equals $O((\Delta t)^2)$ only if $\varphi'(0)=1$.

Now, let us assume that condition (C1) holds for $n<N$ and we prove the $n=N$ case ($2 \leq N \leq p$). The argument is rather similar to the one in the $n=1$ step: in this case $\varphi(\Delta t) = \Delta t + \dfrac{\varphi^{(N)}(0)}{N!} (\Delta t)^N + O((\Delta t)^{N+1})$. By substitution, we get
$$ L(y,t,\Delta t) = \sum_{j=0}^k \left( \alpha_j y(t+j \Delta t) - \left( \Delta t + \dfrac{\varphi^{(N)}(0)}{N!} (\Delta t)^N \right) \beta_j y'(t+j \Delta t) \right) + O((\Delta t)^{N+1}). $$
As before, we apply the \rev{\say{standard}} Taylor expansion to functions $y(t+j\Delta t)$ and $y'(t+j\Delta t)$ at points $x=t+j\Delta t$, resulting in
$$ \begin{aligned}
    L(y,t,\Delta t) &= \sum_{j=0}^k \left( \alpha_j \sum_{q=0}^N \dfrac{\left(j \Delta t\right)^{q}}{q!} y^{(q)}(t) - \right. \\
&-\left.\left( \Delta t + \dfrac{\varphi^{(N)}(0)}{N!} (\Delta t)^N \right) \beta_j \sum_{r=0}^N \dfrac{\left(j \Delta t\right)^{r}}{r!} y^{(r+1)}(t) \right) + O((\Delta t)^{N+1}) =
\end{aligned}$$
$$ = \sum_{q=1}^N \dfrac{y^{(q)}(t)}{q!} \left(\Delta t\right)^{q} \left(\sum_{j=0}^k \alpha_j j^q  -  q\sum_{j=0}^k \beta_j j^{q-1}\right) - \dfrac{\varphi^{(N)}(0)}{N!} (\Delta t)^N \sum_{j=0}^k \beta_j y'(t) + O((\Delta t)^{N+1}) = $$
\begin{equation}\label{eq:orderproof}
    = - \dfrac{\varphi^{(N)}(0)}{N!} (\Delta t)^N \sum_{j=0}^k \beta_j y'(t) + O((\Delta t)^{N+1}),
\end{equation}
where we used that $\displaystyle \sum_{j=0}^k \alpha_j j^q  -  q\sum_{j=0}^k \beta_j j^{q-1}=0$ for $q=1 ,\dots N$ ($N\leq p$) and $\sum_{j=0}^k \alpha_j=0$. Therefore, expression \eqref{eq:orderproof} can only equal $O((\Delta t)^{N+1})$ in general if $\varphi^{(N)}(0)=0$. This concludes this part of the proof. \qedhere
\end{enumerate}
    
\end{proof}

\begin{proof}[Proof of Theorem \ref{th:bounded}]
We prove the statement by induction: let us assume that $u_k^0,\; u_k^1,\; \dots u_k^N\leq M$ and we aim to show that $u_k^{N+1}\leq M$ (here $N+1\leq \mathcal{N}$).     
    
Let us consider the numerical method \eqref{eq:NSLMM2} and \rev{analyze} its part concerning the element $u_k$: 
\begin{equation}\label{eq:boundedness_proof}
    u_k^{N+1} = \sum_{j=1}^s \tilde{\alpha}_j \left(  u_k^{N+1-j} + \varphi(\Delta t) \dfrac{\tilde{\beta}_j}{\tilde{\alpha}_j} f_k(u^{N+1-j})\right),
\end{equation}
where we used the notation $f_k$ for the $k$th element of the function $f:\mathbb{R}^m \rightarrow \mathbb{R}^m$.
Since the boundedness property holds for the forward Euler method if $\Delta t_{FE} \leq B_{FE}$, for $\Delta t_{FE} \leq B_{FE}$ we have
$$ \left( v_k^{0} + \Delta t  f_k(v^{0})\right) \leq M, $$
where $v_k^{0}$ is an arbitrary value for which $v_k^0 \leq M$ holds.
By the assumption of the theorem, $\varphi(\Delta t)\leq \mathcal{C} B_{FE}$ is true, meaning that for every $j$
$$ \varphi(\Delta t) \dfrac{\tilde{\beta}_j}{\tilde{\alpha}_j} \leq \mathcal{C} B_{FE} \dfrac{\tilde{\beta}_j}{\tilde{\alpha}_j} \leq B_{FE}.$$
Therefore, every term of the sum \eqref{eq:boundedness_proof} can be thought of as a Forward Euler step with stepsize $\varphi(\Delta t) \dfrac{\tilde{\beta}_j}{\tilde{\alpha}_j}\leq B_{FE}$, resulting in
$$ \sum_{j=1}^s \tilde{\alpha}_j \left(  u_k^{N+1-j} + \varphi(\Delta t) \dfrac{\tilde{\beta}_j}{\tilde{\alpha}_j} f_k(u^{N+1-j})\right) \leq \sum_{j=1}^s \tilde{\alpha}_j M = M, $$
where we used that $\sum_{j=1}^s \tilde{\alpha}_j=1$ (this holds because the method is consistent, see \cite{sspbook}).
\end{proof}

\begin{proof}[Proof of Theorem \ref{th:monotone}]
We prove the statement by induction: let us assume that 
$$\displaystyle u_k^{n+1} \geq \min_{\ell=n-s+1, \; n-s+2, \; \dots, \; n}\{u_k^{\ell}\} \qquad \text{for} \; n=s-1, \; s, \; \dots, \; N-1$$ 
and our goal is to show that $\displaystyle u_k^{N+1} \geq \min_{\ell=N-s+1, \; N-s+2, \; \dots, \; N}\{u_k^{\ell}\}$ ($N+1\leq \mathcal{N}$).     
    
We again consider the numerical method \eqref{eq:NSLMM2} and \rev{analyze} its part concerning the element $u_k$ (see \eqref{eq:boundedness_proof} in the proof of Theorem \ref{th:bounded}). 
Since the monotone increase holds for the forward Euler method if $\Delta t_{FE} \leq B_{FE}$, for $\Delta t_{FE} \leq B_{FE}$ we have
$$  v_k^{0} + \Delta t  f_k(v^{0}) \geq v_k^0, $$
where $v_k^{0}$ is an arbitrary value for which $v_k^0 \in (A,B)$ holds.
By the assumption of the theorem, $\varphi(\Delta t)\leq \mathcal{C} B_{FE}$ is true, meaning that for every $j$,
$ \varphi(\Delta t) \dfrac{\tilde{\beta}_j}{\tilde{\alpha}_j} \leq B_{FE}$. Moreover, by the statement of Theorem \ref{th:bounded}, $u_k^n \in (A,B)$ for $n=s-1, \; s, \; \dots, \; N-1$.
Therefore, every term of the sum \eqref{eq:boundedness_proof} can be thought of as a Forward Euler step with stepsize $\varphi(\Delta t) \dfrac{\tilde{\beta}_j}{\tilde{\alpha}_j}\leq B_{FE}$, resulting in
$$ u_k^{N+1} = \sum_{j=1}^s \tilde{\alpha}_j \left(  u_k^{N+1-j} + \varphi(\Delta t) \dfrac{\tilde{\beta}_j}{\tilde{\alpha}_j} f_k(u^{N+1-j})\right) \geq \sum_{j=1}^s \tilde{\alpha}_j u_k^{N+1-j} \geq $$
$$\geq \sum_{j=1}^s \tilde{\alpha}_j \min_{\ell=N+1-s, \; N+2-s, \; \dots, \; N}\{u_k^{\ell}\} \geq \min_{\ell=N-s+1, \; N-s+2, \; \dots, \; N}\{u_k^{\ell}\},$$
where we used that $\sum_{j=1}^s \tilde{\alpha}_j=1$.
\end{proof}

\begin{proof}[Proof of Theorem \ref{th:conserv}]
We prove the statement by induction: let us assume that $P(u_1^n, \; u_2^n,\; \dots , u_m^n) = \mathcal{M}$ for \linebreak $n=s, \; s+1, \; \dots, \; N$ and our goal is to show that $P(u_1^{N+1}, \; u_2^{N+1},\; \dots , u_m^{N+1}) = \mathcal{M}$ (\rev{here} $N+1\leq \mathcal{N}$).     
    
Let us consider the numerical method \eqref{eq:NSLMM2} and observe its part concerning the element $u_k$ (see \eqref{eq:boundedness_proof} in the proof of Theorem \ref{th:bounded}).
Since the preservation of the value $\mathcal{M}$ holds for the forward Euler method if $\Delta t_{FE} \leq B_{FE}$, for $\Delta t_{FE} \leq B_{FE}$ we have
$$  P(v_1^{N+1}, \; v_2^{N+1},\; \dots , v_m^{N+1}) = \sum_{\ell =1}^{m} \gamma_{\ell} v_{\ell}^{N+1} = \sum_{\ell =1}^{m} \gamma_{\ell} \left( v_{\ell}^{N} + \Delta t f_{\ell} (v^{N}) \right) = \sum_{\ell =1}^{m} \gamma_{\ell} v_{\ell}^{N} = M, $$
where $\gamma_{\ell} \in \mathbb{R}$ are some given constants.

By the assumption of the theorem, $\varphi(\Delta t)\leq \mathcal{C} B_{FE}$ is true, meaning that for every $j$, $\varphi(\Delta t) \dfrac{\tilde{\beta}_j}{\tilde{\alpha}_j} \leq B_{FE}$. 
Therefore, every term of the sum \eqref{eq:boundedness_proof} can be thought of as a Forward Euler step with stepsize $\varphi(\Delta t) \dfrac{\tilde{\beta}_j}{\tilde{\alpha}_j}\leq B_{FE}$, resulting in
$$ P(u_1^{N+1}, \; u_2^{N+1},\; \dots , u_m^{N+1}) = \sum_{\ell =1}^{m} \gamma_{\ell} u_{\ell}^{N+1} = \sum_{\ell =1}^{m} \gamma_{\ell} \sum_{j=1}^s \tilde{\alpha}_j \left(  u_{\ell}^{N+1-j} + \varphi(\Delta t) \dfrac{\tilde{\beta}_j}{\tilde{\alpha}_j} f_{\ell}(u^{N+1-j})\right) =  $$
$$ =  \sum_{j=1}^s \tilde{\alpha}_j \sum_{\ell =1}^{m} \gamma_{\ell} \left(  u_{\ell}^{N+1-j} + \varphi(\Delta t) \dfrac{\tilde{\beta}_j}{\tilde{\alpha}_j} f_{\ell}(u^{N+1-j})\right) = \sum_{j=1}^s \tilde{\alpha}_j M = M, $$
where we used that $\sum_{j=1}^s \tilde{\alpha}_j=1$.
\end{proof}

\section{Applications and numerical simulations}\label{sec:Num}

In this section, we apply the previously \rev{analyzed} nonstandard multistep method \eqref{eq:NSLMM2} to different initial value problems. 

Prior to showing some specific applications, the main steps of the construction of property-preserving nonstandard multistep methods are outlined. 
\begin{enumerate}
    \item Assume that the original differential equation \eqref{eq:ODE} attains one of the qualitative properties mentioned in Section \ref{sec:pres}, which should be preserved by the numerical method too.
    \item Apply the Forward Euler method to the given equation with $\Delta t$ and calculate such a bound $B_{FE}$, for which the \rev{considered} property is preserved if $\Delta t \leq B_{FE}$ holds.
    \item Choose function $\varphi(x)$ in a way that $\varphi(x)\leq \mathcal{C} B_{FE}$. For this, one of the possible choices listed in Section \ref{sec:phi} can be used with $\mathcal{B} = \mathcal{C} B_{FE}$.
    \item For a method with $s$-many steps, the values $u^1, \; u^2, \; \dots, \; u^{s-1}$ should be calculated. For this, a nonstandard Runge-Kutta method can be used with an appropriate function $\varphi(x)$ \rev{(see Remark \ref{rem:initialvalues})}. 
    \item Thus, the nonstandard multistep method in the form \eqref{eq:NSLMM2} preserves the desired property and has the same order as the standard counterpart.
\end{enumerate}
In Section \ref{sec:log}, the method is applied to a one-dimensional problem, while in Section \ref{sec:seir} a system is considered.

\subsection{One-dimensional logistic growth equation}\label{sec:log}
Consider the one-dimensional ($n=1$) logistic growth equation \rev{with constant $c \in \mathbb{R}^+$}
\begin{equation}\label{eq:logistic}
   \begin{cases}
       y'(t) & = y(t) \left( \rev{c} - y(t) \right),\\
       y(0) & = \tilde{y}.
   \end{cases}
\end{equation}
The solution of the above equation is $y(t) = \rev{\dfrac{c e^{c t}\tilde{y}}{\tilde{y} (e^{ct}-1)+c}}$. It is also apparent that the two equilibria of the equation are $y_1^*=0$ and $y_2^*=c$, where $y_1^*$ is unstable and $y_2^*$ is asymptotically stable. The solution of the equation is also monotone.

To apply the nonstandard multistep method, we should first \rev{analyze} the behavior of the Forward Euler method.

\begin{prop}\label{prop:FE_logistic}
    The Forward Euler method applied to equation \eqref{eq:logistic} preserves the monotonicity and the boundedness of the solution if $\Delta t \leq \min \left\{ \rev{\dfrac{1}{c}}, \dfrac{1}{\tilde{y}} \right\}$ \rev{and} $\tilde{y}\geq 0$. If $\tilde{y}<0$, then the method preserves the properties unconditionally.
\end{prop}

\begin{proof}
    Let us apply the forward Euler method to equation \eqref{eq:logistic}:
    $$ \rev{v^{n+1} = v^n + \Delta t v^n (c-v^n)= v^n (1 + c \Delta t - \Delta t \; v^n).} $$

    We prove the statement by induction: let us assume that the method behaves as expected for $n=0, \dots, N$. Consider four different cases:
    \begin{itemize}
        \item If either $\tilde{y}=0$ or $\tilde{y}=\rev{c}$, the numerical solution remains constant.
        \item Let us assume that $\tilde{y}=\rev{v^0}<0$. In this case, the solution should remain negative, and it should be monotone decreasing. If $\rev{v^n}<0$, then $\rev{v^{n+1}}$ is also negative if $1+\rev{c} \Delta t - \Delta t \; \rev{v^n}>0$, but this is true. For the decreasing property, we need that $1+\rev{c} \Delta t - \Delta t \; \rev{v^n} >1$, but this holds too. Consequently, in this case, the method preserves the boundedness and the monotonicity of the solution unconditionally.
        \item Assume that $0 < \tilde{y} < \rev{c}$. For the increasing property, it is enough to show that \linebreak $1 + \rev{c} \Delta t - \Delta t \; \rev{v^n} >1$, but this is true since $\rev{v^n < c}$. For the boundedness from above, consider the function $g(x)=x (1 + \rev{c} \Delta t - \Delta t \; x)$. This attains its maximum at $ \dfrac{1+\rev{c} \Delta t}{2 \Delta t}$, and its value at this point is $\dfrac{(1+\rev{c} \Delta t)^2}{4 \Delta t}$ which is bigger than $\rev{c}$ for every positive value of $\Delta t$. Because of this, let us choose $\Delta t$ in a way that the maximum point is bigger than $\rev{c}$: $\Delta t < \frac{1}{\rev{c}}$ is sufficient. Therefore, the function $g$ is monotone increasing on $[0,\rev{c}]$, with its maximum point at $x=\rev{c}$ with value $g(\rev{c})=\rev{c}$. Therefore, $\rev{v^{n+1} = g(v^n)}<\rev{c}$ if $\Delta t < \dfrac{1}{\rev{c}}$. The boundedness from below is trivial because of the monotone increasing property. 
        \item Consider the case $\tilde{y}>\rev{c}$. For the decreasing property we need $1 + \rev{c} \Delta t - \Delta t \; \rev{v^n}<1$ but this is true since $\rev{v^n}>\rev{c}$. For the boundedness from below, consider the function $g(x)=x (1 + \rev{c} \Delta t - \Delta t \; x)$. We should prove that if $\tilde{y}>x>\rev{c}$, then $g(x)>\rev{c}$. In the previous case we have shown that if $\Delta t < \frac{1}{\rev{c}}$, then the maximum point is bigger than $\rev{c}$, so since $g(\rev{c})=\rev{c}$, for a sufficiently small $\varepsilon$, $g(x)>\rev{c}$ if $x \in (\rev{c}, \rev{c}+\varepsilon)$. Our goal is to choose $\Delta t$ in a way that $\rev{c}+\varepsilon>\tilde{y}$.

        If we solve $g(x)=\rev{c}$, the solution (apart from the previously mentioned $x_0=\rev{c}$) is $x_1=\frac{1}{\Delta t}$. Consequently, if $\Delta t \leq \frac{1}{\tilde{y}}$, then the scheme preserves the lower bound of the method. \qedhere
    \end{itemize}
\end{proof}

In conclusion, a proper choice for parameter $\mathcal{B}$ is $\mathcal{B}:=\mathcal{C} \min \left\{ \dfrac{1}{\rev{c}}, \dfrac{1}{\tilde{y}} \right\}$, and function $\varphi$ can be constructed in accordance with the ideas mentioned in Section \ref{sec:phi}. Therefore, in line with the results of Theorems \ref{th:bounded} and \ref{th:monotone}, we expect the numerical solution produced by the nonstandard multistep method to be bounded and attain property \eqref{eq:weakmon} for all $\Delta t \geq 0$ step sizes. In this case, Theorem \ref{th:conserv} is not applicable since the equation is one-dimensional (and the preservation of the constant values is trivial).

In the following, we run several numerical tests \rev{analyzing} the performance of the aforementioned numerical methods. For this, we use the strongly stable forms of multistep methods \cite{sspbook} - the coefficients of the methods applied are listed in Appendix A. Moreover, since the solution is known, we used the exact values for those initial values of the numerical method that are needed (for example, for an $s$-step method, one needs the values of $u^0, \; u_1, \; \dots, \; u^{s-1}$). Alternatively, we could have used some further, high-order (possibly nonstandard) numerical method to approximate these values - such an approach is taken in the second example in Section \ref{sec:seir}.

\rev{Initially, we examine the choice $c=2$, and then the case of $c=500$ (resulting in a strict bound for $\varphi$) is also considered.}

\rev{\subsubsection{\texorpdfstring{The non-stiff case of $c=2$}{The non-stiff case of c=2}}}

\begin{figure}[!h]
    \centering
    \includegraphics[width=0.7\linewidth]{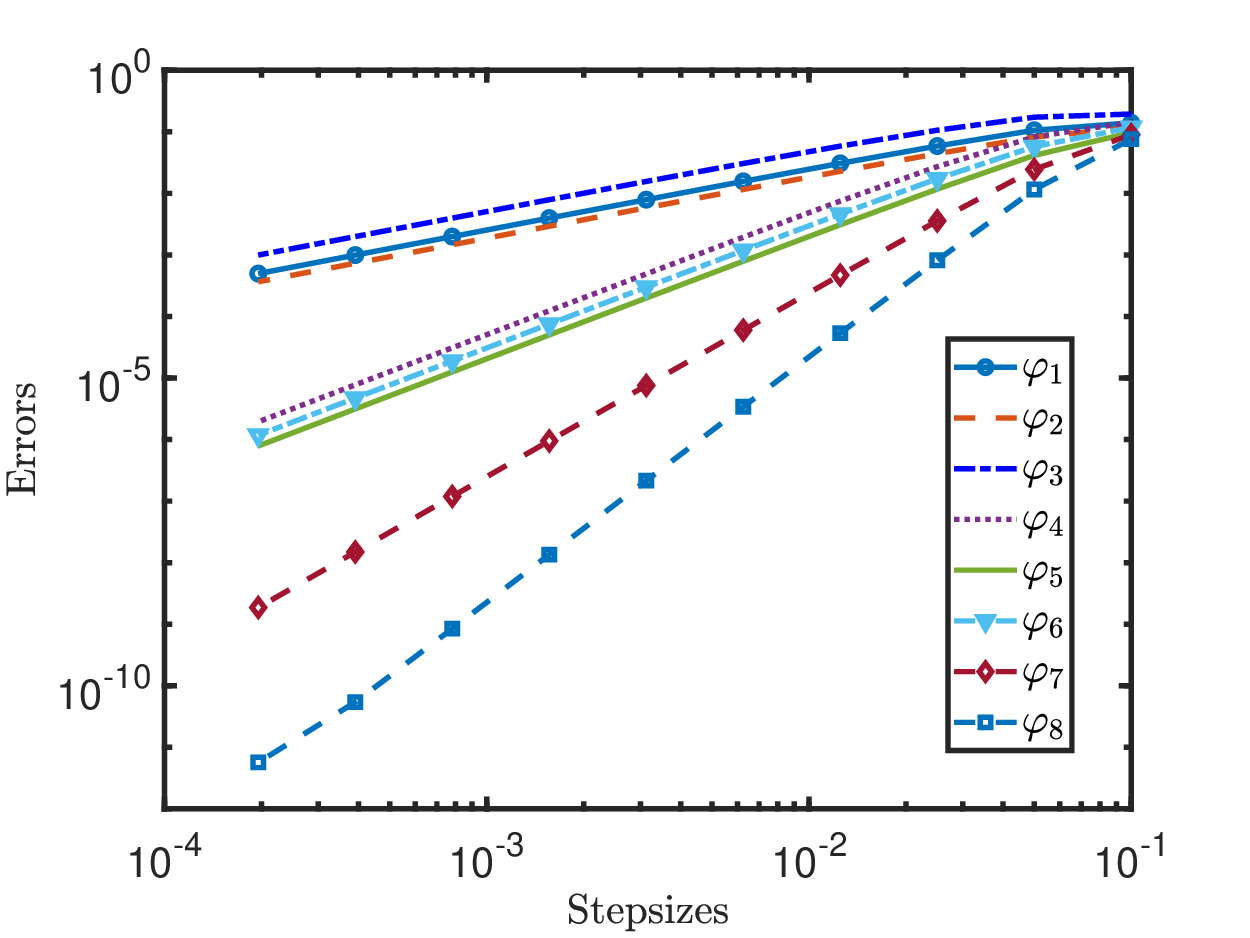}
    \caption{The errors of the different methods with different $\varphi$ functions for different values of $\Delta t$ \rev{for $c=2$}.}
    \label{fig:phi_orders}
\end{figure}

\rev{In this subsection we consider equation \eqref{eq:logistic} with $c=2$.} First, we select the nonstandard version of the fourth-order multistep method SSPMS(6,4) and check the performance of the different choices of the function $\varphi(x)$. We run the methods with timesteps \rev{$\Delta t=0.1/2^k, \; k=0, 1,  \dots, 9$} and initial condition $\tilde{y}=1$, and compare the values of the numerical solutions at $T=1$ to the exact solution. The values of the errors are plotted in Figure \ref{fig:phi_orders}, \rev{while the exact values of the errors and the corresponding orders can be seen in Table \ref{tab:phi_errors} in Appendix B.} From now on, we use the notations introduced in Section \ref{sec:phi}. As we can see, with functions $\varphi_1, \varphi_2$ and $\varphi_3$, the method is only of first order, with $\varphi_4, \varphi_5$ and $\varphi_6$ we attain second order, while the uses of $\varphi_7$ and $\varphi_8$ result in a third- and fourth order methods, respectively. It is also worth mentioning here, that from the first three methods $\varphi_2$ performs the best, while among the second order ones, $\varphi_5$ seems to be the \rev{one with the smallest errors} - this is somewhat anticipated since these are the ones with the smallest error terms among the methods with the same order (see the $\varphi''(0)$ and $\varphi'''(0)$ values in Section \ref{sec:phi}). 

Next, we test the order\rev{s} of the nonstandard multistep methods: for this, we choose $\varphi_8$ \rev{(or $\varphi_7$, see the next paragraph)} and run the nonstandard versions of methods SSPMS(4,2), SSPMS(4,3) and SSPMS(6,4) (\rev{denoted by NSSPMS(4,2), NSSPMS(4,3) and NSSPMS(6,4)}) \rev{along with the nonstandard versions of the SSP Runge-Kutta methods SSPRK(2,2), SSPRK(3,3) and SSPRK(10,4) (denoted by NSSPRK(2,2), NSSPRK(3,3) and NSSPRK(10,4), respectively). The coefficients of these methods can be found in Appendix A}. Here we run the methods with timesteps \linebreak $\Delta t = \rev{0.05/2^k, \; k = 0, 1, \dots, 8}$ (with $\tilde{y}=1$ as before) and compare the numerical solutions to the exact solution at $T=1$. Note that here we have three different values of $\mathcal{C}$, meaning that parameter $\mathcal{B}$ might have different values depending on the method. In Figure \ref{fig:LMM_orders} we plot the errors, while Table \ref{tab:LMM_errors} \rev{in Appendix B.} contains the errors along with the orders. \rev{Apart from the third order ones}, the methods perform as expected. \rev{In general, the multistep and the Runge-Kutta methods have a similar order, while the latter produce smaller errors.}

\begin{figure}[!h]
    \centering
     \includegraphics[width=0.7\linewidth]{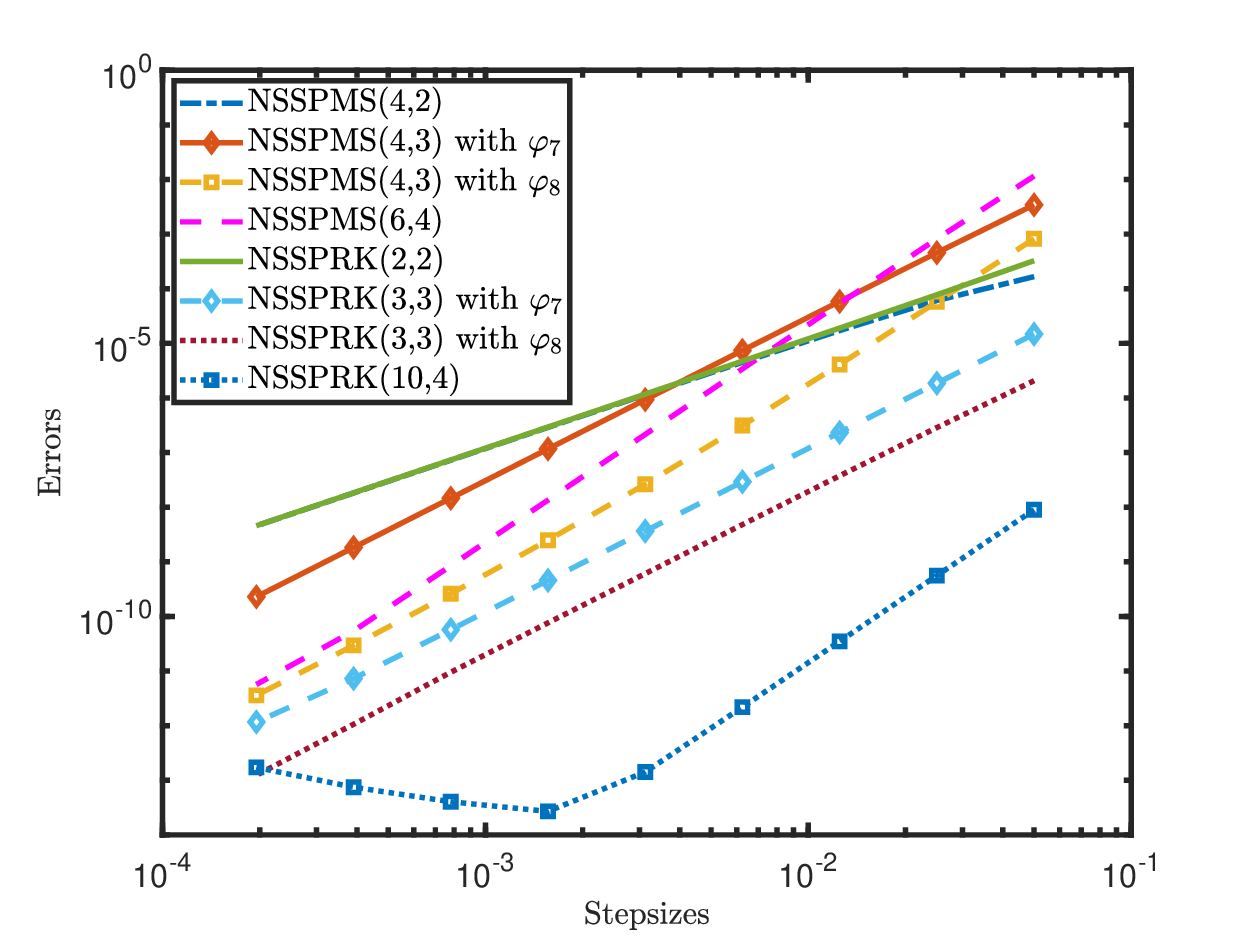}
    \caption{The errors of the different \rev{nonstandard} methods for different values of $\Delta t$ \rev{with $c=2$}.}
    \label{fig:LMM_orders}
\end{figure}

One surprising property is the \rev{order} of the convergence of the third-order method \rev{NSSPMS(4,3)}\rev{: its order is close to $4$ for larger values of $\Delta t$}. This might be the result of the effect of the function $\varphi_8$, since this function enables higher order (see the cases of NSSPMS(6,4) and NSSPRK(10,4)). \rev{This is also proven by the plot of the same method run with the choice $\varphi_7$: in this case, since this function can only enable an order of $3$, the order of the method is decreased by one, resulting in the expected order of $3$.} This might indicate that for carefully chosen $\varphi$ functions, nonstandard multistep methods might have larger orders than their standard counterparts. \rev{Consequently, it is usually advised to use such $\varphi$ functions that enable higher orders (four or even bigger). Note that this effect is not present in the Runge-Kutta case: for NSSPRK(3,3), both of the choices of $\varphi_7$ and $\varphi_8$ result in an order of $3$, although the $\varphi_8$ choice results in smaller errors.}

\begin{figure}[!ht]
    \centering
    \includegraphics[width=0.45\linewidth]{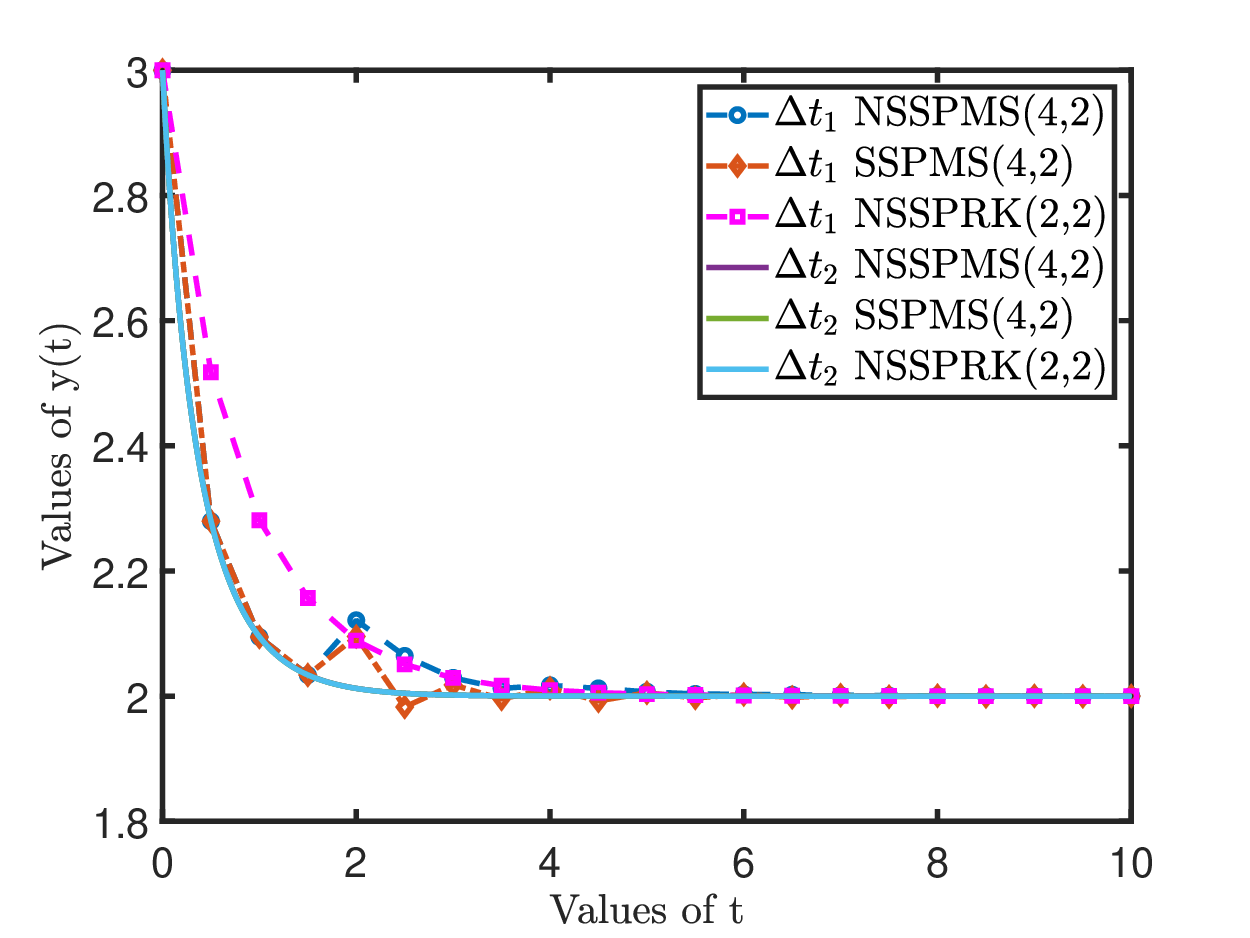}
    \includegraphics[width=0.45\linewidth]{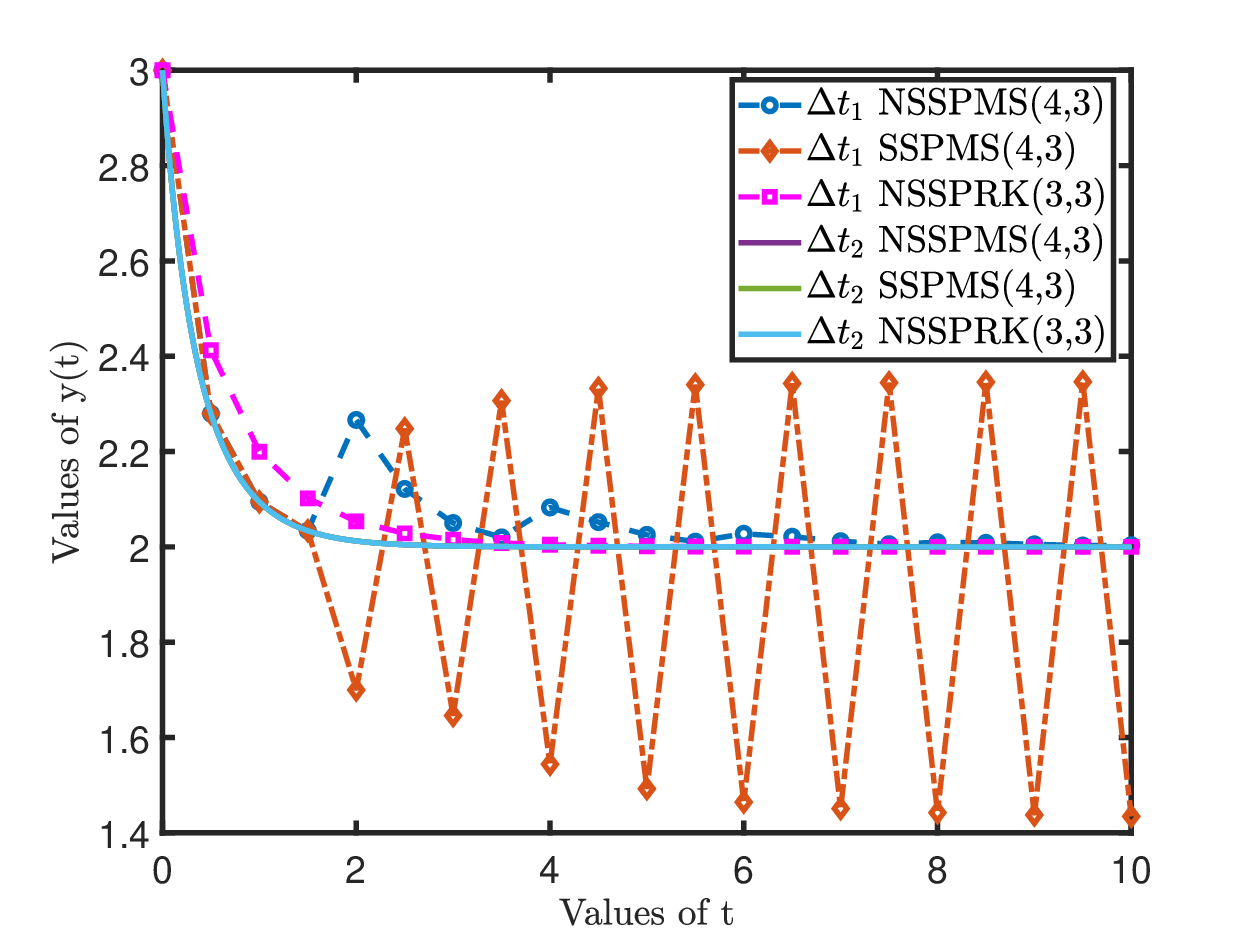}
    \includegraphics[width=0.45\linewidth]{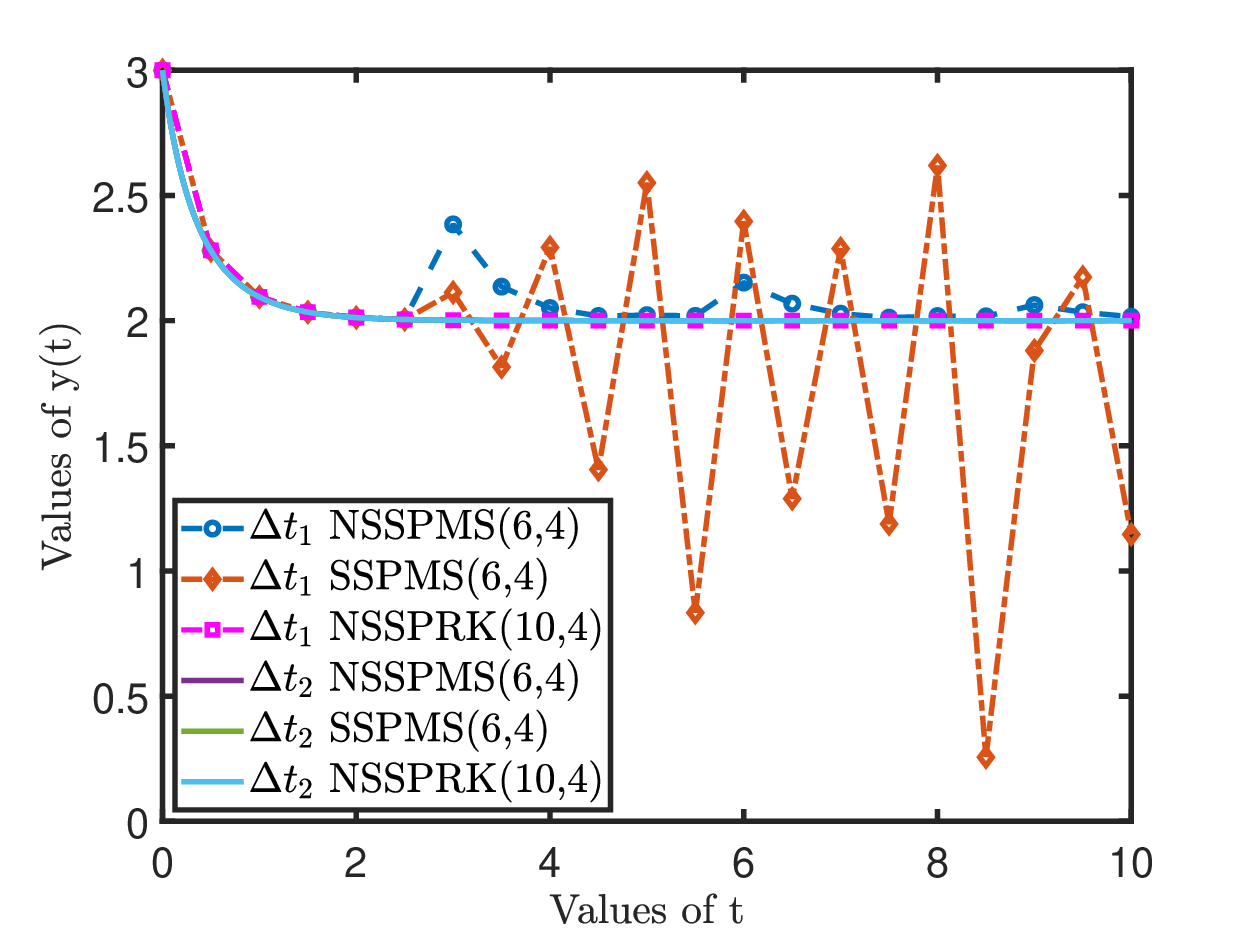}
    \caption{Methods SSPMS(4,2), NSSPMS(4,2) and NSSPRK(2,2) (upper left), SSPMS(4,3), NSSPMS(4,3) and NSSPRK(3,3) (upper right) and SSPMS(6,4), NSSPMS(6,4) and NSSPRK(10,4) (lower) with timesteps $\Delta t_1=0.5$ and $\Delta t_2=0.02$. For $\Delta t_2$, the schemes are very close to each other, therefore, only one of them is visible.}
    \label{fig:conserv}
\end{figure}

Moreover, we \rev{consider} the way the nonstandard methods preserve the qualitative properties of the original differential equation. By Theorems \ref{th:bounded} and \ref{th:monotone}, we know that the nonstandard multistep method preserves the boundedness property of the original solution (since the Forward Euler method has this property for $\Delta t \leq \mathcal{B}$, see Proposition \ref{prop:FE_logistic}), and condition \eqref{eq:weakmon} also holds. In Figure \ref{fig:conserv}, we compare the results of the nonstandard \rev{multistep} methods with their standard counterparts applied with the same step sizes, initial value $\tilde{y}=3$, and function choices $\varphi_5$, $\varphi_7$, and $\varphi_8$, respectively. \rev{On the same figures, we also plot the solutions produced by the nonstandard Runge-Kutta methods of the same order.} As it can be seen, for $\Delta t_1 =0.5$ the standard methods are highly unstable: they are not bounded from below, and \eqref{eq:weakmon} is also violated. However, the nonstandard ones behave as expected, even for a larger timestep: their solutions are bounded from below, and the property \eqref{eq:weakmon} also holds with $s=4$ (in the first two cases) and with $s=6$ (in the last case). For $\Delta t_2=0.02$, \rev{all} methods behave as expected. \rev{In general, the nonstandard multistep methods perform better than their standard versions, but fail to preserve e.g. monotonicity for $\Delta t_1$, while the nonstandard Runge-Kutta ones achieve even this goal.}

\begin{figure}[!ht]
    \centering
    \includegraphics[width=0.45\linewidth]{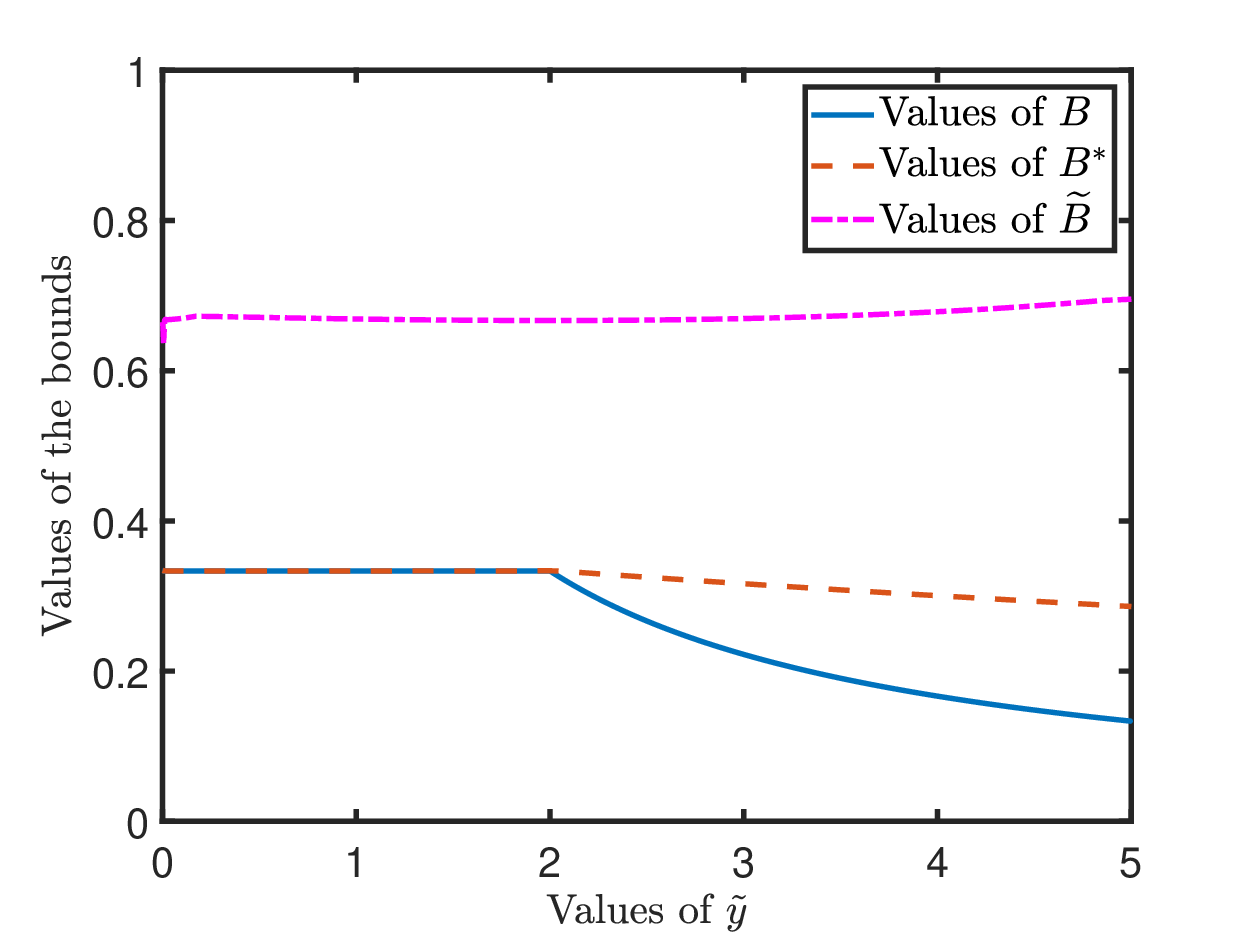}
    \includegraphics[width=0.45\linewidth]{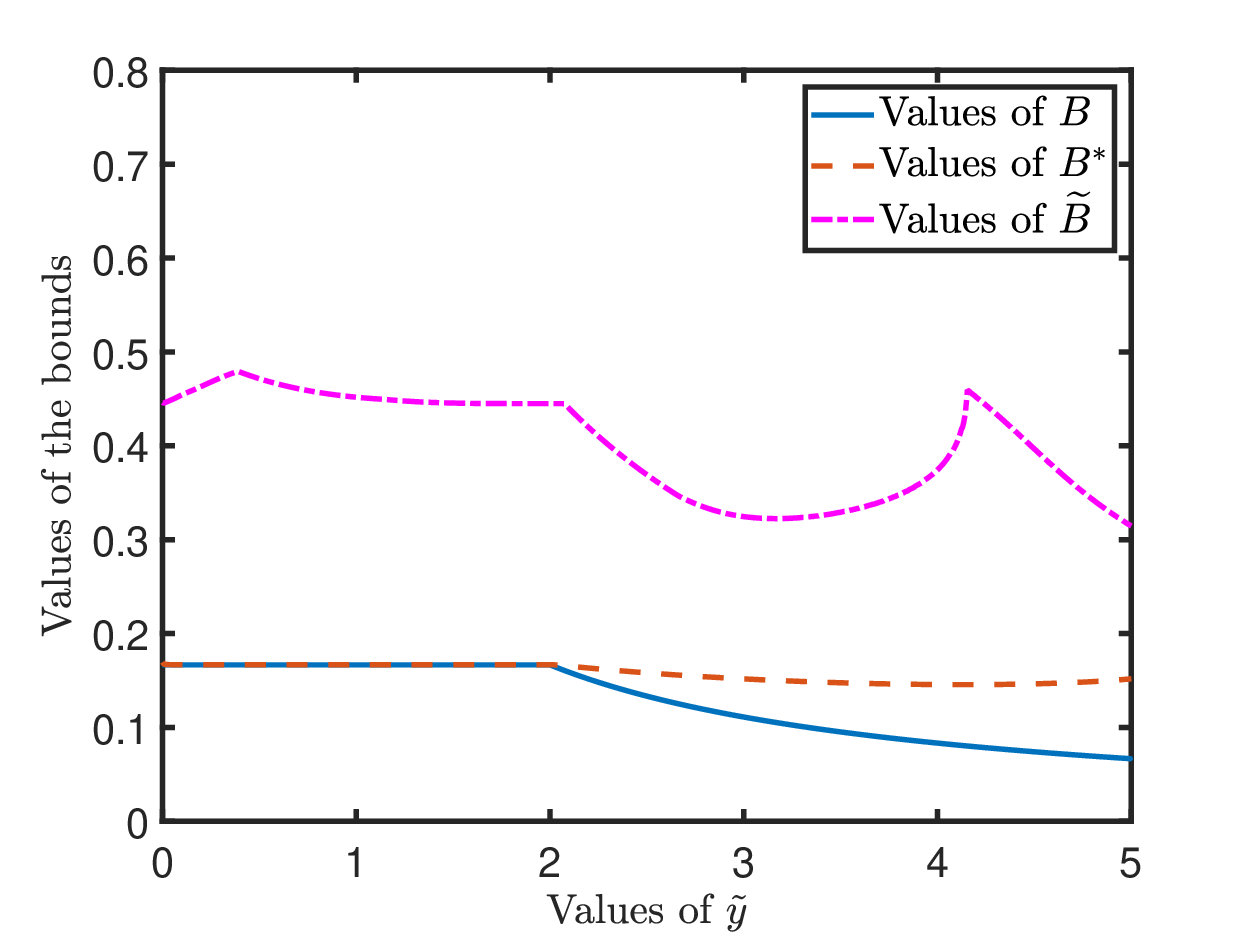}
    \includegraphics[width=0.45\linewidth]{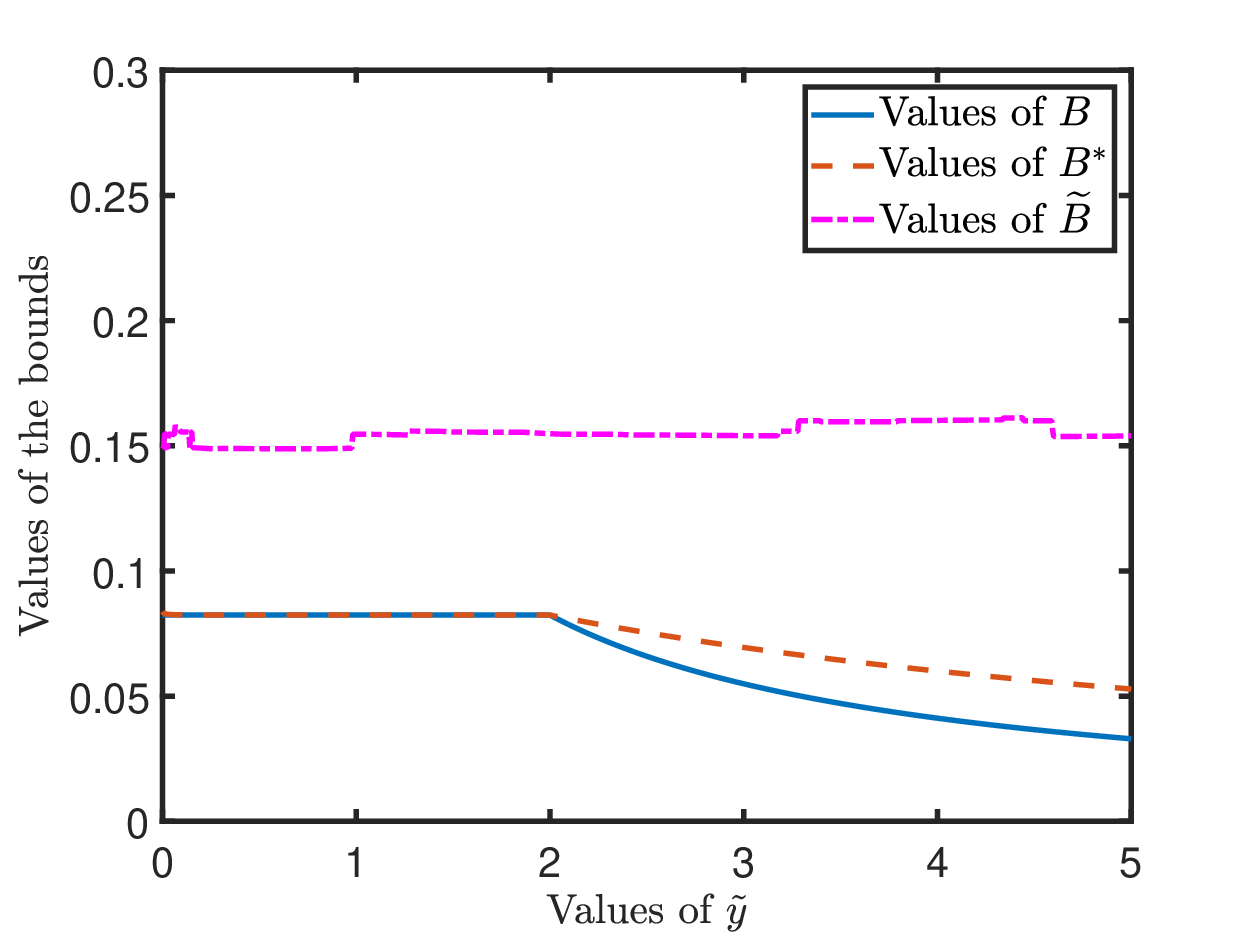}
    \caption{\rev{The sufficient bound $\mathcal{B}$ and the ``real'' bounds for methods NSSPMS(4,2) (upper left), NSSPMS(4,3) (upper right) and NSSPMS(6,4) (lower) for the boundedness property ($\mathcal{B}^*$) and the weak monotonicity ($\widetilde{\mathcal{B}}$) with $c=2$.}}
    \label{fig:dahl_boundsharp}
\end{figure}

\rev{Lastly, we observe the sharpness of bound $\mathcal{B}$. The conditions in Theorems \ref{th:bounded}, \ref{th:monotone} and \ref{th:conserv} are sufficient, but might not be necessary. Because of this, it might happen that the use of functions $\varphi(x)$ for which only $\varphi(x)\leq \mathcal{B}^*$ holds (where $\mathcal{B}<\mathcal{B}^*$) results in numerical solutions that behave in a qualitatively reasonable way too. In the following, we test the sharpness of the bound $\mathcal{B}$ calculated previously, namely, how close the sufficient bound $\mathcal{B}$ and the \say{real} bound (determined by numerical experiments) are. For this, we change the bound for $\varphi$ to an arbitrary positive constant, and test the methods with $1000$ values of $\tilde{y}$ from the interval $[10^{-3}, \; 5]$ and consider $1000$ possible values of timesteps in the interval $[0.5,\; 3]$. We say that the method behaves as expected for a fixed value of $\tilde{y}$ and with this given bound for $\varphi$, if the given property is preserved for every possible timestep on the interval $t \in [0,100]$. By a bisection method, we determine the bound $\mathcal{B}^*$ for the boundedness property (i.e., $u^n \leq 0$ if $\tilde{y}\geq 2$ and $0 \leq u^n \leq 2$ if $0 \leq \tilde{y} \leq 2$ for every $n$) and bound $\widetilde{\mathcal{B}}$ for the weak monotonicity property \eqref{eq:weakmon}. For methods NSSPMS(4,2), NSSPMS(4,3) and NSSPMS(6,4), functions $\varphi_5$, $\varphi_7$ and $\varphi_8$ were used, respectively. As we can see in Figure \ref{fig:dahl_boundsharp}, the bound for the boundedness property is very close to the sufficient bound given by Theorem \ref{th:bounded} when $\tilde{y} \in (0,2)$, while the bound for the weak monotonicity property \eqref{eq:weakmon} is considerably bigger than the bound given by Theorem \ref{th:monotone}. However, since we would like to preserve both of these during our numerical method, we can say that for $\tilde{y} \in (0,2)$, the very small difference between $\mathcal{B}$ and $\mathcal{B}^*$ makes the bound relatively sharp, while for $\tilde{y}>2$, the bound $\mathcal{B}$ is also not that far from the ``real'' bound. It should also be mentioned that surprisingly, for the ``strong'' monotinicity property (i.e. $u^{n+1}\leq u^n$ if $\tilde{y}>2$ and $u^{n+1}\geq u^n$ if $\tilde{y} \in ]0, \; 2[$) we can always find a timestep $\Delta t \in [0.5, \; 3]$ for which this condition is violated, meaning that this property is never satisfied for every possible value of $\Delta t$, even if we decrease the bound of $\varphi(x)$ drastically.}

\rev{\subsubsection{\texorpdfstring{The stiff case of $c=500$}{The stiff case of c=500}}}

\rev{In this section, we analyze the performance of the nonstandard multistep method applied to a stiff equation. For this, we consider the previous equation \eqref{eq:logistic} with $c=500$. By Proposition \ref{prop:FE_logistic}, if $\varphi(x)\leq \mathcal{C} \min \left\{ \frac{1}{500}, \frac{1}{\tilde{y}} \right\}$, then the method preserves the boundedness and property \eqref{eq:weakmon} for any choice of $\Delta t$.}

\begin{figure}[!ht]
    \centering
    \includegraphics[width=0.7\linewidth]{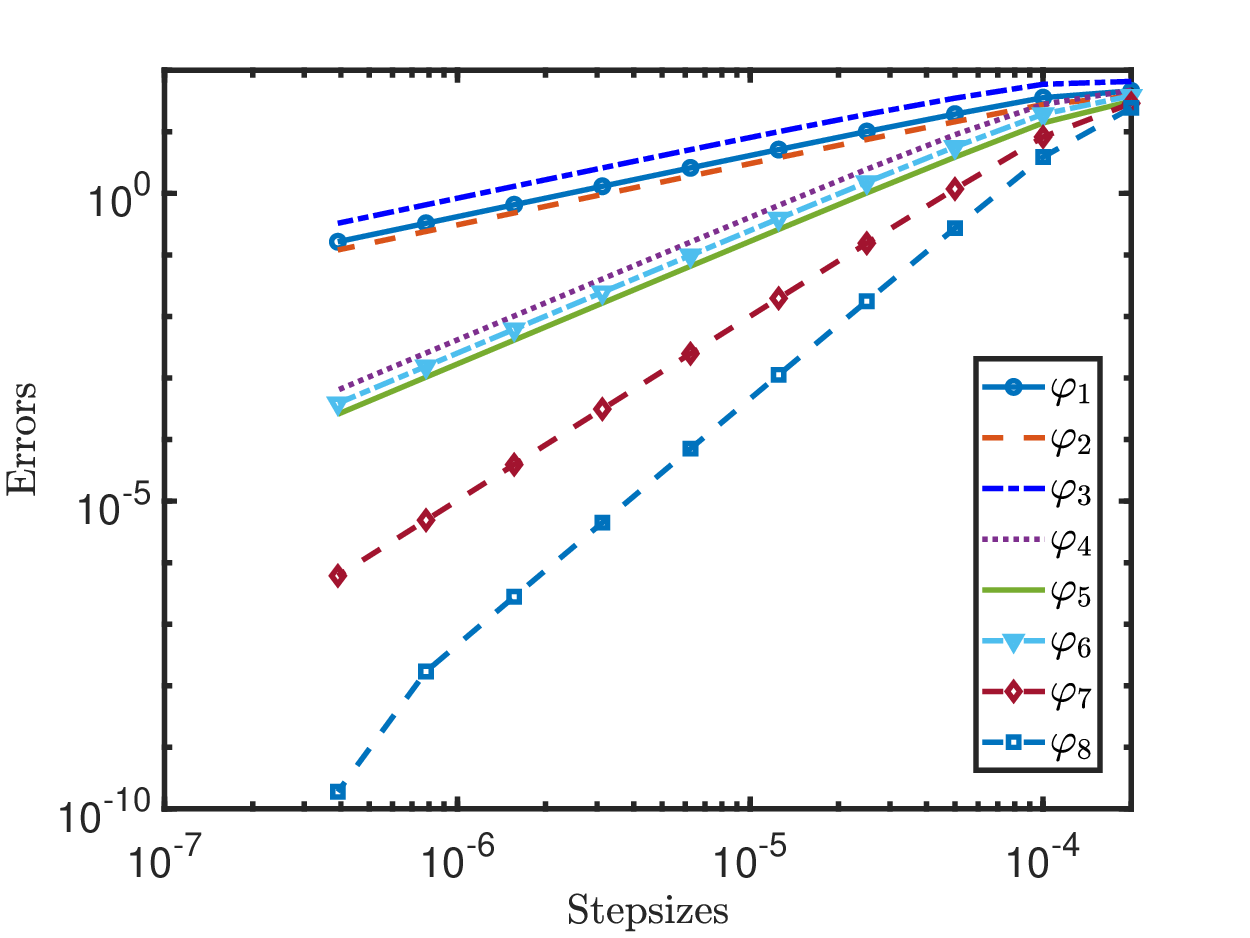}
    \caption{\rev{The errors of the NSSPMS(5,4) method with different $\varphi$ functions for different values of $\Delta t$ for $c=500$.}}
    \label{fig:phi_orders_c500}
\end{figure}

\rev{We run similar tests as in the case of $c=2$ - namely, we first fix the numerical method as the fourth order NSSPMS(6,4) and observe the errors of the different choices of function $\varphi$. For this, we choose $\tilde{y} = 2 c = 1000$ and compare the errors of the different choices at $T=\frac{1}{c} = \frac{1}{500}$ with time steps $\Delta t=2^{-k} \frac{1}{10 c}, \; k = 0, \; 1, \; \dots \; 9$. The reason for the seemingly small value of $T$ is that for bigger ones, all the numerical methods are close to the exact solution (since the function is almost constant for bigger values of $t$, meaning that its change is smaller than machine precision). In Figure \ref{fig:phi_orders_c500}, we can see that we got a plot very similar to Figure \ref{fig:phi_orders}, so the increase of $c$ does not affect the order of the method, but the errors are significantly bigger. The corresponding values of the errors and the orders can be seen in Table \ref{tab:phi_errors_c500} in Appendix B.}

\begin{figure}[!ht]
    \centering
     \includegraphics[width=0.7\linewidth]{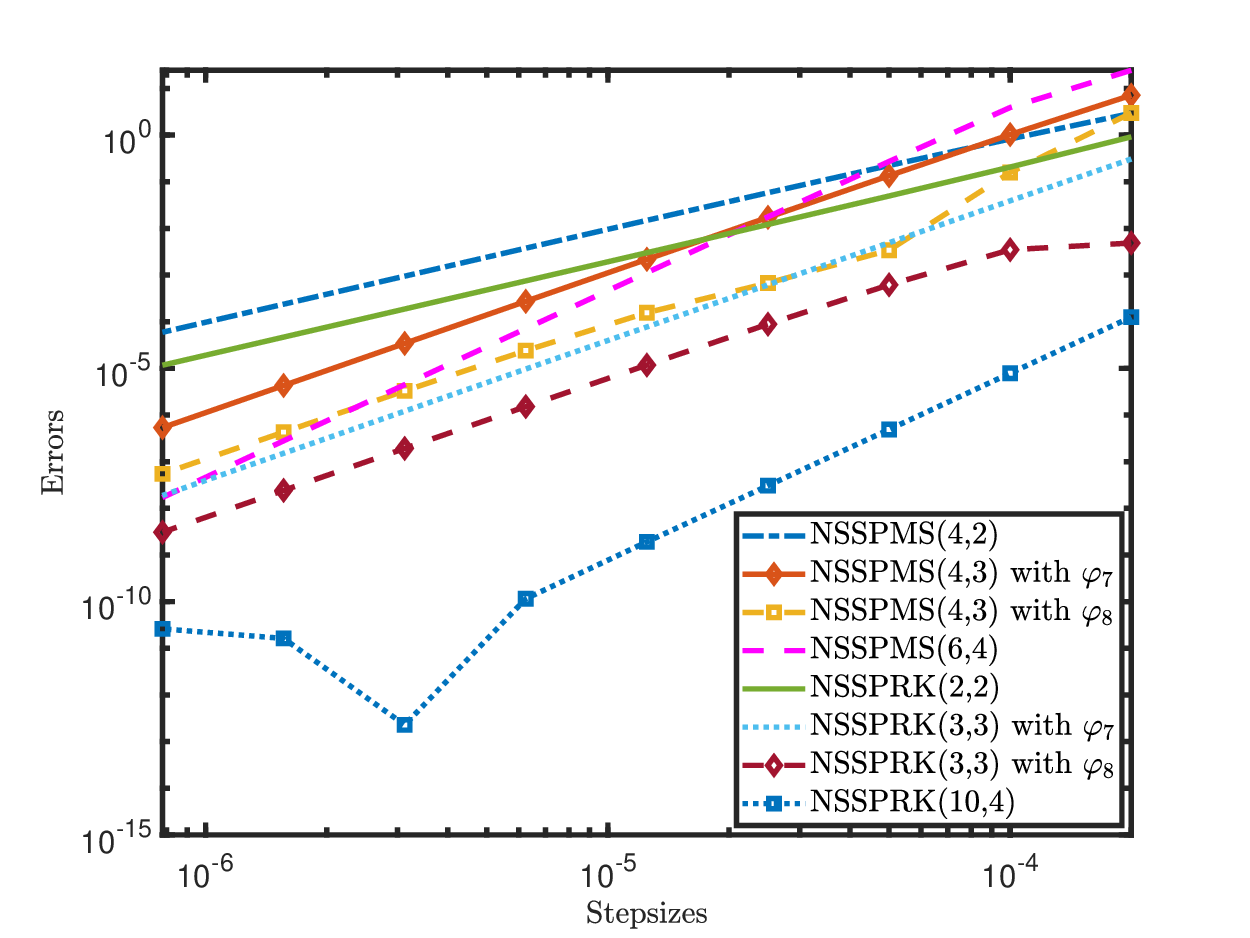}
    \caption{\rev{The errors of the different nonstandard methods for different values of $\Delta t$ with $c=500$}.}
    \label{fig:LMM_orders_c500}
\end{figure}

\rev{Moreover, we also examine the order of the different multistep methods when we use $\varphi_8$, along with the nonstandard Runge-Kutta ones. For this, we use $\tilde{y}=2 c =1000$ and compare the errors of methods NSSPMS(4,2), NSSPMS(4,3) and NSSPMS(6,4) along with Runge-Kutta ones NSSPRK(2,2), NSSPRK(3,3) and NSSPRK(10,4) with timesteps $\Delta t=2^{-k} \frac{1}{10 c}, \; k = 0, \; 1, \; \dots \; 9$ at $T=\frac{1}{c} = \frac{1}{500}$. In Figure \ref{fig:LMM_orders_c500}, it can be seen again that the methods perform as expected, even for $c=500$, but with larger errors. The corresponding values of the errors and the orders can be found in Table \ref{tab:LMM_errors_c500} in Appendix B. The orders are similar to the ones of the case $c=2$: the order of the method NSSPMS(4,3) is considerably larger for bigger values of $\Delta t$ when we use $\varphi_8$, but it is reduced when we use $\varphi_7$, while for the Runge-Kutta ones the order is three in both cases.}

\begin{figure}[!ht]
    \centering
    \includegraphics[width=0.45\linewidth]{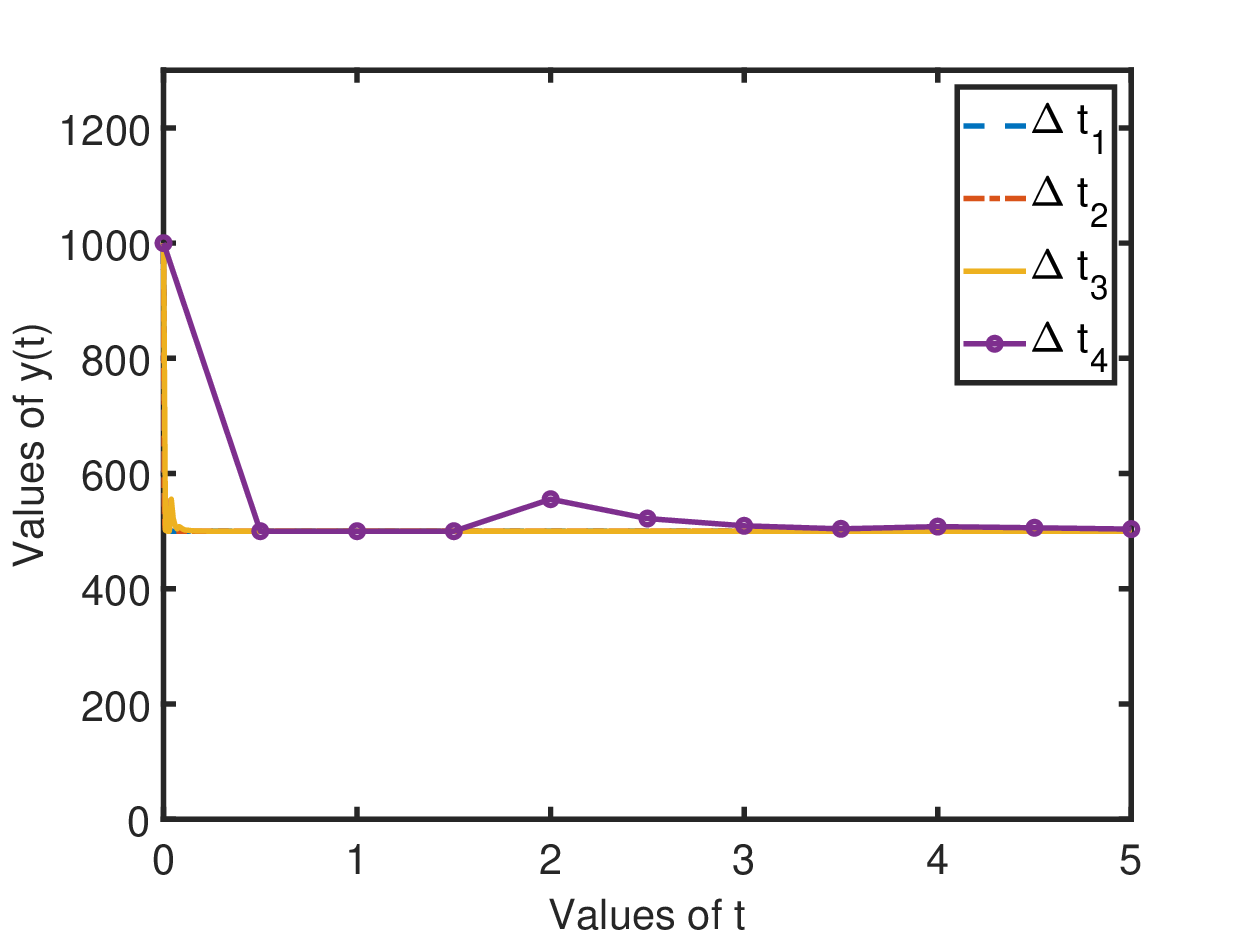}
    \includegraphics[width=0.45\linewidth]{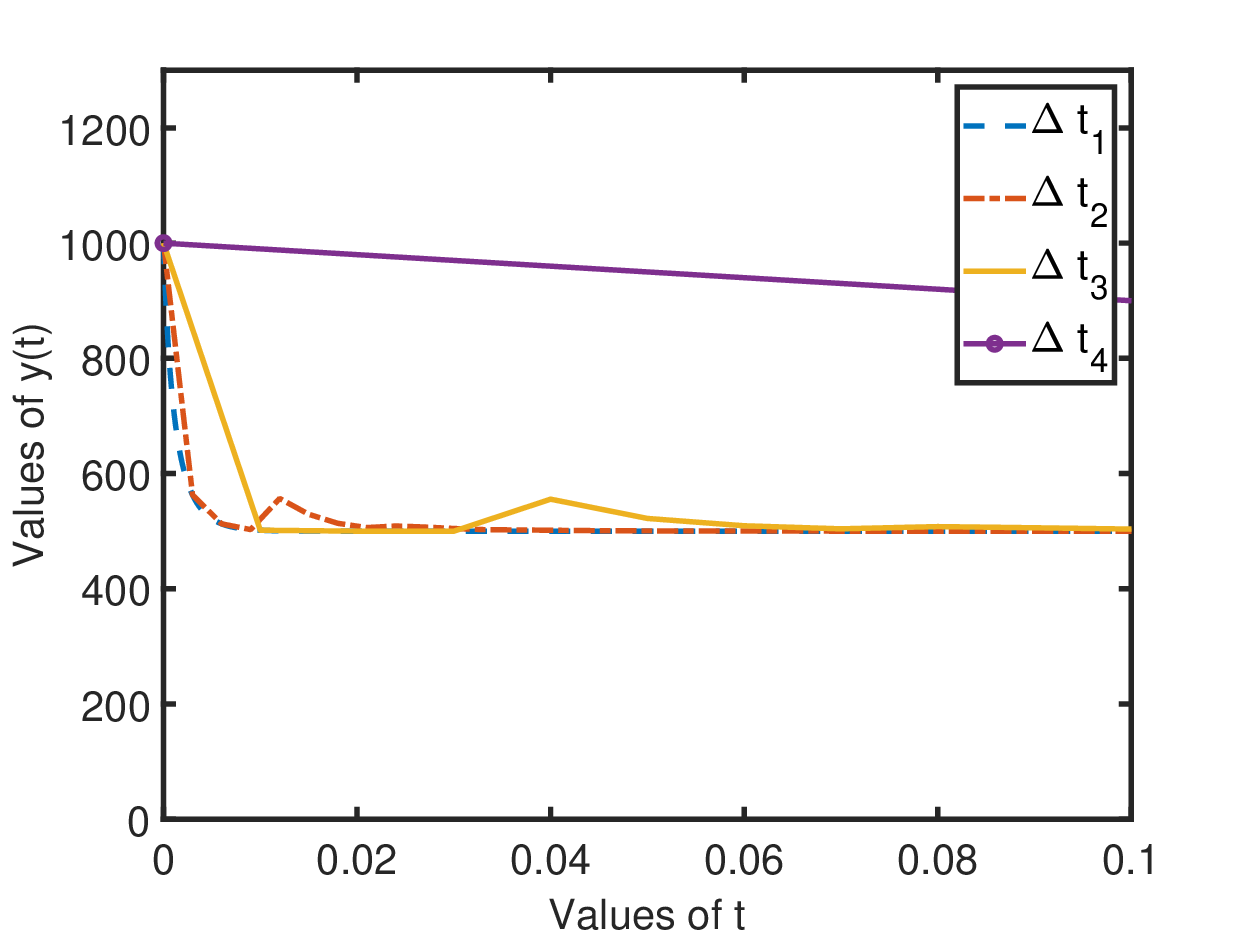}
    \includegraphics[width=0.45\linewidth]{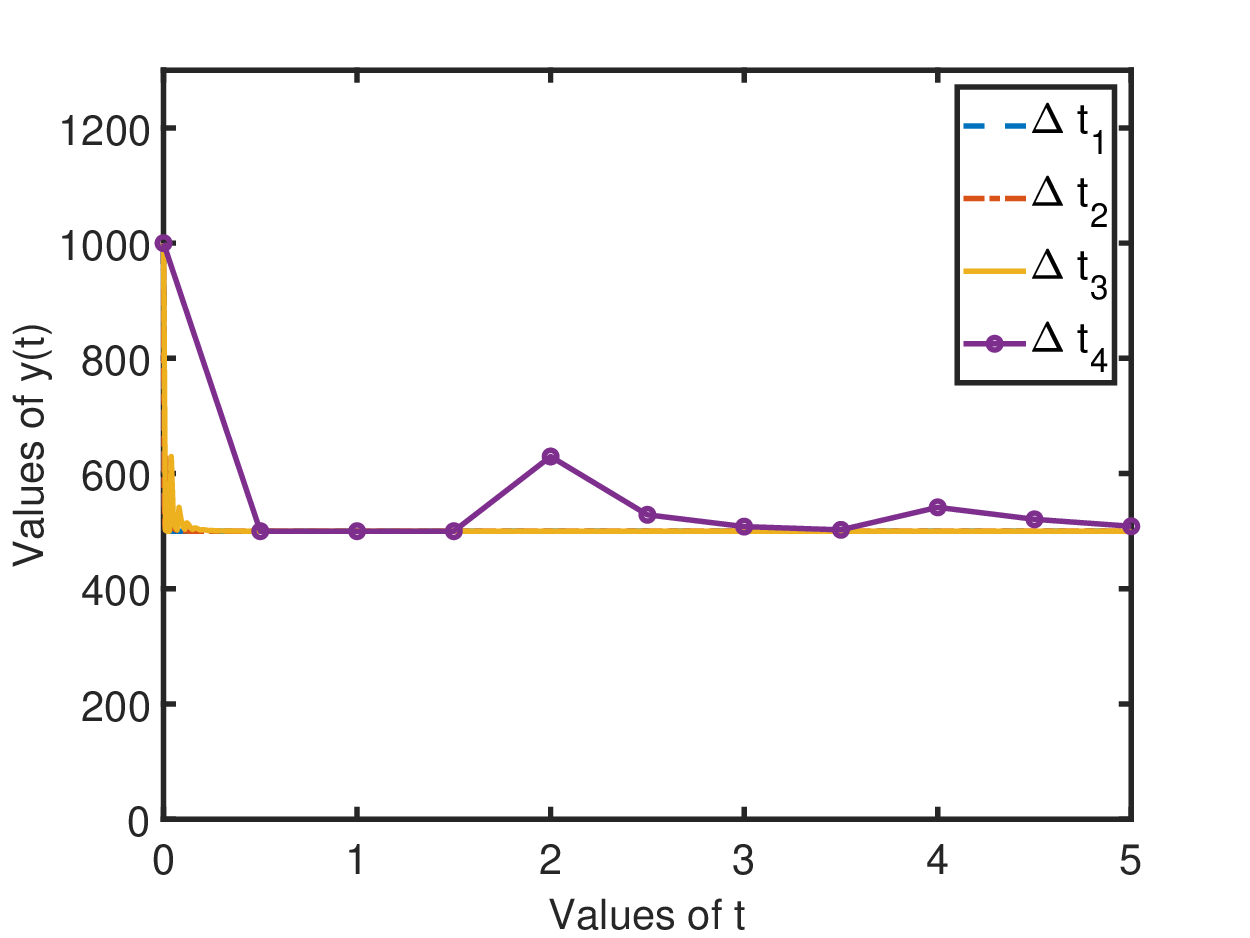}
    \includegraphics[width=0.45\linewidth]{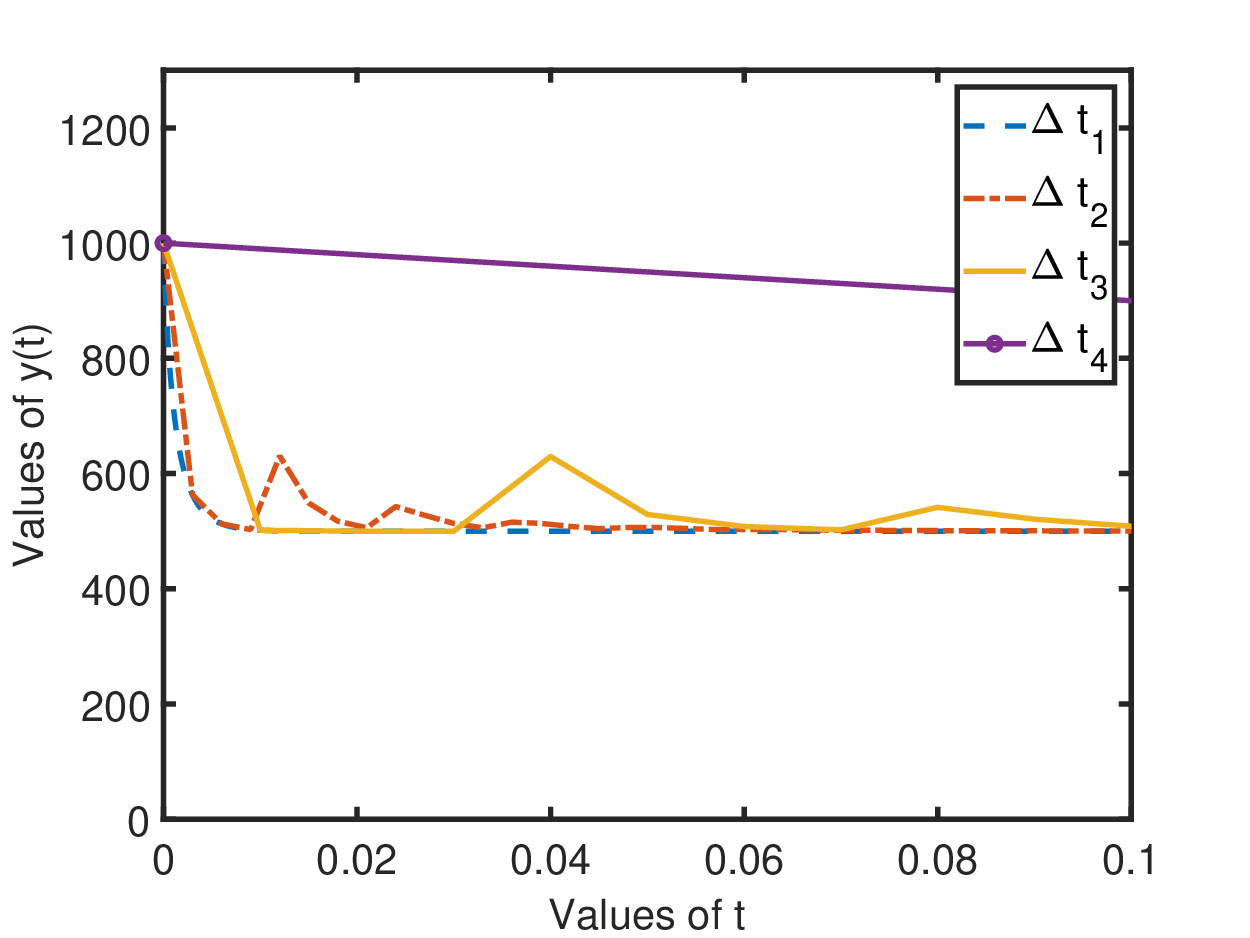}
    \includegraphics[width=0.45\linewidth]{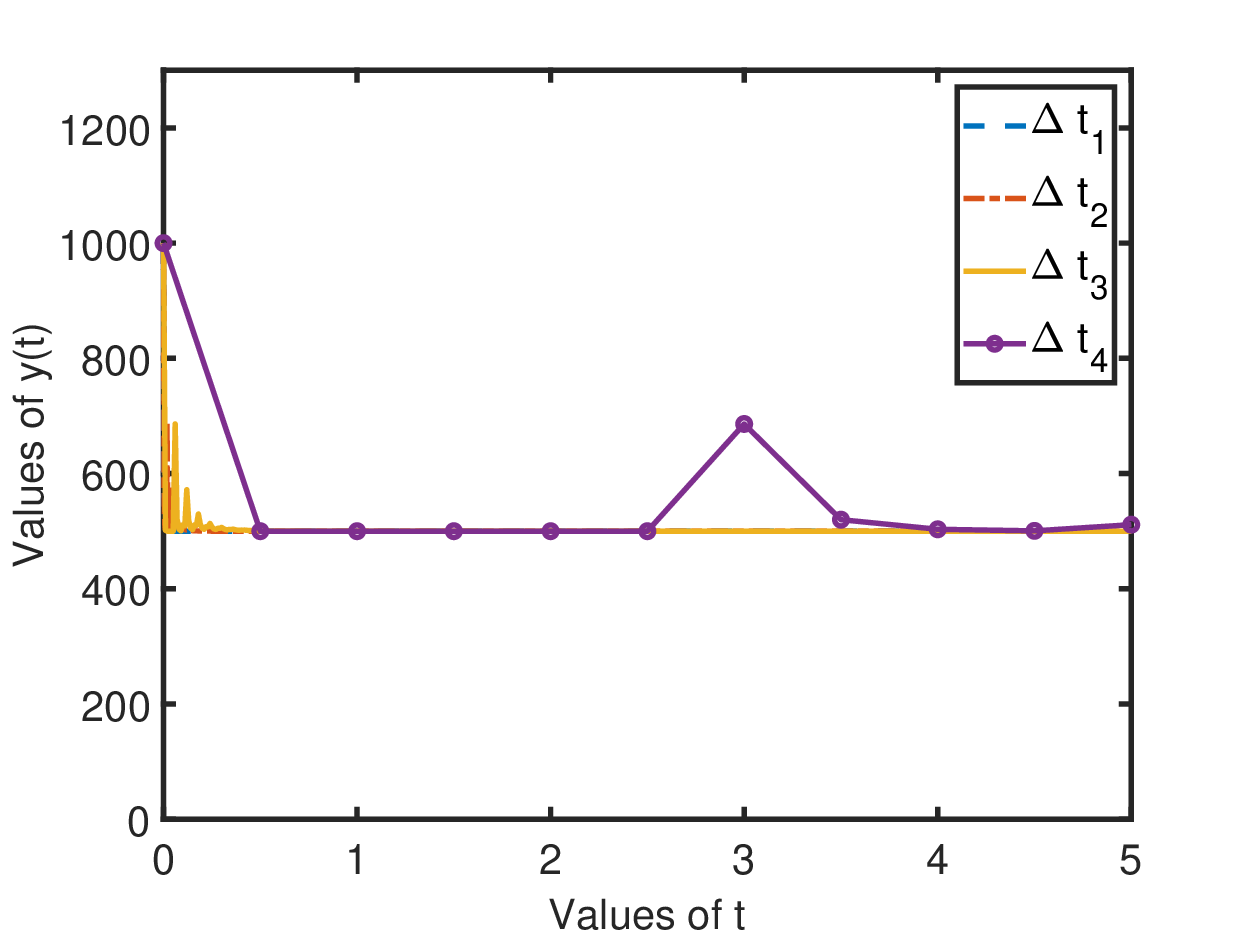}
    \includegraphics[width=0.45\linewidth]{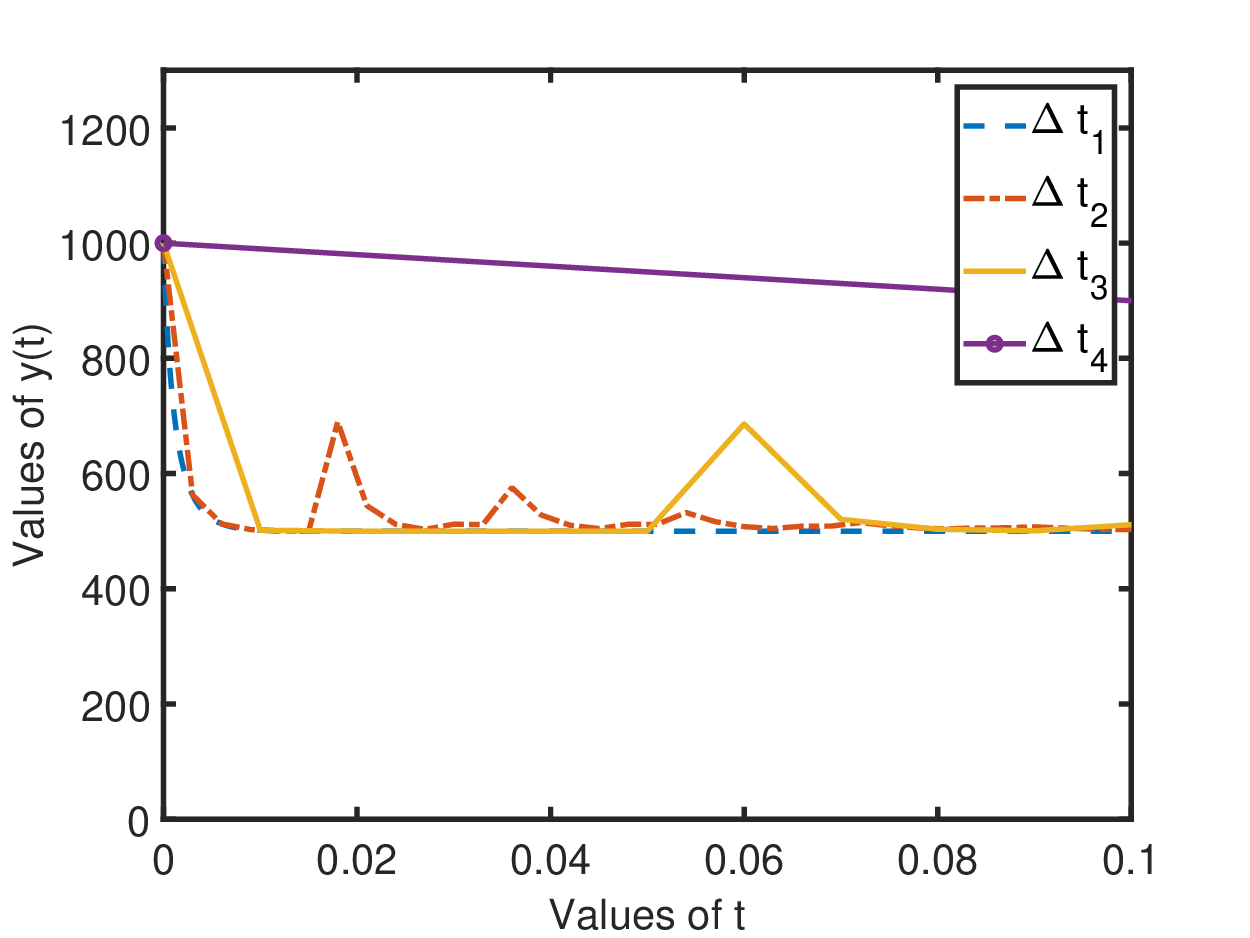}
    \caption{\rev{The numerical solutions of the logistic equation produced by methods NSSPMS(4,2) (top), NSSPMS(4,3) (middle) and NSSPMS(6,4) (bottom) with timesteps $\Delta t_1=10^{-5}$, $\Delta t_2 =3 \cdot 10^{-3}$, $\Delta t_3 = 10^{-3}$ and $\Delta t_4 = 0.5$ with $c=500$. The plots on the right are zoomed-in versions of the left ones.}}
    \label{fig:stiff_compare}
\end{figure}

\rev{As it was mentioned before, the methods preserve the desired qualitative properties even for larger timesteps, but these would result in bigger errors. In Figure \ref{fig:stiff_compare}, we compare the three nonstandard multistep methods examined before. As it can be seen, while the qualitative properties are preserved in all cases, the methods produce significant errors when we use larger timesteps. Note that in this case the ``standard'' multistep methods would only preserve these properties when $\Delta t \leq 2 \cdot 10^{-3}$, while in the nonstandard case these are still preserved for much larger timesteps ($\Delta t_4$ is 100 times bigger than the bound mentioned previously). In general, one can say that the use of bigger timesteps is not efficient in neither the \say{standard} nor the nonstandard cases, but for the nonstandard ones, the solutions at least behave as expected, which could give us some ideas about some properties of the continuous solution.}

\begin{figure}[!ht]
    \centering
    \includegraphics[width=0.45\linewidth]{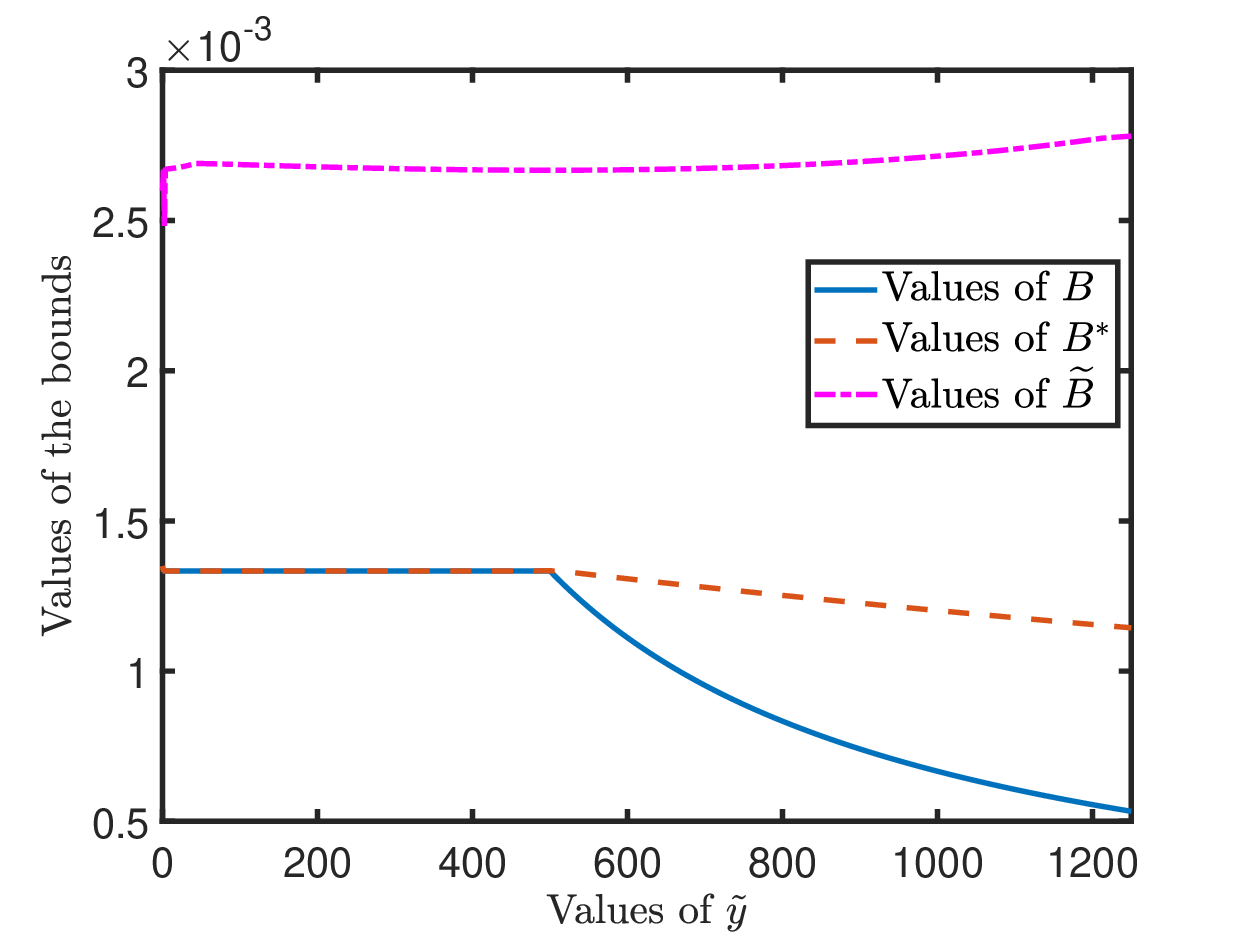}
    \includegraphics[width=0.45\linewidth]{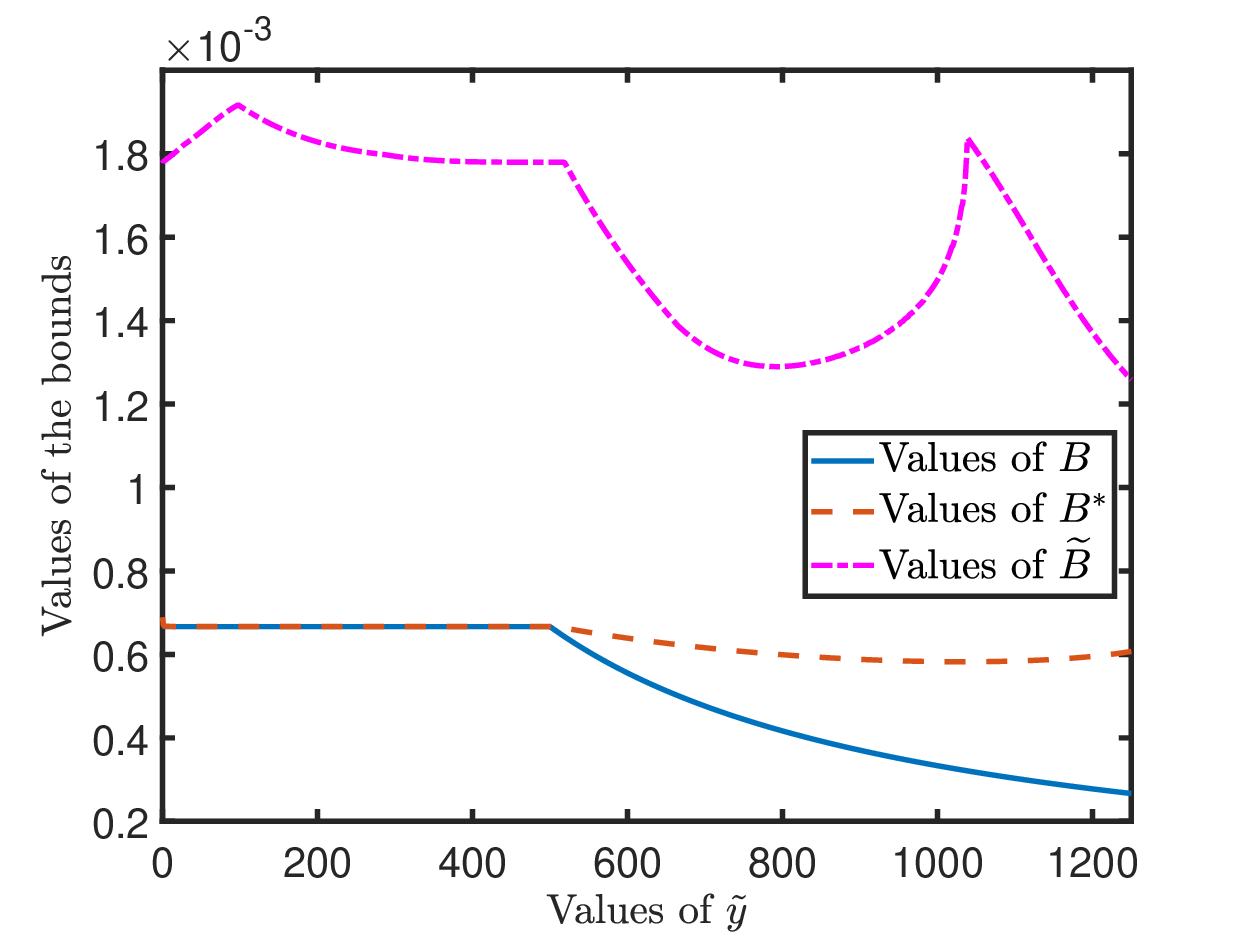}
    \includegraphics[width=0.45\linewidth]{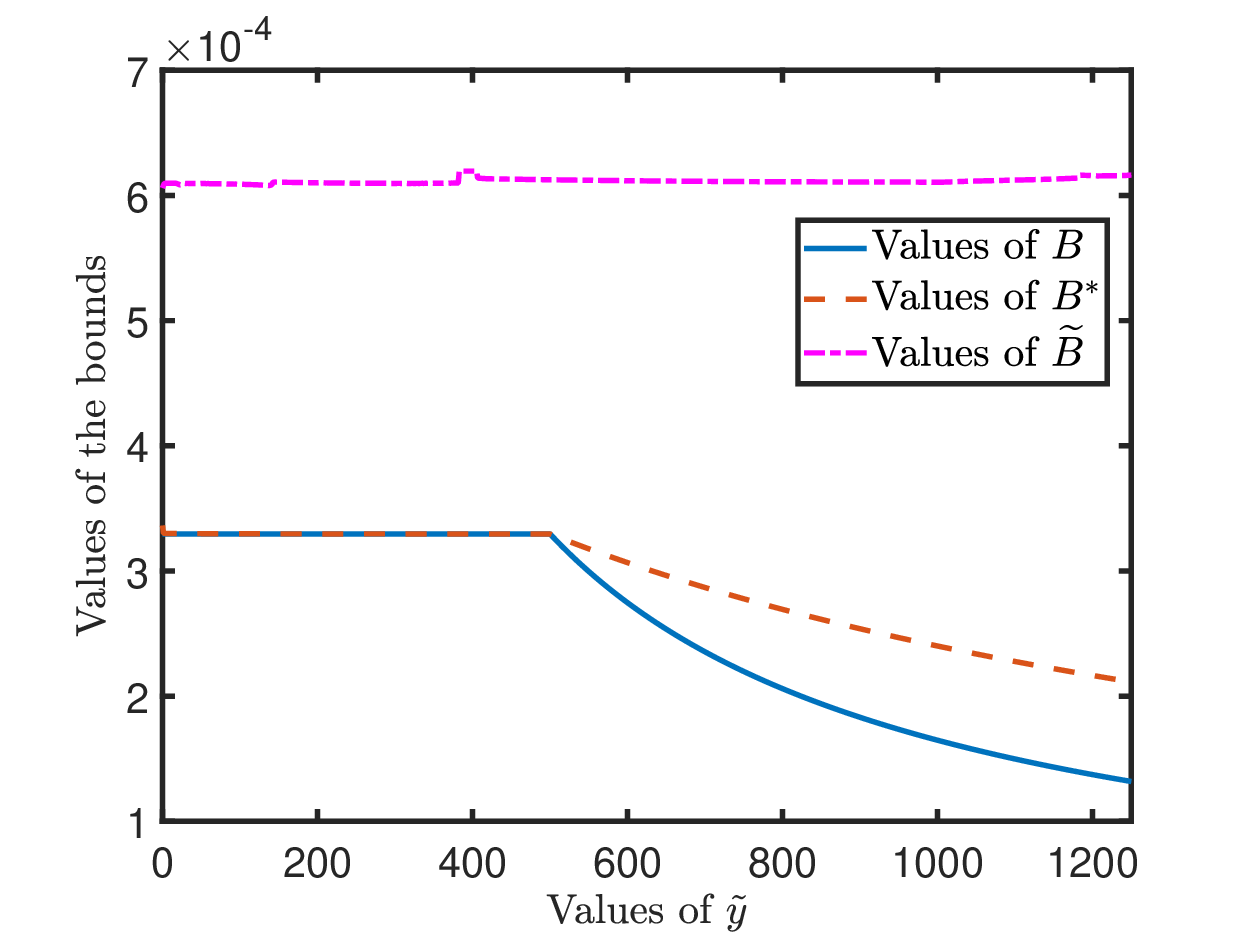}
    \caption{\rev{The sufficient bound $\mathcal{B}$ and the ``real'' bounds for methods NSSPMS(4,2) (upper left), NSSPMS(4,3) (upper right) and NSSPMS(6,4) (lower) for the boundedness property ($\mathcal{B}^*$) and the weak monotonicity ($\widetilde{\mathcal{B}}$) with $c=500$.}}
    \label{fig:dahl_boundsharp_c500}
\end{figure}

\rev{Lastly, we also examine how close the sufficient bound $\mathcal{B}$ is to the bound under which the methods preserve the desired properties. Like in the $c=2$ case, we consider $1000$ values of $\tilde{y}$, but in this case we use the interval $[10^{-3}, \; 1250]$, and test the methods with $1000$ timesteps from the interval $\left[\frac{1}{500}, \frac{6}{500}\right]$, and we use $T=10$. As it can be seen in Figure \ref{fig:dahl_boundsharp_c500}, the plots are very similar to the bounds of the case $c=2$, namely, the bound for the boundedness property can be considered sharp when $\tilde{y} \in (0, 500)$, while one could also use a less strict bound for $\varphi$ in the case of $\tilde{y}>500$. Also, the weak monotonicity property \eqref{eq:weakmon} is preserved for $\varphi$ functions with much higher bounds.}

\subsection{The SEIR system}\label{sec:seir}
The following form of the Kermack-McKendrick model \cite{kermack} describes the propagation of an illness among a given population \rev{with a constant influx of healthy people $\Pi \in \mathbb{R}^+$}:
\begin{equation}\label{sir}
\begin{aligned}
    S'(t) &= \rev{\Pi}-5 S(t) I(t),\\
    E'(t) &= 5 S(t) I(t) - E(t),\\
    I'(t) &= E(t) - I(t),\\
    R'(t) &= I(t),
\end{aligned} 
\end{equation}
where $S(t), E(t), I(t)$ and $R(t)$ describe the number of susceptible, exposed (infected but not yet infectious), ill (infected and infectious), and recovered people, respectively. It is easy to see that if $S(t_0), E(t_0), I(t_0), R(t_0)\geq 0$, then $S(t), E(t), I(t), R(t)\geq 0$ holds for every $t \in [t_0, T]$. Moreover, $S(t)+E(t)+I(t)+R(t) = \rev{ \Pi t + \mathcal{M}}$ is also true \rev{with constant $S(t_0) + E(t_0) + I(t_0) + R(t_0) = \mathcal{M} \in \mathbb{R}^+$} for every $t\in [t_0, T]$ (an interested reader might consult \cite{capasso}). \rev{From now on, we use $t_0=0$.}

Next, let us state a result about the behavior of the forward Euler method, which will be used in the choice of the function $\varphi(x)$.

\begin{prop}\label{prop:seir_bound}
    Let us apply the forward Euler method to equation \eqref{sir} with stepsize $\Delta t_{FE}\leq \min\left\{ \rev{\dfrac{1}{5\mathcal{M}}}, 1 \right\}$. Then, the forward Euler method applied to equation \eqref{sir} preserves the non-negative property of the numerical solution, and \rev{if $\Pi = 0$, then} $S^n+E^n+I^n+R^n=\mathcal{M}$ also holds.
\end{prop}
The proof of a more general case can be found in \cite{takacs}. Consequently, the choice $\mathcal{B}=\mathcal{C} \min\left\{ \rev{\dfrac{1}{5\mathcal{M}}}, 1 \right\}$ is appropriate.


\rev{In \cite{takacs} it was shown that \say{standard} Runge-Kutta methods written in the strongly stability preserving form preserve positivity and the sum of the components if \linebreak $\Delta t \leq \mathcal{C} \min\left\{ \dfrac{1}{5\mathcal{M}}, 1 \right\}$ where $\mathcal{C}$ is the SSP coefficient of the Runge-Kutta method. Therefore, for the initial $1, 2, \dots, s-1$ values we can use a nonstandard Runge-Kutta method with a suitable order. The forms of the SSP-Runge-Kutta methods used during the experiments are listed in Appendix A.}

\rev{First, the case $\Pi=0$ is considered, and then the case of nonzero influx is also observed.}

\rev{\subsubsection{\texorpdfstring{The case $\Pi=0$}{The case Pi=0}}}

In the following, we apply the same tests to this system as we did with the previous example\rev{s in Section \ref{sec:log} with choice $\Pi=0$}. First, we construct a reference solution by applying the (standard) fourth-order Runge-Kutta method with a very small timestep ($\Delta t=0.1 \cdot 2^{-10} \cdot 10^{-3}$). Then, we run the nonstandard multistep method NSSPMS(6,4) with $(S^0, E^0, I^0, R^0)=(0.8, 0, 0.2, 0)$ and with the different choices of function $\varphi$ (see Section \ref{sec:phi}) with timesteps $\Delta t = \rev{0.1/2^k, \; k=0, 1, \dots, 10}$. For the values at points \linebreak $t=\Delta t, \; 2 \Delta t, \; \dots, \; 5 \Delta t$, we used the nonstandard version of the fourth-order SSPRK(10,4) method with function $\varphi_8$. Finally, we \rev{analyze} the difference between these solutions and the reference solution at $T=5$. The errors and the corresponding orders can be seen in Figure \ref{fig:seir_phiconv} and the values are listed in Table \ref{tab:seir_phi_errors} \rev{in Appendix B}. As we can see, the methods behave as expected, and as in the case of the first example, the first-order method \rev{with the smallest error} is the one produced by $\varphi_2$, and the second-order one \rev{with the smallest error} is $\varphi_5$.

\begin{figure}[!ht]
    \centering
     \includegraphics[width=0.7\linewidth]{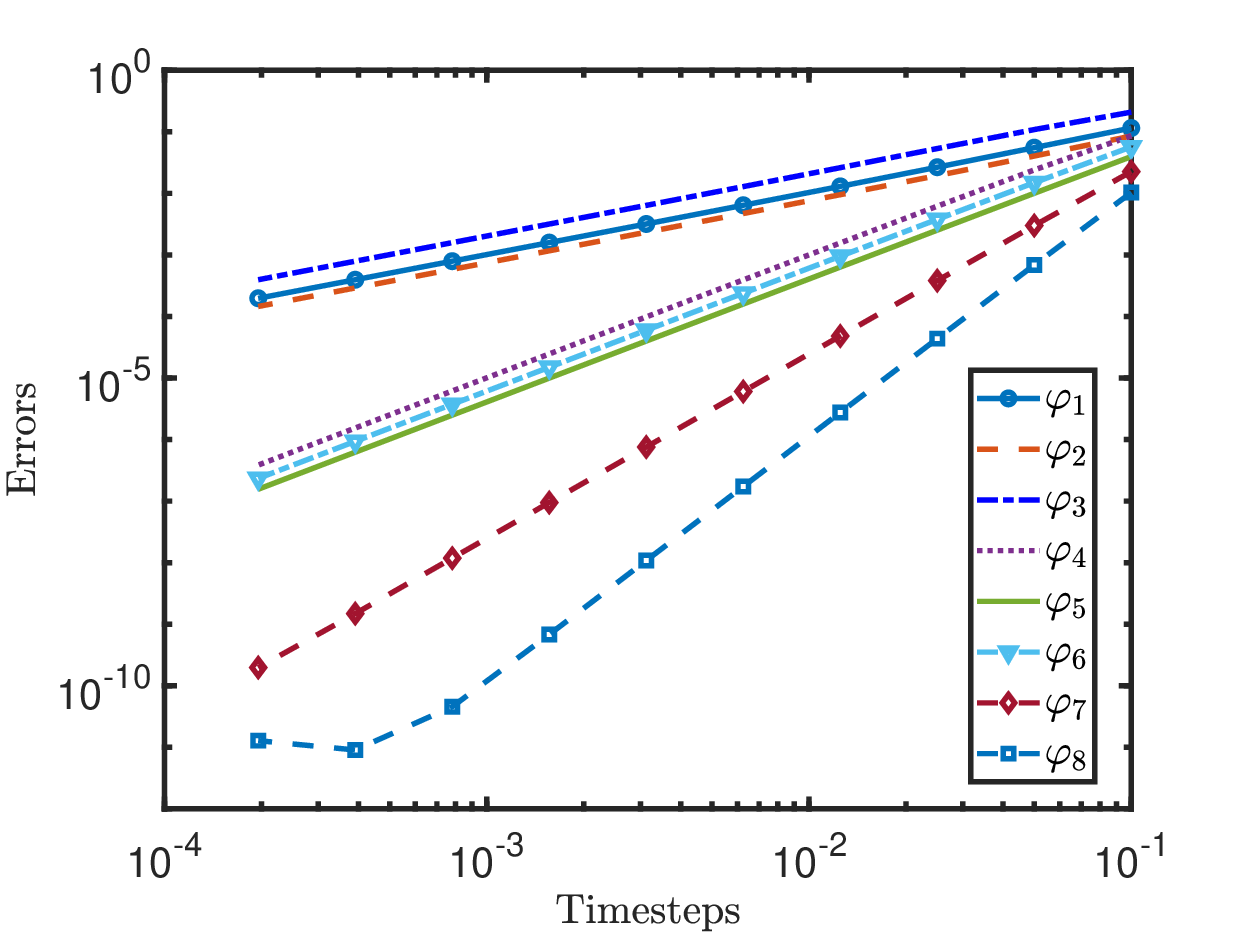}
    \caption{The errors for different choices of function $\varphi$ for $\Pi=0$.}
    \label{fig:seir_phiconv}
\end{figure}

\begin{figure}[!ht]
    \centering
    \includegraphics[width=0.7\linewidth]{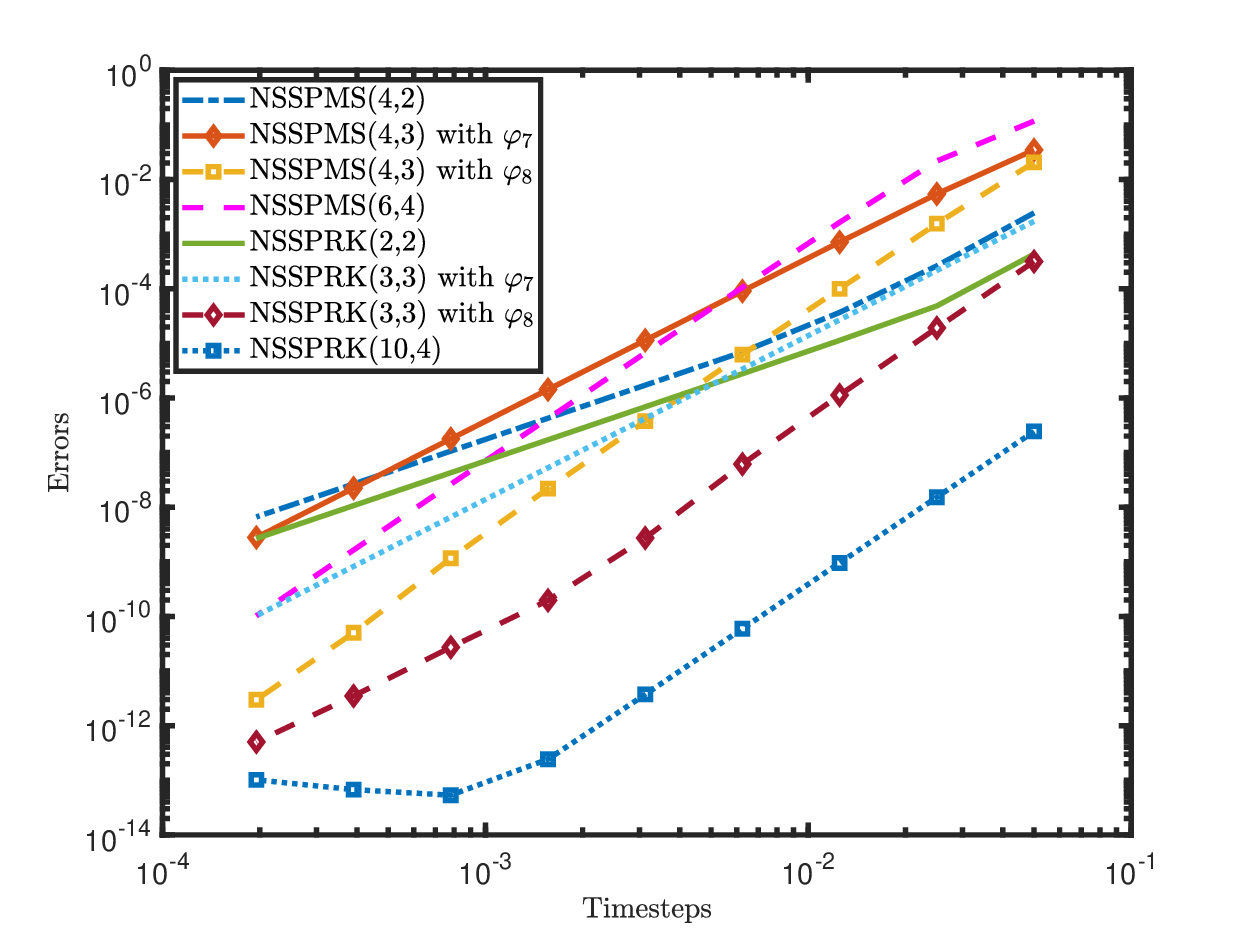}
    \caption{The errors of the different nonstandard methods for different values of $\Delta t$ for $\Pi=0$.}
    \label{fig:seir_meth}
\end{figure}

Now, we \rev{consider} the convergence of the multistep methods themselves: we run the tests with the same parameters, reference solution (with timestep $\Delta t = 0.05 \cdot 2 ^{-9} \cdot 10^{-3}$), value of $T$ ($T=1$), and initial values as before, but fix the choice $\varphi_8$ and compare the errors of NSSPMS(4,2), NSSPMS(4,3) and NSSPMS(6,4), \rev{along with the nonstandard Runge-Kutta methods NSSPRK(2,2), NSSPRK(3,3) and NSSPRK(10,4)} with timesteps $\Delta t=\rev{0.05/2^k, \; k=0, 1, \dots, 9}$. For the initial values \rev{of the multistep methods}, we used NSSPRK(2,2) with $\varphi_5$, NSSPRK(3,3) with $\varphi_7$, and NSSPRK(10,4) with $\varphi_8$, which are of order $2$, $3$, and $4$, respectively. The errors are plotted in Figure \ref{fig:seir_meth} and the values of the errors and orders are listed in Table \ref{tab:seir_method_errors} \rev{in Appendix B}. \rev{As it can be seen, unlike the case of the logistic equation, here both third order methods attain an order of $4$, even for smaller timesteps when we use $\varphi_8$ (although for the Runge-Kutta methods, this effect is less significant for smaller timesteps) - however, upon choosing $\varphi_7$, the order reduces to three. Moreover, the errors of the Runge-Kutta methods are smaller than the ones of the multistep schemes.}

\begin{figure}[!ht]
    \centering
    \includegraphics[width=0.45\linewidth]{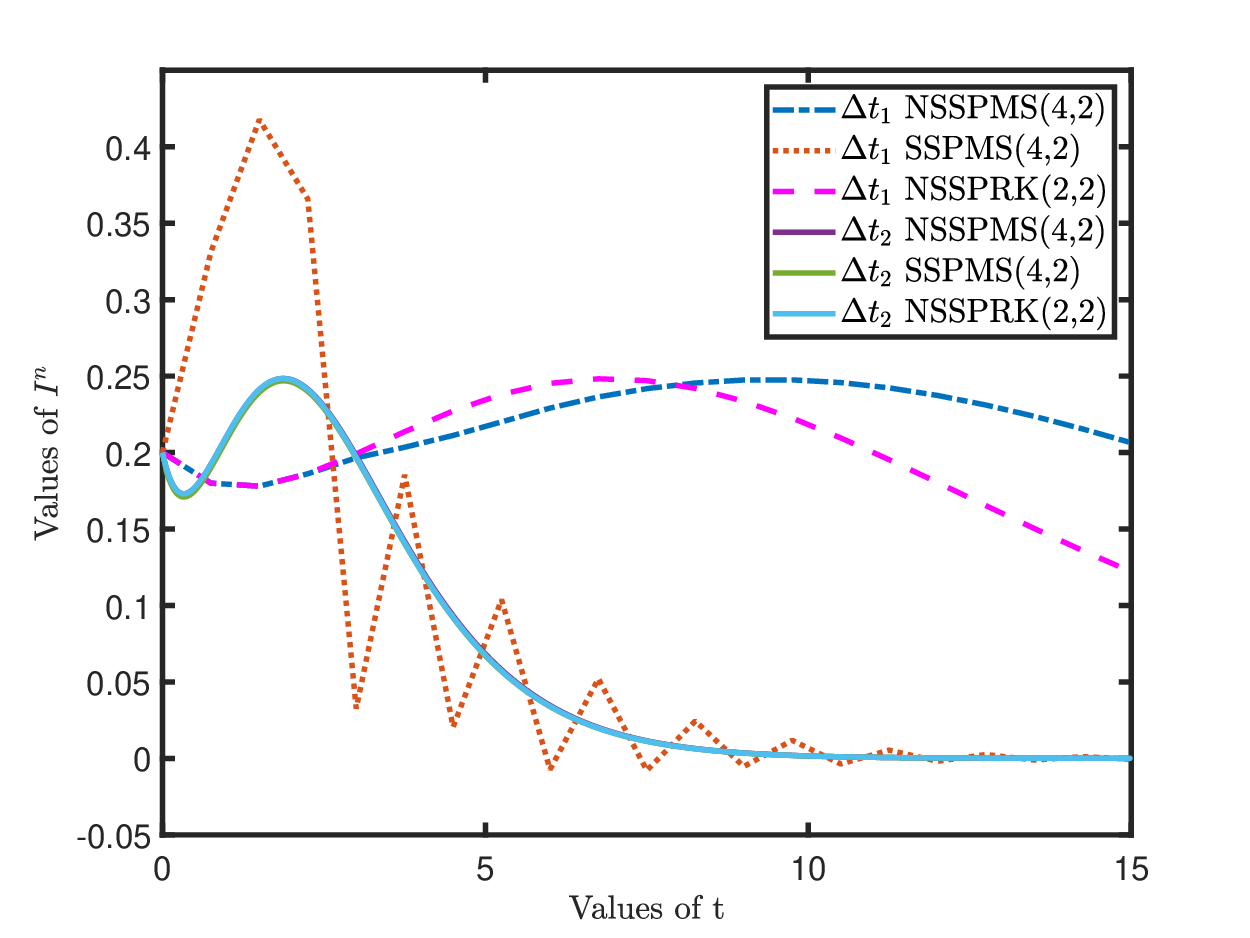}
    \includegraphics[width=0.45\linewidth]{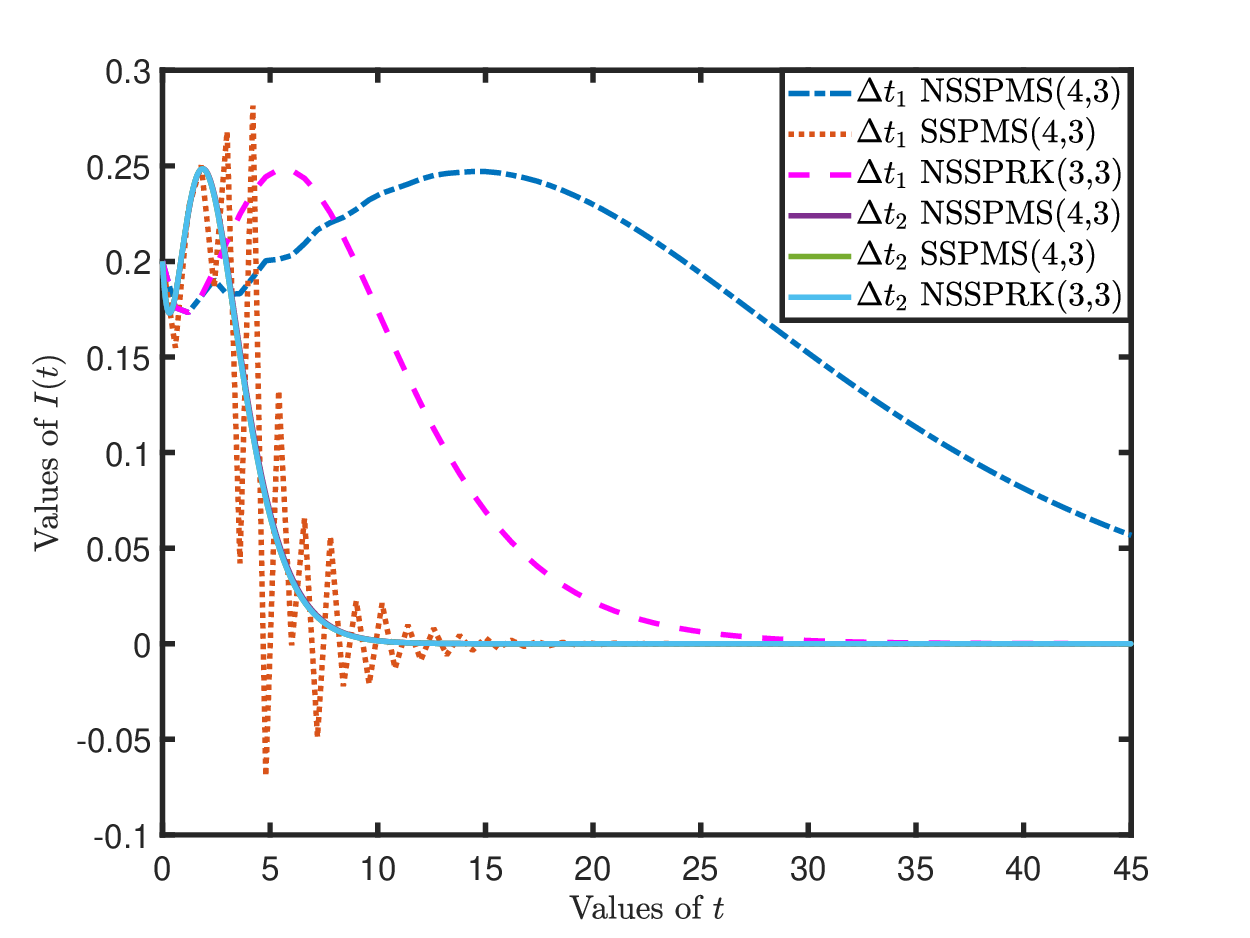}
    \includegraphics[width=0.45\linewidth]{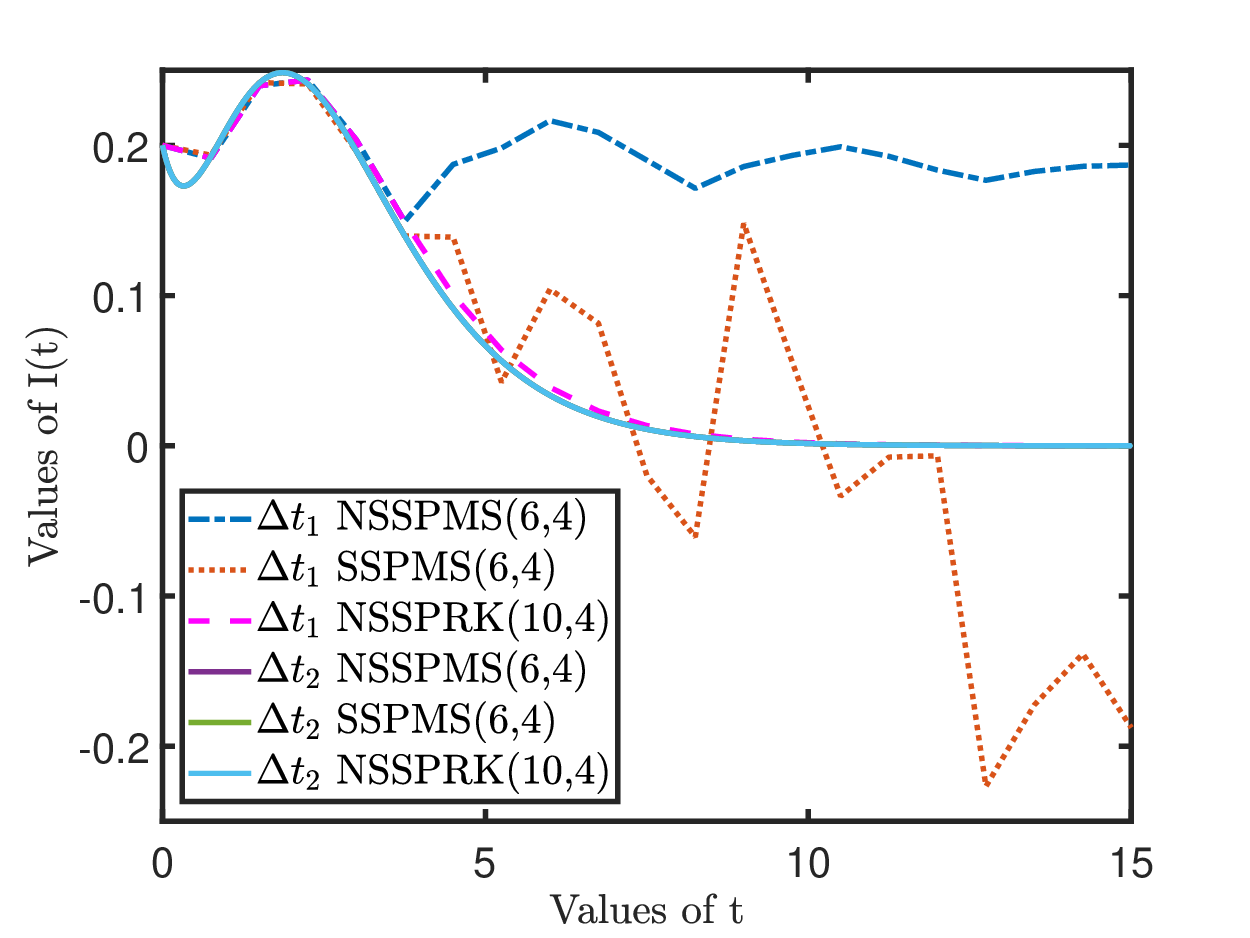}
    \caption{The comparison of the standard and nonstandard methods. Upper left: SSPMS(4,2), NSSPMS(4,2) \rev{and NSSPRK(2,2)} with $\Delta t_1=15/20$ and $\Delta t_2=0.02$. Upper right: SSPMS(4,3), NSSPMS(4,3) \rev{and NSSPRK(3,3)} with $\Delta t_1=15/25$ and $\Delta t_2=0.02$. Lower: SSPMS(6,4), NSSPMS(6,4) \rev{NSSPRK(10,4)} with $\Delta t_1=15/20$ and $\Delta t_2=0.02$. The solutions \rev{corresponding to time-step $\Delta t_2$} are very close to each other in all three cases.}
    \label{fig:seir_consmethods}
\end{figure}

For the \rev{next} test concerning the \rev{conservative} SEIR model \rev{($\Pi=0$)}, we compare the behavior of the nonstandard methods and the standard counterparts for two different timesteps (the unmentioned parameters remain the same) for \rev{$t \in [0,15]$ (in the cases of first order and third order), and $t \in [0,45]$ in the case of second order methods}. For the initial values of the multistep methods, Runge-Kutta methods with an appropriate order are used: standard ones for the standard multistep methods and nonstandard ones (with appropriate choices of $\varphi$) for the nonstandard versions. \rev{We also plot the graphs produced by the nonstandard Runge-Kutta methods.} In all three cases, we only plot the numerical solutions of $I^n$, but the others also show signs of similar behavior. In Figure \ref{fig:seir_consmethods}, we can see that the standard methods behave in an unstable way: they become negative, while the nonstandard ones remain in the upper half plane for the whole process, even for larger timesteps. Also, for a small timestep, the solutions are very close to each other and \rev{all three} of them behave as expected. \rev{It should also be mentioned that the errors are significantly smaller in the case of the Runge-Kutta methods, especially in the fourth order case.}

\begin{figure}[!ht]
    \centering
    \includegraphics[width=0.45\linewidth]{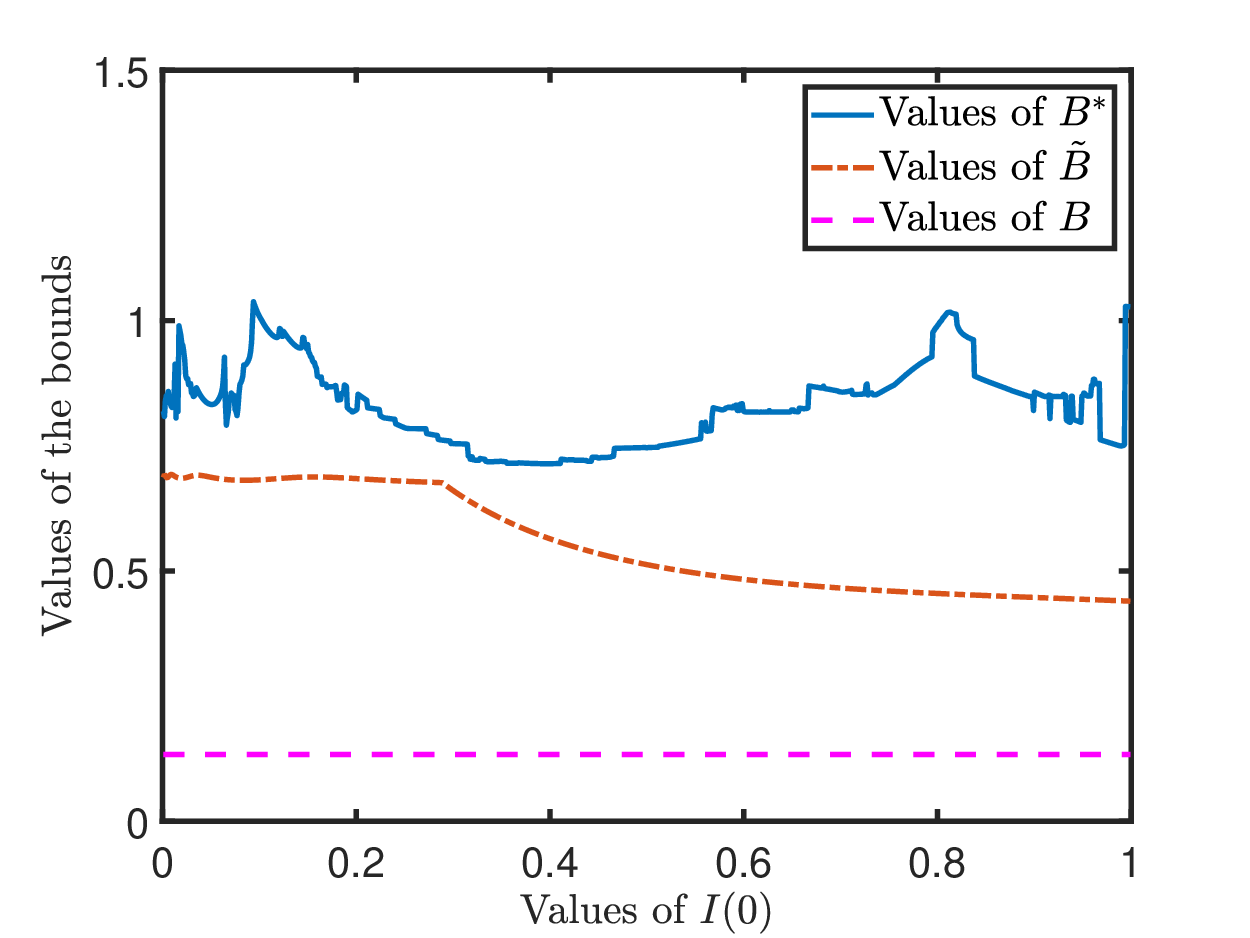}
    \includegraphics[width=0.45\linewidth]{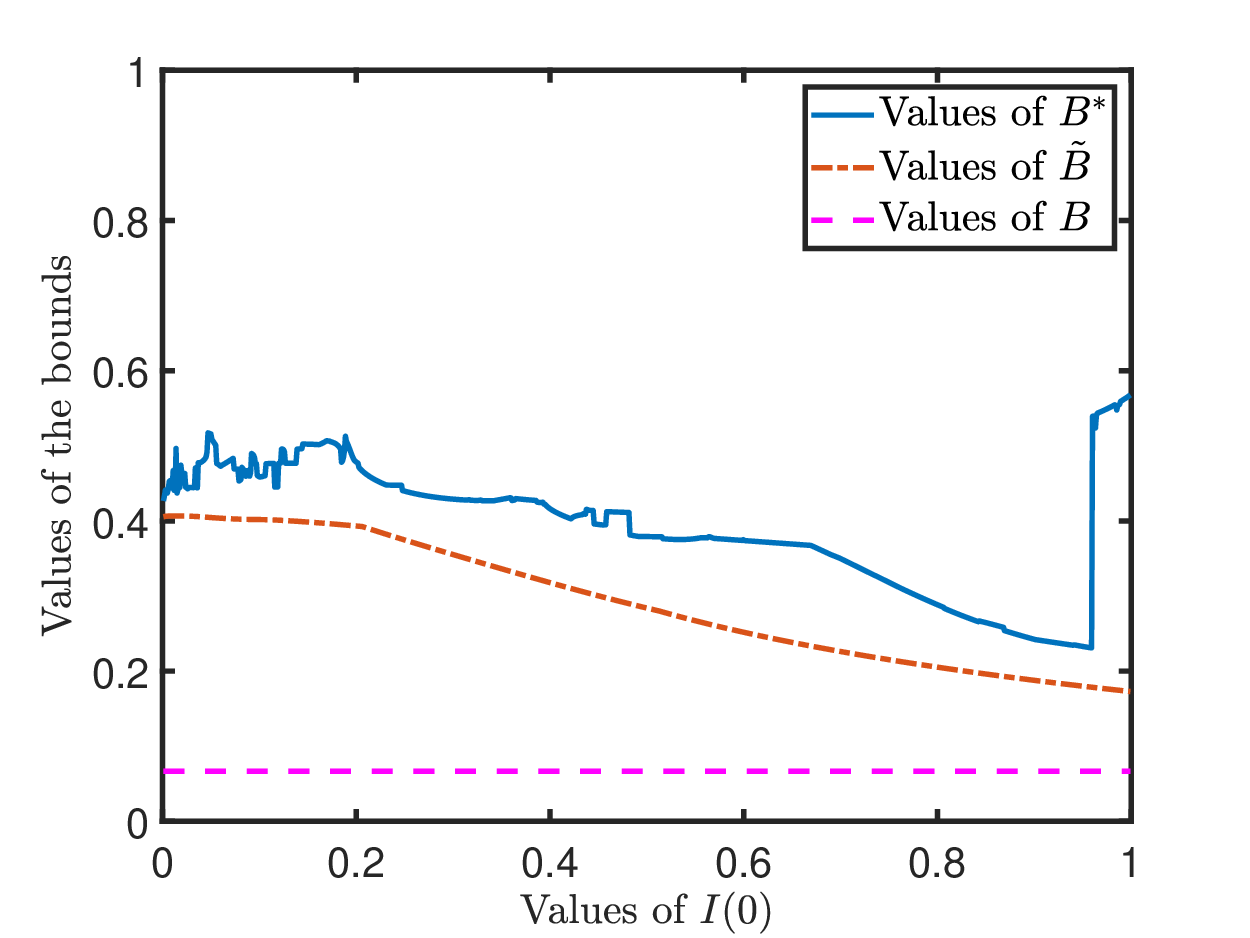}
    \includegraphics[width=0.45\linewidth]{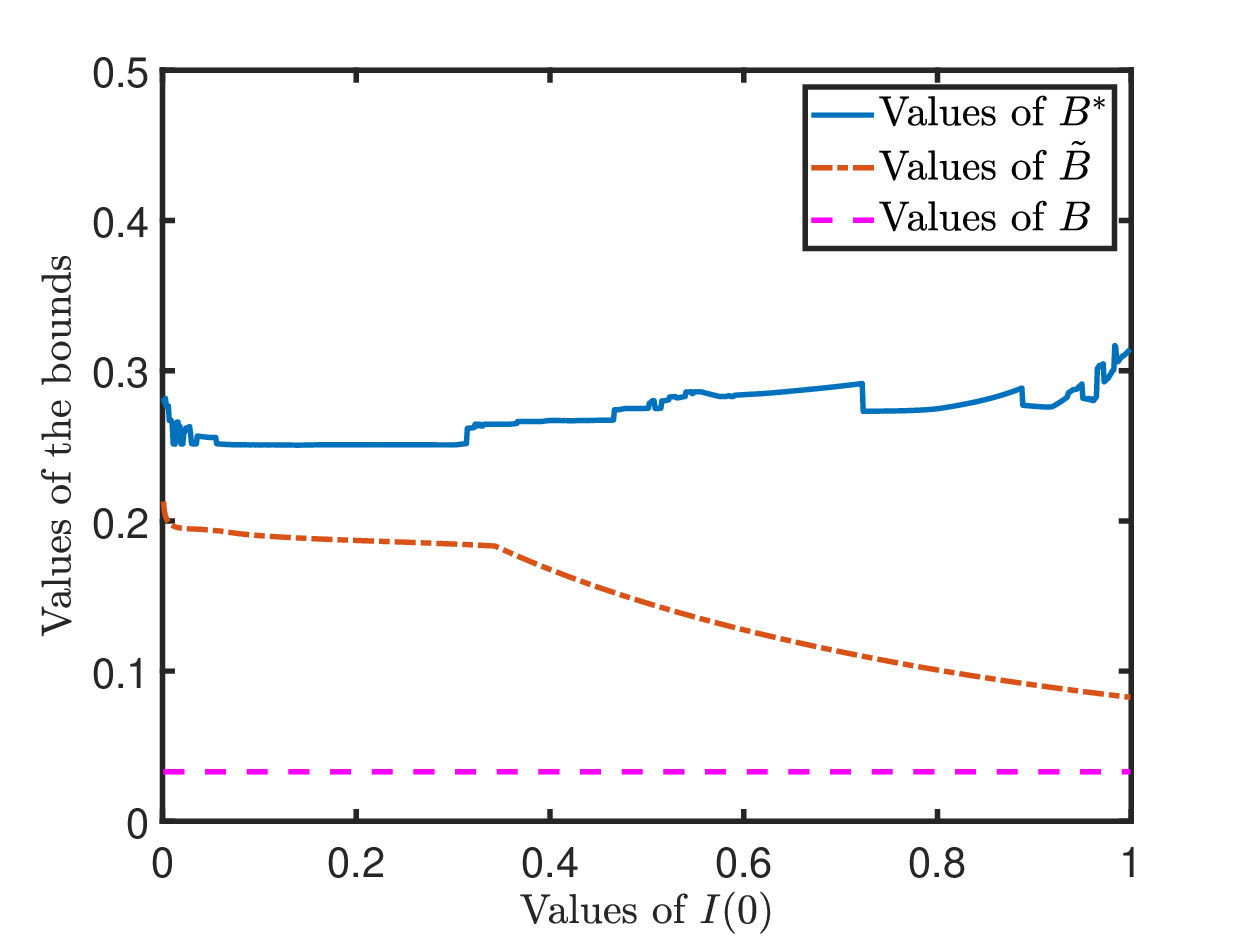}
    \caption{\rev{The sufficient bound $\mathcal{B}$ and the ``real'' bounds for methods NSSPMS(4,2) (upper left), NSSPMS(4,3) (upper right) and NSSPMS(6,4) (lower) for the boundedness property ($\mathcal{B}^*$) and the weak monotonicity ($\widetilde{\mathcal{B}}$) with $\Pi=0$.}}
    \label{fig:seir_boundsharp}
\end{figure}

\rev{Lastly, the sharpness of the bound for function $\varphi$ is tested in these cases too. For this, we used the initial condition $(S^0, E^0, I^0, R^0) = (1-I^0, 0, I^0, 0)$, and the methods were observed for $1000$ different values of $I^0$ in the interval $[0.001, 0.999]$. For a given value of $I^0$, the method is considered to behave as expected if the given property is preserved for $1000$ different timesteps from the interval $\Delta t \in [0.5, 3]$ on the time interval $t \in [0, 100]$. By using the bisection method, we determine the bound $\mathcal{B}^*$ for the boundedness of the solution (meaning that $S^n, I^n, E^n, R^n \in [0,1])$ and bound $\widetilde{\mathcal{B}}$ for the weak monotonicity property \eqref{eq:weakmon}. The bounds are plotted in Figure \ref{fig:seir_boundsharp}. As it can be seen, considerably larger bounds for function $\varphi$ could be used - however, one might be able to find such initial conditions where the sufficient bound proved in Proposition \ref{prop:seir_bound} is closer to the ``real'' one.}

\rev{\subsubsection{\texorpdfstring{The case of nonzero $\Pi$}{The case of nonzero Pi}}}

\rev{In this section, we consider the case when the SEIR model has an influx of healthy individuals, i.e. $\Pi >0$ in equation \eqref{sir}. It turns out that the order of method NSSPMS(6,4) upon using different $\varphi$ functions is very similar to the one we got in the case $\Pi=0$ (Figure \ref{fig:seir_phiconv} and Table \ref{tab:seir_phi_errors}). Moreover, the orders of the different nonstandard multistep and Runge-Kutta methods are also similar: the third order NSSPMS(4,3) and NSSPRK(3,3) (when used with $\varphi_8$) attain an order of four in these cases too, although numerical experiments indicate that the order of these methods decrease as $\Pi$ is increased, and gets close to $3$ when $\Pi$ is large (e.g. around $\Pi=100$).}

\begin{figure}[!ht]
    \centering
    \includegraphics[width=0.45\linewidth]{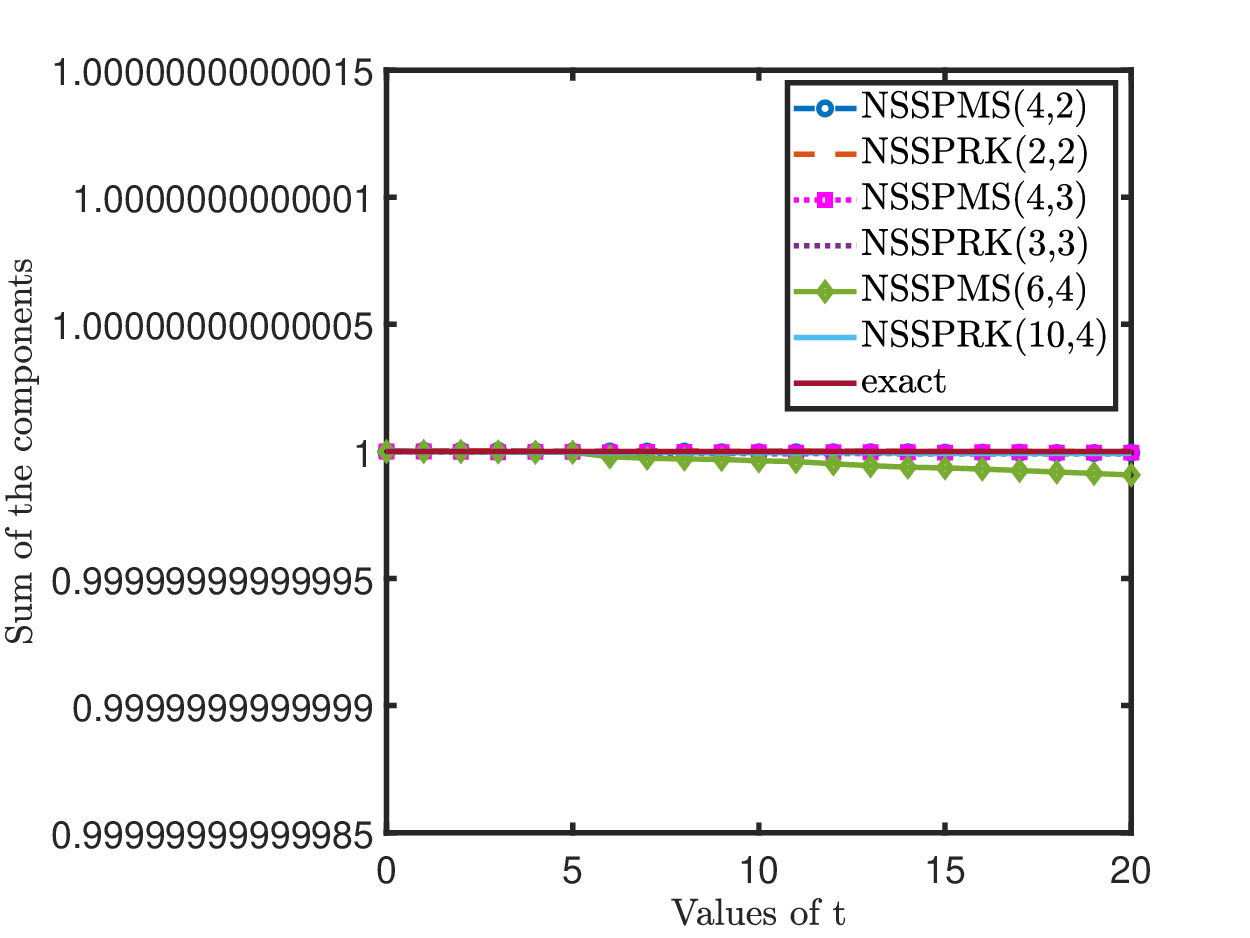}
    \includegraphics[width=0.45\linewidth]{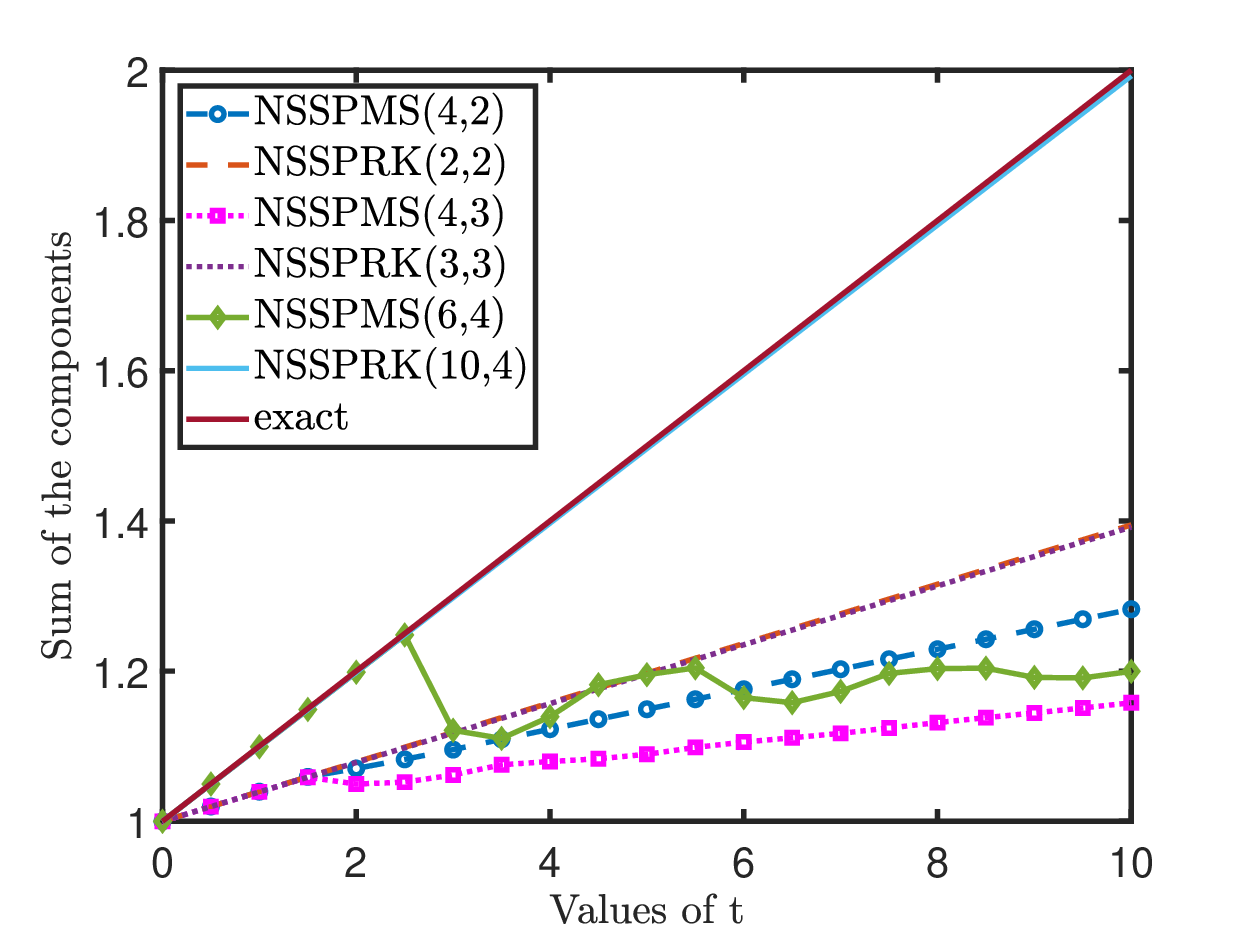}
    \includegraphics[width=0.45\linewidth]{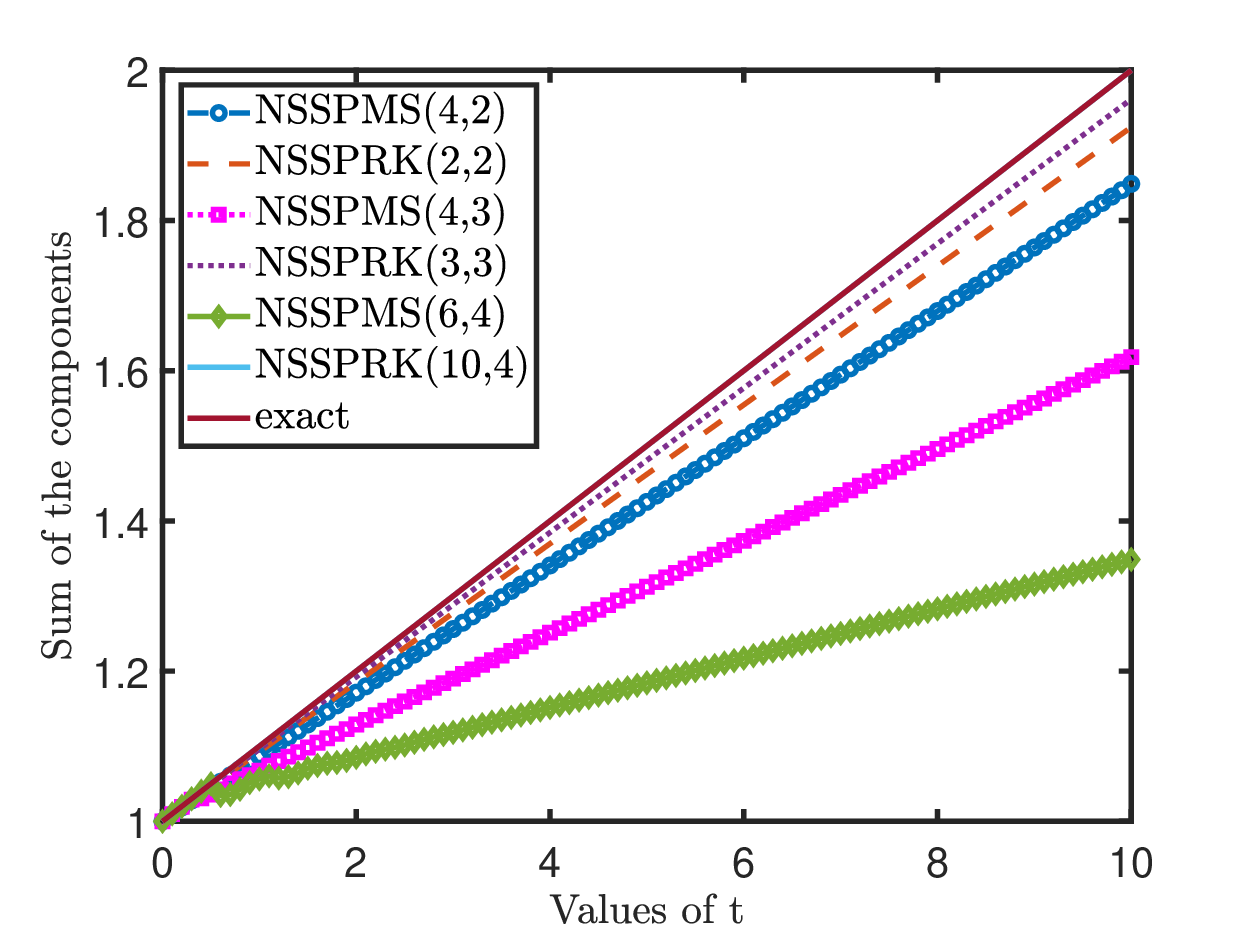}
    \includegraphics[width=0.45\linewidth]{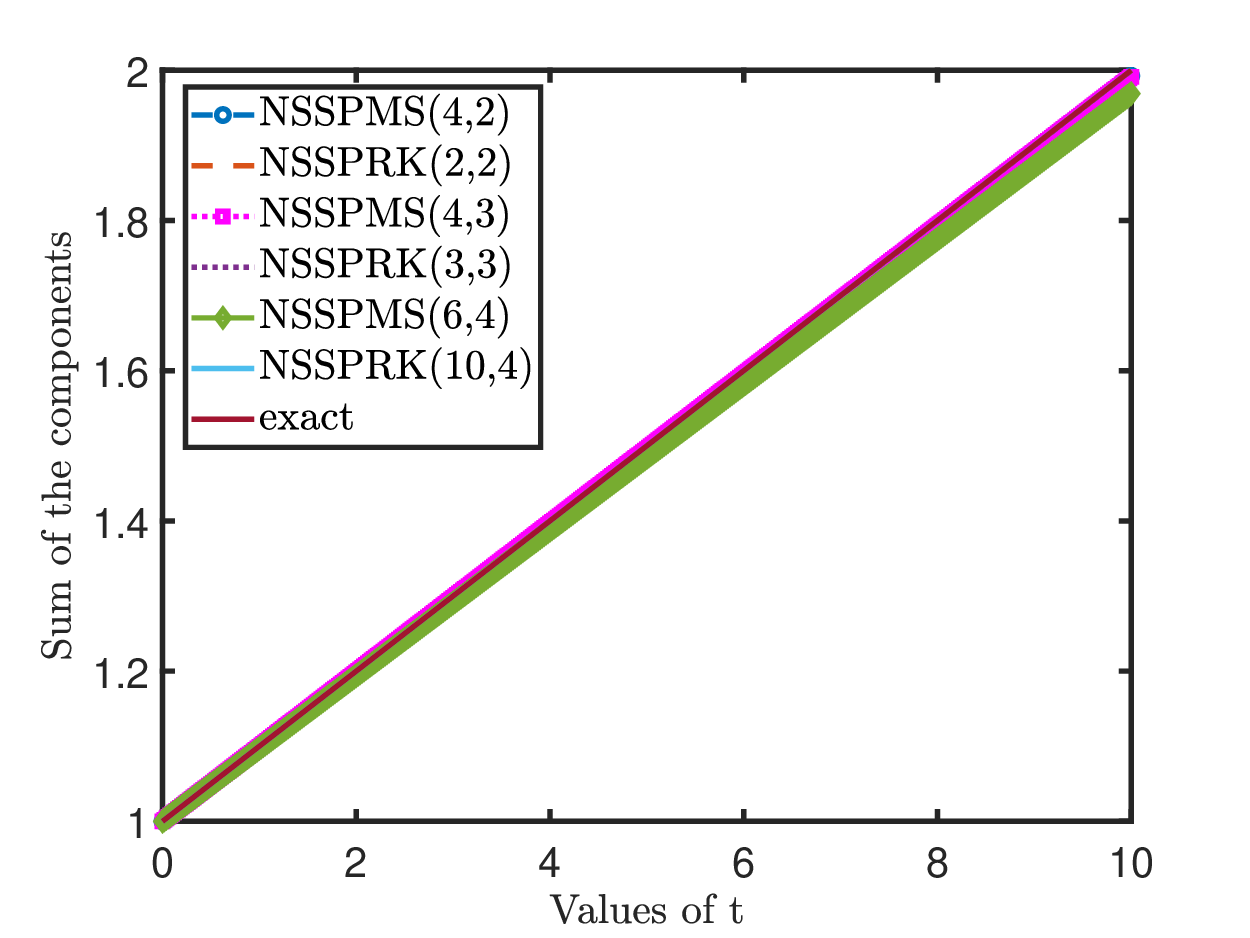}
    \caption{\rev{The sum $S^n+E^n+I^n+R^n$ for $\Pi=0$ with $\Delta t=1$ (upper left), and for $\Pi=0.1$ with $\Delta t=0.5$ (upper right), $\Delta t=0.1$ (lower left) and $\Delta t=0.02$ (lower right) with different nonstandard methods. Method NSSPRK(10,4) is always very close to the exact solution, therefore not visible on the graphs.}}
    \label{fig:seir_sum}
\end{figure}

\rev{Because of the similarities with the previous section, the main question we would like to focus on in the case of nonzero $\Pi$ is the behavior of the sum $S^n + E^n + I^n + R^n$. It is easy to see that from \eqref{sir}, for the function $N(t)=S(t)+E(t)+I(t)+R(t)$ (which can be thought of as the total population), we get the equation
\begin{equation*}
    N'(t) = \Pi,
\end{equation*}
which has a solution $N(t)=\Pi t + N(0)$ (again, $t_0=0$). Here we observe how well the nonstandard multistep (and Runge-Kutta) methods approximate this function when applied to equation \eqref{sir}. In Figure \ref{fig:seir_sum}, the sum $S^n + E^n + I^n + R^n$ is plotted for different nonstandard methods in the cases $\Pi=0$ and $\Pi=0.1$ (by using different timesteps) with initial condition $(S^0, E^0, I^0, R^0)=(0.8, 0, 0.2, 0)$. As it can be seen, the methods preserve the sum (up to machine precision) in the case of $\Pi=0$ even with a bigger timestep ($\Delta t=1$). For the case $\Pi=0.1$, the method NSSPRK(10,4) is always close to the exact value, while the others are furher, but as we decrease the timestep, they get closer and closer to the exact value. The nonstandard Runge-Kutta methods perform better in all of the cases, while the fourth-order multistep method NSSPMS(6,4) is the furthest from the exact value. }

\newpage
\section{Conclusions and discussion} \label{sec:conc}
In this work, we considered nonstandard multistep methods. By using a nonstandard version of the classical Taylor series (Lemma \ref{lem:taylor}), we were able to show that (under some conditions for the denominator function $\varphi$) the nonstandard versions of these methods have the same order as their standard counterparts. Note that Lemma \ref{lem:taylor} might also be useful in the \rev{analysis} of other, high-order nonstandard methods that are constructed in a similar fashion, i.e., involving a generalized denominator function.

In Section \ref{sec:pres}, it was shown that (under some further conditions for function $\varphi$) the nonstandard methods preserve the boundedness \rev{properties} (and as a special case, positivity) of the original, continuous problem. To achieve this goal, a bound for the timestep should be calculated under which the Euler method preserves these properties. We note that the elementary stable property was not considered, but theorems similar to those in \cite{dang} might be proven in the future for nonstandard multistep methods too. These results would involve the \rev{analysis} of the elementary stability property of standard multistep methods, which might turn out to be harder (in general, these are not even true \cite{anguelov}, as mentioned in Section \ref{sec:prelim}) -  this could also be a topic for further research. On the other hand, the \rev{analysis} of forward Euler methods might turn out to be easier, and for some systems, the preservation of positivity ensures elementary stability (see, for example, the case of an integro-differential equation \cite{takacs22} and its delayed version \cite{takacs222}).

In Section \ref{sec:phi} we also presented several possible choices for the function $\varphi$, for which the conditions of the theorems of the previous sections are met. We also mentioned there that in \cite{dang}, the authors gave a different construction of similar functions - in a further work, one could also compare the performance of the functions of this article and those defined in \cite{dang}.

In Section \ref{sec:Num}, we performed some tests for two test problems and observed that the nonstandard linear multistep methods behave as expected. In further articles, the same methods could also be applied to more robust ordinary systems and also to partial differential equations. One could also compare the performance of these multistep methods with other, higher-order methods (for example, the ones mentioned in Section \ref{sec:intro}). \rev{In the present work the nonstandard multistep methods were compared with nonstandard Runge-Kutta methods \cite{dang}. It turned out that although nonstandard multistep methods are able to achieve the same order as the nonstandard Runge-Kutta ones, the errors of the former are usually higher.}

Another interesting fact was the high order of convergence in the case of the method \linebreak NSSPMS(4,3). \rev{It turned out that upon applying a denominator function $\varphi$ which might be able achieve an order of four, these methods (which are usually only of third order) were also able to achieve this higher order of convergence. Because of this, it is usually advised to use such denominator functions that enable such high orders.}

\rev{As it was mentioned in Section \ref{sec:prelim}, from the usual two guiding principles of nonstandard methods, only one of them (the use of denominator functions) was adopted, while the implementation of nonlocal approximations was omitted. In a future work it would be interesting to incorporate the latter ideas too, which might result in smaller errors, while still keeping the high orders of the methods discussed in this paper.}

\bigskip

\noindent\textbf{Acknowledgements} 
The author would like to thank Prof. István Faragó for his insightful comments during the preparation of the article.

\bigskip

\noindent\textbf{Author Contributions} 
The author is the sole contributor to this paper.

\bigskip

\noindent\textbf{Funding:} This research has been supported by the National Research, Development and Innovation Office – NKFIH, grant no. K137699.
The research reported in this paper is also part of project no. BME-NVA-02, implemented with the support provided by the Ministry of Innovation and Technology of Hungary from the National Research, Development and Innovation Fund, financed under the TKP2021 funding scheme.

\bigskip

\noindent \textbf{Data Availability} The data used to support the findings of this study are available from the corresponding
author upon request.


\section*{Declarations}
\textbf{Conflict of interest:} The authors declare no competing interests.


\section*{Appendix A. The coefficients of the different methods}\label{secA1}

Here we list the coefficients of the different strongly-stable multistep and Runge-Kutta methods used in the paper, along with their SSP coefficients.
\begin{itemize}
    \item Explicit SSPMS(4,2) \cite{shu}, with $\mathcal{C}=\frac{2}{3}$:
    $$ \tilde{\alpha}_1=\dfrac{8}{9},\quad \tilde{\beta}_1=\dfrac{4}{3}, \quad \tilde{\alpha}_4=\dfrac{1}{9}. $$
    \item Explicit SSPMS(4,3) \cite{shu}, with $\mathcal{C}=\frac{1}{3}$:
    $$ \tilde{\alpha}_1=\dfrac{16}{27},\quad \tilde{\beta}_1=\dfrac{16}{81}, \quad \tilde{\alpha}_4=\dfrac{11}{27}, \quad \tilde{\beta}_4=\dfrac{4}{9}. $$
    \item Explicit SSPMS(6,4) \cite{ketcheson1,ketcheson2}, with $\mathcal{C}=0.1648$:
    $$ \tilde{\alpha}_1=0.342460855717007,\quad \tilde{\beta}_1=2.078553105578060,$$
    $$ \tilde{\alpha}_4=0.191798259434736,\quad \tilde{\beta}_4=1.164112222279710, $$
    $$ \tilde{\alpha}_5=0.093562124939008,\quad \tilde{\beta}_5=0.567871749748709, $$
    $$ \tilde{\alpha}_6=0.372178759909247. $$
    \item Explicit SSPRK(2,2) with $\mathcal{C}=1$ \cite{gottlieb98}:
    \begin{align*}
        u^{(1)} &= u^n + \Delta t f(u^n),\\
        u^{n+1} &= \dfrac{1}{2} u^n + \dfrac{1}{2} \left( u^{(1)} + \Delta t f(u^{(1)}) \right).
    \end{align*}
    \item Explicit SSPRK(3,3) with $\mathcal{C}=1$ \cite{sspbook}:
    \begin{align*}
        u^{(1)} &= u^n + \Delta t f(u^n),\\
        u^{(2)} &= \dfrac{3}{4} u^n + \dfrac{1}{4} u^{(1)} + \dfrac{1}{4} \Delta t f(u^{(1)}),\\
        u^{n+1} &= \dfrac{1}{3} u^n + \dfrac{2}{3} \left( u^{(2)} + \Delta t f(u^{(2)}) \right).
    \end{align*}
    \item Explicit SSPRK(10,4) with $\mathcal{C}=6$ \cite{ketcheson2}:
    \begin{align*}
        u^{(1)} &= u^n + \dfrac{1}{6}\Delta t f(u^n),\\
        u^{(j+1)} &= u^{(j)} + \dfrac{1}{6}\Delta t f(u^{(j)}), \qquad j=1, 2, 3, \\ 
        u^{(5)} &= \dfrac{3}{5} u^n + \dfrac{2}{5}u^{(4)} + \dfrac{1}{15} \Delta t f(u^{(4)}),\\
        u^{(j+1)} &= u^{(j)} + \dfrac{1}{6}\Delta t f(u^{(j)}), \qquad j=5,\dots, 8, \\ 
        u^{n+1} &= \dfrac{1}{25} u^n + \dfrac{9}{25} u^{(4)} + \dfrac{3}{5} u^{(9)} + \dfrac{3}{50} \Delta t f(u^{(4)}) + \dfrac{1}{10} \Delta t f(u^{(9)}) .
    \end{align*}
\end{itemize}

\rev{\section*{Appendix B. Tables of errors and orders}\label{sec:A2}
On the next pages, we list the errors and orders of the different methods discussed before.}

\begin{table}
    \ssmall
    \centering
    \rev{\begin{tabular}{c||c|c||c|c||c|c}
        $\Delta t$ & $\varphi_1$ errors & $\varphi_1$ orders & $\varphi_2$ errors & $\varphi_2$ orders & $\varphi_3$ errors & $\varphi_3$ orders  \\ \hline
        $0.1$ & $1.4009 \cdot 10^{-1}$ & - & $1.1611 \cdot 10^{-1}$ & -& $1.9461 \cdot 10^{-1}$ & - \\ \hline
        $0.1 \cdot 2^{-1}$ & $1.0611 \cdot 10^{-1}$ & $0.4008$ & $8.2200 \cdot 10^{-2}$ & $0.4983$ & $1.7290 \cdot 10^{-1}$ & $0.1706$ \\ \hline
        $0.1 \cdot 2^{-2}$ & $5.8780 \cdot 10^{-2}$ & $0.8521$ & $4.44178 \cdot 10^{-2}$ & $0.8958$ & $1.0622 \cdot 10^{-1}$ & $0.7028$ \\ \hline
        $0.1 \cdot 2^{-3}$  & $3.0750 \cdot 10^{-2}$ & $0.9347$ & $2.2833 \cdot 10^{-2}$ & $0.9522$ & $5.8599 \cdot 10^{-2}$ & $0.8582$ \\ \hline
        $0.1 \cdot 2^{-4}$ & $1.5669 \cdot 10^{-2}$ & $0.9727$ & $1.1576 \cdot 10^{-2}$ & $0.9799$ & $3.0621 \cdot 10^{-2}$ & $0.9364$ \\ \hline
        $0.1 \cdot 2^{-5}$ & $7.9013 \cdot 10^{-3}$ & $0.9877$ & $5.8249 \cdot 10^{-3}$ & $0.9909$ & $1.5626 \cdot 10^{-2}$ & $0.9706$ \\ \hline
        $0.1 \cdot 2^{-6}$ & $3.9666 \cdot 10^{-3}$ & $0.9942$ & $2.9212 \cdot 10^{-3}$ & $0.9957$ & $7.8894 \cdot 10^{-3}$ & $0.9860$ \\ \hline
        $0.1 \cdot 2^{-7}$ & $1.9871 \cdot 10^{-3}$ & $0.9972$ & $1.4627 \cdot 10^{-3}$ & $0.9979$ & $3.9634 \cdot 10^{-3}$ & $0.9932$ \\ \hline
        $0.1 \cdot 2^{-8}$ & $9.9452 \cdot 10^{-4}$ & $0.9986$ & $7.3190 \cdot 10^{-4}$ & $0.9990$ & $1.9863 \cdot 10^{-3}$ & $0.9966$ \\ \hline
        $0.1 \cdot 2^{-9}$ & $4.9750 \cdot 10^{-4}$ & $0.9993$ & $3.6608 \cdot 10^{-4}$ & $0.9995$ & $9.9431 \cdot 10^{-4}$ & $0.9983$ \\ \hline \hline 
    \end{tabular}
    \begin{tabular}{c||c|c||c|c||c|c}
        $\Delta t$ & $\varphi_4$ errors & $\varphi_4$ orders & $\varphi_5$ errors & $\varphi_5$ orders & $\varphi_6$ errors & $\varphi_6$ orders  \\ \hline
        $0.1$ & $1.4359 \cdot 10^{-1}$ & - & $9.7188 \cdot 10^{-2}$ & -& $1.1764 \cdot 10^{-1}$ & - \\ \hline
        $0.1 \cdot 2^{-1}$ & $8.2705 \cdot 10^{-2}$ & $0.7959$ & $4.1449 \cdot 10^{-2}$ & $1.2294$ & $5.7578 \cdot 10^{-2}$ & $1.0308$ \\ \hline
        $0.1 \cdot 2^{-2}$ & $2.7204 \cdot 10^{-2}$ & $1.6041$ & $1.1739 \cdot 10^{-2}$ & $1.8200$ & $1.7248 \cdot 10^{-2}$ & $1.7390$ \\ \hline
        $0.1 \cdot 2^{-3}$ & $7.5017 \cdot 10^{-3}$ & $1.8586$ & $3.0902 \cdot 10^{-3}$ & $1.9255$ & $4.6113 \cdot 10^{-3}$ & $1.9032$ \\ \hline
        $0.1 \cdot 2^{-4}$ & $1.9405 \cdot 10^{-3}$ & $1.9508$ & $7.8967 \cdot 10^{-4}$ & $1.9684$ & $1.1830 \cdot 10^{-3}$ & $1.9628$ \\ \hline
        $0.1 \cdot 2^{-5}$ & $4.9156 \cdot 10^{-4}$ & $1.9810$ & $1.994 \cdot 10^{-4}$ & $1.9854$ & $2.9904 \cdot 10^{-4}$ & $1.9840$ \\ \hline
        $0.1 \cdot 2^{-6}$ & $1.2358 \cdot 10^{-4}$ & $1.9919$ & $5.0099 \cdot 10^{-5}$ & $1.9930$ & $7.5143 \cdot 10^{-5}$ & $1.9926$ \\ \hline
        $0.1 \cdot 2^{-7}$ & $3.0976 \cdot 10^{-5}$ & $1.9963$ & $1.2555 \cdot 10^{-5}$ & $1.9966$ & $1.8832 \cdot 10^{-5}$ & $1.9965$ \\ \hline
        $0.1 \cdot 2^{-8}$ & $7.7534 \cdot 10^{-6}$ & $1.9982$ & $3.1424 \cdot 10^{-6}$ & $1.9983$ & $4.7135 \cdot 10^{-6}$ & $1.9983$ \\ \hline
        $0.1 \cdot 2^{-9}$ & $1.9395 \cdot 10^{-6}$ & $1.9991$ & $7.8606 \cdot 10^{-7}$ & $1.9992$ & $1.1791 \cdot 10^{-6}$ & $1.9991$ \\ \hline \hline
    \end{tabular}
    \begin{tabular}{c||c|c||c|c}
        $\Delta t$ & $\varphi_7$ errors & $\varphi_7$ orders & $\varphi_8$ errors & $\varphi_8$ orders  \\ \hline
        $0.1$ & $8.9836 \cdot 10^{-2}$ & - & $7.6103 \cdot 10^{-2}$ & - \\ \hline
        $0.1 \cdot 2^{-1}$ & $2.4513 \cdot 10^{-2}$ & $1.8738$ & $1.1542 \cdot 10^{-2}$ & $2.7211$ \\ \hline
        $0.1 \cdot 2^{-2}$ & $3.5736 \cdot 10^{-3}$ & $2.7781$ & $8.1974 \cdot 10^{-4}$ & $3.8155$ \\ \hline
        $0.1 \cdot 2^{-3}$ & $4.6978 \cdot 10^{-4}$ & $2.9273$ & $5.3510 \cdot 10^{-5}$ & $3.9373$ \\ \hline
        $0.1 \cdot 2^{-4}$ & $5.9937 \cdot 10^{-5}$ & $2.9705$ & $3.4099 \cdot 10^{-6}$ & $3.9720$ \\ \hline
        $0.1 \cdot 2^{-5}$ & $7.5646 \cdot 10^{-6}$ & $2.9861$ & $2.1515 \cdot 10^{-7}$ & $3.9863$ \\ \hline
        $0.1 \cdot 2^{-6}$ & $9.5005 \cdot 10^{-7}$ & $2.9932$ & $1.3511 \cdot 10^{-8}$ & $3.9932$ \\ \hline
        $0.1 \cdot 2^{-7}$ & $1.1904 \cdot 10^{-7}$ & $2.9966$ & $8.4697 \cdot 10^{-10}$ & $3.9956$ \\ \hline
        $0.1 \cdot 2^{-8}$ & $1.4898 \cdot 10^{-8}$ & $2.9982$ & $5.4143 \cdot 10^{-11}$ & $3.9675$ \\ \hline
        $0.1 \cdot 2^{-9}$ & $1.8655 \cdot 10^{-9}$ & $2.9974$ & $5.673 \cdot 10^{-12}$ & $3.2545$ \\ \hline
    \end{tabular}}
    \caption{\rev{The errors and the corresponding orders of method NSSPMS(6,4) using different $\varphi$ functions while solving the logistic equation with $c=2$, plotted in Figure \ref{fig:phi_orders}.}}
    \label{tab:phi_errors}
\end{table}

\begin{table}
    \ssmall
    \centering
    \rev{\begin{tabular}{c||c|c||c|c||c|c}
         & \multicolumn{2}{c||}{NSSPMS(4,2)} & \multicolumn{2}{c||}{NSSPMS(4,3) with $\varphi_8$} & \multicolumn{2}{c}{NSSPMS(4,3) with $\varphi_7$}  \\ \hline
        $\Delta t$ & errors & orders & errors & orders & errors & orders  \\ \hline
        $0.05$ & $1.6660 \cdot 10^{-4}$ & - & $8.2145 \cdot 10^{-4}$ & -& $3.4349 \cdot 10^{-3}$ & - \\ \hline
        $0.05 \cdot 2^{-1}$ & $6.0870 \cdot 10^{-5}$ & $1.4527$ & $5.7502 \cdot 10^{-5}$ & $3.8365$ & $4.5630 \cdot 10^{-4}$ & $2.9122$ \\ \hline
        $0.05 \cdot 2^{-2}$ & $1.7144 \cdot 10^{-5}$ & $1.8280$ & $4.1033 \cdot 10^{-6}$ & $3.8087$ & $5.8507 \cdot 10^{-5}$ & $2.9633$ \\ \hline
        $0.05 \cdot 2^{-3}$  & $4.4918 \cdot 10^{-6}$ & $1.9324$ & $3.1262 \cdot 10^{-7}$ & $3.7143$ & $7.4020 \cdot 10^{-6}$ & $2.9826$ \\ \hline
        $0.05 \cdot 2^{-4}$ & $1.1463 \cdot 10^{-6}$ & $1.9703$ & $2.6326 \cdot 10^{-8}$ & $3.5698$ & $9.3074 \cdot 10^{-7}$ & $2.9915$ \\ \hline
        $0.05 \cdot 2^{-5}$ & $2.8934 \cdot 10^{-7}$ & $1.9862$ & $2.4865 \cdot 10^{-9}$ & $3.4043$ & $1.1668 \cdot 10^{-7}$ & $2.9958$ \\ \hline
        $0.05 \cdot 2^{-6}$ & $7.2670 \cdot 10^{-8}$ & $1.9933$ & $2.6035 \cdot 10^{-10}$ & $3.2556$ & $1.4607 \cdot 10^{-8}$ & $2.9979$ \\ \hline
        $0.05 \cdot 2^{-7}$ & $1.8208 \cdot 10^{-8}$ & $1.9967$ & $2.9433 \cdot 10^{-11}$ & $3.1450$ & $1.8273 \cdot 10^{-9}$ & $2.9990$ \\ \hline
        $0.05 \cdot 2^{-8}$ & $4.5571 \cdot 10^{-9}$ & $1.9984$ & $3.5845 \cdot 10^{-12}$ & $3.0376$ & $2.2860 \cdot 10^{-10}$ & $2.9988$ \\ \hline \hline
    \end{tabular}}
    \rev{\begin{tabular}{c||c|c||c|c||c|c}
         & \multicolumn{2}{c||}{NSSPRK(2,2)} & \multicolumn{2}{c||}{NSSPRK(3,3) with $\varphi_8$} & \multicolumn{2}{c}{NSSPRK(3,3) with $\varphi_7$}  \\ \hline
        $\Delta t$ & errors & orders & errors & orders & errors & orders  \\ \hline
        $0.05$ & $3.2621 \cdot 10^{-4}$ & - & $2.0710 \cdot 10^{-6}$ & -& $1.4771 \cdot 10^{-5}$ & - \\ \hline
        $0.05 \cdot 2^{-1}$ & $7.7614 \cdot 10^{-5}$ & $2.0714$ & $2.8654 \cdot 10^{-7}$ & $2.8535$ & $1.8598 \cdot 10^{-6}$ & $2.9896$ \\ \hline
        $0.05 \cdot 2^{-2}$ & $1.9039 \cdot 10^{-5}$ & $2.0274$ & $3.7559 \cdot 10^{-8}$ & $2.9315$ & $2.3330 \cdot 10^{-7}$ & $2.9949$ \\ \hline
        $0.05 \cdot 2^{-3}$  & $4.7220 \cdot 10^{-6}$ & $2.0114$ & $4.8041 \cdot 10^{-9}$ & $2.9668$ & $2.9213 \cdot 10^{-8}$ & $2.9975$ \\ \hline
        $0.05 \cdot 2^{-4}$ & $1.1763 \cdot 10^{-6}$ & $2.0051$ & $6.0734 \cdot 10^{-10}$ & $2.9837$ & $3.6548 \cdot 10^{-9}$ & $2.9987$ \\ \hline
        $0.05 \cdot 2^{-5}$ & $2.9358 \cdot 10^{-7}$ & $2.0024$ & $7.6316 \cdot 10^{-11}$ & $2.9925$ & $4.5709 \cdot 10^{-10}$ & $2.9993$ \\ \hline
        $0.05 \cdot 2^{-6}$ & $7.3336 \cdot 10^{-8}$ & $2.0012$ & $9.5060 \cdot 10^{-12}$ & $3.0051$ & $5.7211 \cdot 10^{-11}$ & $2.9981$ \\ \hline
        $0.05 \cdot 2^{-7}$ & $1.8327 \cdot 10^{-8}$ & $2.0006$ & $1.0607 \cdot 10^{-12}$ & $3.1638$ & $7.2802 \cdot 10^{-12}$ & $2.9743$ \\ \hline
        $0.05 \cdot 2^{-8}$ & $4.5807 \cdot 10^{-9}$ & $2.0003$ & $1.2346 \cdot 10^{-13}$ & $3.1029$ & $1.1662 \cdot 10^{-12}$ & $2.6422$ \\ \hline
    \end{tabular}} 
   \rev{\begin{tabular}{c||c|c||c|c}
         & \multicolumn{2}{c||}{NSSPRK(6,4)} & \multicolumn{2}{c}{NSSPRK(10,4)}  \\ \hline
        $\Delta t$ & errors & orders & errors & orders \\ \hline
        $0.05$ & $1.1541 \cdot 10^{-2}$ & - & $8.9811 \cdot 10^{-9}$ & - \\ \hline
        $0.05 \cdot 2^{-1}$ & $8.1974 \cdot 10^{-4}$ & $3.8155$ & $5.5896 \cdot 10^{-10}$ & $4.0061$ \\ \hline
        $0.05 \cdot 2^{-2}$ & $5.3510 \cdot 10^{-6}$ & $3.9373$ & $3.4863 \cdot 10^{-11}$ & $4.0030$ \\ \hline
        $0.05 \cdot 2^{-3}$  & $3.4099 \cdot 10^{-6}$ & $3.9720$ & $2.1829 \cdot 10^{-12}$ & $3.9974$ \\ \hline
        $0.05 \cdot 2^{-4}$ & $2.1515 \cdot 10^{-7}$ & $3.9863$ & $1.4166 \cdot 10^{-13}$ & $3.9457$ \\ \hline
        $0.05 \cdot 2^{-5}$ & $1.3511 \cdot 10^{-8}$ & $3.9932$  & $2.6867 \cdot 10^{-14}$ & $2.3985$ \\ \hline
        $0.05 \cdot 2^{-6}$ & $8.4697 \cdot 10^{-10}$ & $3.9956$ & $4.0412 \cdot 10^{-14}$ & - \\ \hline
        $0.05 \cdot 2^{-7}$ & $5.4143 \cdot 10^{-11}$ & $3.9675$  & $7.4607 \cdot 10^{-14}$ & -\\ \hline
        $0.05 \cdot 2^{-8}$ & $5.673 \cdot 10^{-12}$ & $3.2545$ & $1.7186 \cdot 10^{-13}$ & - \\ \hline
    \end{tabular}}
    \caption{\rev{The errors and the corresponding orders of the different multistep methods while solving the logistic equation with $c=2$, plotted in Figure \ref{fig:LMM_orders}.}}
    \label{tab:LMM_errors}
\end{table}

\begin{table}
    \ssmall
    \centering
    \rev{\begin{tabular}{c||c|c||c|c||c|c}
        $\Delta t$ & $\varphi_1$ errors & $\varphi_1$ orders & $\varphi_2$ errors & $\varphi_2$ orders & $\varphi_3$ errors & $\varphi_3$ orders  \\ \hline
        $2 \cdot 10^{-4}$ & $4.6308 \cdot 10^{1}$ & - & $3.7984 \cdot 10^{1}$ & -& $6.6018 \cdot 10^{1}$ & - \\ \hline
        $2 \cdot 10^{-4} \cdot 2^{-1}$ & $3.5903 \cdot 10^{1}$ & $0.3672$ & $2.7737 \cdot 10^{1}$ & $0.4536$ & $6.6018 \cdot 10^{1}$ & $0.1555$ \\ \hline
        $2 \cdot 10^{-4} \cdot 2^{-2}$ & $1.9406 \cdot 10^{1}$ & $0.8876$ & $1.4576 \cdot 10^{1}$ & $0.9283$ & $3.5249 \cdot 10^{1}$ & $0.7498$ \\ \hline
        $2 \cdot 10^{-4} \cdot 2^{-3}$  & $1.0111 \cdot 10^{1}$ & $0.9405$ & $7.5087 \cdot 10^{0}$ & $0.9570$ & $1.9283 \cdot 10^{1}$ & $0.8702$ \\ \hline
        $2 \cdot 10^{-4} \cdot 2^{-4}$ & $5.1489 \cdot 10^{0}$ & $0.9736$ & $3.8049 \cdot 10^{0}$ & $0.9807$ & $1.0059 \cdot 10^{0}$ & $0.9389$ \\ \hline
        $2 \cdot 10^{-4} \cdot 2^{-5}$ & $2.5969 \cdot 10^{0}$ & $0.9875$ & $1.9147 \cdot 10^{0}$ & $0.9907$ & $5.1337 \cdot 10^{0}$ & $0.9704$ \\ \hline
        $2 \cdot 10^{-4} \cdot 2^{-6}$ & $1.3040 \cdot 10^{0}$ & $0.9939$ & $9.6039 \cdot 10^{-1}$ & $0.9955$ & $2.5928 \cdot 10^{0}$ & $0.9855$ \\ \hline
        $2 \cdot 10^{-4} \cdot 2^{-7}$ & $6.5334 \cdot 10^{-1}$ & $0.9970$ & $4.8094 \cdot 10^{-1}$ & $0.9977$ & $1.3029 \cdot 10^{0}$ & $0.9928$ \\ \hline
        $2 \cdot 10^{-4} \cdot 2^{-8}$ & $3.2701 \cdot 10^{-1}$ & $0.9985$ & $2.4066 \cdot 10^{-1}$ & $0.9989$ & $6.5307 \cdot 10^{-1}$ & $0.9964$ \\ \hline
        $2 \cdot 10^{-4} \cdot 2^{-9}$ & $1.6359 \cdot 10^{-1}$ & $0.9993$ & $1.2038 \cdot 10^{-1}$ & $0.9994$ & $3.2694 \cdot 10^{-1}$ & $0.9982$ \\ \hline \hline 
    \end{tabular}
    \begin{tabular}{c||c|c||c|c||c|c}
        $\Delta t$ & $\varphi_4$ errors & $\varphi_4$ orders & $\varphi_5$ errors & $\varphi_5$ orders & $\varphi_6$ errors & $\varphi_6$ orders  \\ \hline
        $2 \cdot 10^{-4}$ & $4.7539 \cdot 10^{1}$ & - & $3.1548 \cdot 10^{1}$ & -& $3.8509 \cdot 10^{1}$ & - \\ \hline
        $2 \cdot 10^{-4} \cdot 2^{-1}$ & $2.7909 \cdot 10^{1}$ & $0.7684$ & $1.3964 \cdot 10^{1}$ & $1.1759$ & $1.9400 \cdot 10^{1}$ & $0.9891$ \\ \hline
        $2 \cdot 10^{-4} \cdot 2^{-2}$ & $8.9750 \cdot 10^{0}$ & $1.6368$ & $3.8755 \cdot 10^{0}$ & $1.8492$ & $5.6923 \cdot 10^{0}$ & $1.7690$ \\ \hline
        $2 \cdot 10^{-4} \cdot 2^{-3}$ & $2.4687 \cdot 10^{0}$ & $1.8622$ & $1.0173 \cdot 10^{0}$ & $1.9296$ & $1.5179 \cdot 10^{0}$ & $1.9070$ \\ \hline
        $2 \cdot 10^{-4} \cdot 2^{-4}$ & $6.3821 \cdot 10^{-1}$ & $1.9517$ & $2.5975 \cdot 10^{-1}$ & $1.9696$ & $3.8910 \cdot 10^{-1}$ & $1.9638$ \\ \hline
        $2 \cdot 10^{-4} \cdot 2^{-5}$ & $1.6165 \cdot 10^{-1}$ & $1.9812$ & $6.5583 \cdot 10^{-2}$ & $1.9857$ & $9.8342 \cdot 10^{-2}$ & $1.9843$ \\ \hline
        $2 \cdot 10^{-4} \cdot 2^{-6}$ & $4.0639 \cdot 10^{-2}$ & $1.9919$ & $1.6475 \cdot 10^{-2}$ & $1.9931$ & $2.4710 \cdot 10^{-2}$ & $1.9927$ \\ \hline
        $2 \cdot 10^{-4}\cdot 2^{-7}$ & $1.0186 \cdot 10^{-2}$ & $1.9963$ & $4.1284 \cdot 10^{-3}$ & $1.9966$ & $6.1925 \cdot 10^{-3}$ & $1.9965$ \\ \hline
        $2 \cdot 10^{-4} \cdot 2^{-8}$ & $2.5496 \cdot 10^{-3}$ & $1.9982$ & $1.0333 \cdot 10^{-3}$ & $1.9983$ & $1.5500 \cdot 10^{-3}$ & $1.9983$ \\ \hline
        $2 \cdot 10^{-4} \cdot 2^{-9}$ & $6.3778 \cdot 10^{-4}$ & $1.9991$ & $2.5848 \cdot 10^{-4}$ & $1.9992$ & $3.8772 \cdot 10^{-4}$ & $1.9992$ \\ \hline \hline
    \end{tabular}
    \begin{tabular}{c||c|c||c|c}
        $\Delta t$ & $\varphi_7$ errors & $\varphi_7$ orders & $\varphi_8$ errors & $\varphi_8$ orders  \\ \hline
        $2 \cdot 10^{-4}$ & $2.9078 \cdot 10^{1}$ & - & $2.4506 \cdot 10^{1}$ & - \\ \hline
        $2 \cdot 10^{-4} \cdot 2^{-1}$ & $8.2661 \cdot 10^{0}$ & $1.8146$ & $3.9029 \cdot 10^{0}$ & $2.6505$ \\ \hline
        $2 \cdot 10^{-4} \cdot 2^{-2}$ & $1.1812 \cdot 10^{0}$ & $2.8070$ & $2.7194 \cdot 10^{-1}$ & $3.8432$ \\ \hline
        $2 \cdot 10^{-4} \cdot 2^{-3}$ & $1.5476 \cdot 10^{-1}$ & $2.9321$ & $1.7704 \cdot 10^{-2}$ & $3.9411$ \\ \hline
        $2 \cdot 10^{-4} \cdot 2^{-4}$ & $1.9721 \cdot 10^{-2}$ & $2.9722$ & $1.1274 \cdot 10^{-3}$ & $3.9730$ \\ \hline
        $2 \cdot 10^{-4} \cdot 2^{-5}$ & $2.4881 \cdot 10^{-3}$ & $2.9867$ & $7.1124 \cdot 10^{-5}$ & $3.9865$ \\ \hline
        $2 \cdot 10^{-4} \cdot 2^{-6}$ & $3.1244 \cdot 10^{-4}$ & $2.9934$ & $4.4661 \cdot 10^{-6}$ & $3.9932$ \\ \hline
        $2 \cdot 10^{-4} \cdot 2^{-7}$ & $3.9144 \cdot 10^{-5}$ & $2.9967$ & $2.7957 \cdot 10^{-7}$ & $3.9977$ \\ \hline
        $2 \cdot 10^{-4} \cdot 2^{-8}$ & $4.8982 \cdot 10^{-6}$ & $2.9985$ & $1.7056 \cdot 10^{-8}$ & $4.0349$ \\ \hline
        $2 \cdot 10^{-4} \cdot 2^{-9}$ & $6.1178 \cdot 10^{-7}$ & $3.0012$ & $1.8974 \cdot 10^{-10}$ & $6.4901$ \\ \hline
    \end{tabular}}
    \caption{\rev{The errors and the corresponding orders of method NSSPMS(6,4) using different $\varphi$ functions while solving the logistic equation with $c=500$, plotted in Figure \ref{fig:phi_orders_c500}.}}
    \label{tab:phi_errors_c500}
\end{table}

\begin{table}
    \ssmall
    \centering
    \rev{\begin{tabular}{c||c|c||c|c||c|c}
         & \multicolumn{2}{c||}{NSSPMS(4,2)} & \multicolumn{2}{c||}{NSSPMS(4,3) with $\varphi_8$} & \multicolumn{2}{c}{NSSPMS(4,3) with $\varphi_7$}  \\ \hline
        $\Delta t$ & errors & orders & errors & orders & errors & orders  \\ \hline
        $2 \cdot 10^{-4}$ & $2.9338 \cdot 10^{0}$ & - & $2.9774 \cdot 10^{0}$ & -& $7.1664 \cdot 10^{0}$ & - \\ \hline
        $2 \cdot 10^{-4} \cdot 2^{-1}$ & $8.2578 \cdot 10^{-1}$ & $1.8289$ & $1.5806 \cdot 10^{-1}$ & $4.2355$ & $1.0250 \cdot 10^{0}$ & $2.8057$ \\ \hline
        $2 \cdot 10^{-4} \cdot 2^{-2}$ & $2.2313 \cdot 10^{-1}$ & $1.8879$ & $3.4056 \cdot 10^{-3}$ & $5.5364$ & $1.3489 \cdot 10^{-1}$ & $2.9257$  \\ \hline
        $2 \cdot 10^{-4} \cdot 2^{-3}$  & $5.8223 \cdot 10^{-2}$ & $1.9382$ & $6.7605 \cdot 10^{-4}$ & $2.3327$ & $1.7226 \cdot 10^{-2}$ & $2.9692$  \\ \hline
        $2 \cdot 10^{-4} \cdot 2^{-4}$ & $1.4886 \cdot 10^{-2}$ & $1.9677$ & $1.5597 \cdot 10^{-4}$ & $2.1159$ & $2.1757 \cdot 10^{-3}$ & $2.9850$  \\ \hline
        $2 \cdot 10^{-4} \cdot 2^{-5}$ & $3.7643 \cdot 10^{-3}$ & $1.9835$ & $2.4047 \cdot 10^{-5}$ & $2.6973$ & $2.7337 \cdot 10^{-4}$ &  $2.9925$ \\ \hline
        $2 \cdot 10^{-4} \cdot 2^{-6}$ & $9.4654 \cdot 10^{-4}$ & $1.9916$ & $3.2931 \cdot 10^{-6}$ & $2.8684$ & $3.4260 \cdot 10^{-5}$ & $2.9963$  \\ \hline
        $2 \cdot 10^{-4} \cdot 2^{-7}$ & $2.3733 \cdot 10^{-4}$ & $1.9958$ & $4.2969 \cdot 10^{-7}$ & $2.9381$ & $4.2880 \cdot 10^{-6}$ & $2.9981$  \\ \hline
        $2 \cdot 10^{-4} \cdot 2^{-8}$ & $5.9418 \cdot 10^{-5}$ & $1.9979$ & $5.4860 \cdot 10^{-8}$ & $2.9694$ & $5.3633 \cdot 10^{-7}$ & $2.9991$  \\ \hline
    \end{tabular}}
    \rev{\begin{tabular}{c||c|c||c|c||c|c}
         & \multicolumn{2}{c||}{NSSPRK(2,2)} & \multicolumn{2}{c||}{NSSPRK(3,3) with $\varphi_8$} & \multicolumn{2}{c}{NSSPRK(3,3) with $\varphi_7$}  \\ \hline
        $\Delta t$ & errors & orders & errors & orders & errors & orders  \\ \hline
        $2 \cdot 10^{-4}$ & $9.1774 \cdot 10^{-1}$ & - & $4.8634 \cdot 10^{-3}$ & -& $3.0723 \cdot 10^{-1}$ & -  \\ \hline
        $2 \cdot 10^{-4} \cdot 2^{-1}$ & $2.0672 \cdot 10^{-1}$ & $2.1503$ & $3.4867 \cdot 10^{-3}$ & $4.8010$ & $3.9080 \cdot 10^{-2}$ & $2.9749$  \\ \hline
        $2 \cdot 10^{-4} \cdot 2^{-2}$ & $4.9515 \cdot 10^{-2}$ & $2.0618$ & $6.1653 \cdot 10^{-4}$ & $2.4996$ & $4.9217 \cdot 10^{-3}$ & $2.9892$  \\ \hline
        $2 \cdot 10^{-4} \cdot 2^{-3}$  & $1.2152 \cdot 10^{-2}$ & $2.0267$ & $8.8392 \cdot 10^{-5}$ & $2.8022$ & $6.1740 \cdot 10^{-4}$ & $2.9949$  \\ \hline
        $2 \cdot 10^{-4} \cdot 2^{-4}$ & $3.0122 \cdot 10^{-3}$ & $2.0122$ & $1.1758 \cdot 10^{-5}$ & $2.9103$ & $7.7309 \cdot 10^{-5}$ &  $2.9975$ \\ \hline
        $2 \cdot 10^{-4} \cdot 2^{-5}$ & $7.5003 \cdot 10^{-4}$ & $2.0058$ & $1.5141 \cdot 10^{-6}$ & $2.9571$ & $ 9.6720 \cdot 10^{-6}$ & $2.9988$  \\ \hline
        $2 \cdot 10^{-4} \cdot 2^{-6}$ & $1.8714 \cdot 10^{-4}$ & $2.0028$ & $1.9205 \cdot 10^{-7}$ & $2.9789$ & $ 1.2095 \cdot 10^{-6}$ & $2.9994$  \\ \hline
        $2 \cdot 10^{-4} \cdot 2^{-7}$ & $4.6739 \cdot 10^{-5}$ & $2.0014$ & $2.4199 \cdot 10^{-8}$ & $2.9884$ & $ 1.5120 \cdot 10^{-7}$ & $2.9999$  \\ \hline
        $2 \cdot 10^{-4} \cdot 2^{-8}$ & $1.1679 \cdot 10^{-5}$ & $2.0007$ & $3.0757 \cdot 10^{-9}$ & $2.9760$ & $ 1.8859 \cdot 10^{-8}$ & $3.0031$  \\ \hline
    \end{tabular}}
    \rev{\begin{tabular}{c||c|c||c|c}
         & \multicolumn{2}{c||}{NSSPRK(6,4)} & \multicolumn{2}{c}{NSSPRK(10,4)}  \\ \hline
        $\Delta t$ & errors & orders & errors & orders \\ \hline
        $2 \cdot 10^{-4}$ & $2.4506 \cdot 10^{1}$ & - & $1.2566 \cdot 10^{-4}$ & - \\ \hline
        $2 \cdot 10^{-4} \cdot 2^{-1}$ & $3.9029 \cdot 10^{0}$ & $2.6505$ & $7.8481 \cdot 10^{-6}$ & $4.0010$ \\ \hline
        $2 \cdot 10^{-4} \cdot 2^{-2}$ & $2.7194 \cdot 10^{-1}$ & $3.8432$ & $4.8991 \cdot 10^{-7}$ & $4.0017$\\ \hline
        $2 \cdot 10^{-4} \cdot 2^{-3}$ & $1.7704 \cdot 10^{-2}$ & $3.9411$ & $3.0594 \cdot 10^{-8}$ & $4.0012$ \\ \hline
        $2 \cdot 10^{-4} \cdot 2^{-4}$ & $1.1274 \cdot 10^{-3}$ & $3.9730$ & $1.9088 \cdot 10^{-9}$ & $4.0025$\\ \hline
        $2 \cdot 10^{-4} \cdot 2^{-5}$ & $7.1124 \cdot 10^{-5}$ & $3.9865$ & $1.1539 \cdot 10^{-10}$ & $4.0481$ \\ \hline
        $2 \cdot 10^{-4} \cdot 2^{-6}$ & $4.4661 \cdot 10^{-6}$ & $3.9932$ & $2.2737 \cdot 10^{-13}$ & $8.9873$ \\ \hline
        $2 \cdot 10^{-4} \cdot 2^{-7}$ & $2.7957 \cdot 10^{-7}$ & $3.9977$ & $1.6257 \cdot 10^{-11}$ & -\\ \hline
        $2 \cdot 10^{-4} \cdot 2^{-8}$ & $1.7056 \cdot 10^{-8}$ & $4.0349$ & $2.6034 \cdot 10^{-11}$ & - \\ \hline
    \end{tabular}}
    \caption{\rev{The errors and the corresponding orders of the different multistep methods while solving the logistic equation with $c=500$, plotted in Figure \ref{fig:LMM_orders_c500}.}}
    \label{tab:LMM_errors_c500}
\end{table}

\begin{table}
    \centering
    \ssmall
    \rev{\begin{tabular}{c||c|c||c|c||c|c}
        $\Delta t$ & $\varphi_1$ errors & $\varphi_1$ orders & $\varphi_2$ errors & $\varphi_2$ orders & $\varphi_3$ errors & $\varphi_3$ orders  \\ \hline
        $0.1$ &  $5.3467 \cdot 10^{-1}$ & - & $5.1861 \cdot 10^{-1}$ & -& $6.1021 \cdot 10^{-1}$ & -  \\ \hline
        $0.1 \cdot 2^{-1}$ & $3.2527 \cdot 10^{-1}$ & $0.7170$ & $2.6609 \cdot 10^{-2}$ & $0.9627$ & $4.5905 \cdot 10^{-1}$ & $0.4106$ \\ \hline
        $0.1 \cdot 2^{-2}$ & $1.5828 \cdot 10^{-1}$ & $1.0391$ & $1.1823 \cdot 10^{-1}$ & $1.1703$ & $2.7619 \cdot 10^{-1}$ & $0.7330$ \\ \hline
        $0.1 \cdot 2^{-3}$  & $7.2915 \cdot 10^{-2}$ & $1.1181$ & $5.3255 \cdot 10^{-2}$ & $1.1506$ & $1.4172 \cdot 10^{-1}$ & $0.9625$ \\ \hline
        $0.1 \cdot 2^{-4}$ & $3.4220 \cdot 10^{-2}$ & $1.0913$ & $2.4977 \cdot 10^{-2}$ & $1.0923$ & $6.8670 \cdot 10^{-2}$ & $1.0453$ \\ \hline
        $0.1 \cdot 2^{-5}$ & $1.6476 \cdot 10^{-2}$ & $1.0544$ & $1.2060 \cdot 10^{-2}$ & $1.0503$ & $3.3198 \cdot 10^{-2}$ & $1.0486$ \\ \hline
        $0.1 \cdot 2^{-6}$ & $8.0718 \cdot 10^{-3}$ & $1.0294$ & $5.9220 \cdot 10^{-3}$ & $1.0262$ & $1.6230 \cdot 10^{-2}$ & $1.0324$ \\ \hline
        $0.1 \cdot 2^{-7}$ & $3.9934 \cdot 10^{-3}$ & $1.0153$ & $2.9338 \cdot 10^{-3}$ & $1.0133$ & $8.0116 \cdot 10^{-3}$ & $1.0185$ \\ \hline
        $0.1 \cdot 2^{-8}$ & $1.9859 \cdot 10^{-3}$ & $1.0078$ & $1.4601 \cdot 10^{-3}$ & $1.0067$ & $3.9785 \cdot 10^{-3}$ & $1.0099$ \\ \hline
        $0.1 \cdot 2^{-9}$ & $9.9029 \cdot 10^{-4}$ & $1.0039$ & $7.2833 \cdot 10^{-4}$ & $1.0034$ & $1.9823 \cdot 10^{-3}$ & $1.0051$ \\ \hline \hline 
    \end{tabular}
    \begin{tabular}{c||c|c||c|c||c|c}
        $\Delta t$ & $\varphi_4$ errors & $\varphi_4$ orders & $\varphi_5$ errors & $\varphi_5$ orders & $\varphi_6$ errors & $\varphi_6$ orders  \\ \hline
        $0.1$ & $5.6656 \cdot 10^{-1}$ & - & $5.1799 \cdot 10^{-1}$ & -& $5.3549 \cdot 10^{-1}$ & - \\ \hline
        $0.1 \cdot 2^{-1}$ & $3.4930 \cdot 10^{-1}$ & $0.6977$ & $2.4090 \cdot 10^{-1}$ & $1.1044$ & $2.8795 \cdot 10^{-1}$ & $0.8950$ \\ \hline
        $0.1 \cdot 2^{-2}$ & $1.3496 \cdot 10^{-1}$ & $1.3720$ & $6.5892 \cdot 10^{-2}$ & $1.8703$ & $9.2618 \cdot 10^{-2}$ & $1.6364$ \\ \hline
        $0.1 \cdot 2^{-3}$ & $3.8084 \cdot 10^{-2}$ & $1.8253$ & $1.6146 \cdot 10^{-2}$ & $2.0289$ & $2.3895 \cdot 10^{-2}$ & $1.9546$ \\ \hline
        $0.1 \cdot 2^{-4}$ & $9.7659 \cdot 10^{-3}$ & $1.9634$ & $4.0000 \cdot 10^{-3}$ & $2.0131$ & $5.9821 \cdot 10^{-3}$ & $1.9980$ \\ \hline
        $0.1 \cdot 2^{-5}$ & $2.4566 \cdot 10^{-3}$ & $1.9911$ & $9.9818 \cdot 10^{-4}$ & $2.0026$ & $1.4962 \cdot 10^{-3}$ & $1.9993$ \\ \hline
        $0.1 \cdot 2^{-6}$ & $6.1532 \cdot 10^{-4}$ & $1.9972$ & $2.4954 \cdot 10^{-4}$ & $2.0000$ & $3.7424 \cdot 10^{-4}$ & $1.9992$ \\ \hline
        $0.1 \cdot 2^{-7}$ & $1.5393 \cdot 10^{-4}$ & $1.9990$ & $6.2399 \cdot 10^{-5}$ & $1.9997$ & $9.3594 \cdot 10^{-5}$ & $1.9994$ \\ \hline
        $0.1 \cdot 2^{-8}$ & $3.8495 \cdot 10^{-5}$ & $1.9996$ & $1.5602 \cdot 10^{-5}$ & $1.9998$ & $2.3403 \cdot 10^{-5}$ & $1.9997$ \\ \hline
        $0.1 \cdot 2^{-9}$ & $9.6252 \cdot 10^{-6}$ & $1.9998$ & $3.9010 \cdot 10^{-6}$ & $1.9998$ & $5.8515 \cdot 10^{-6}$ & $1.9998$ \\ \hline \hline
    \end{tabular}
    \begin{tabular}{c||c|c||c|c}
        $\Delta t$ & $\varphi_7$ errors & $\varphi_7$ orders & $\varphi_8$ errors & $\varphi_8$ orders  \\ \hline
        $0.1$ & $5.2069 \cdot 10^{-1}$ & - & $5.1734 \cdot 10^{-1}$ & - \\ \hline
        $0.1 \cdot 2^{-1}$ & $2.3401 \cdot 10^{-1}$ & $1.1539$ & $2.1140 \cdot 10^{-1}$ & $1.2911$ \\ \hline
        $0.1 \cdot 2^{-2}$ & $4.4870 \cdot 10^{-2}$ & $2.9092$ & $2.5323 \cdot 10^{-2}$ & $3.8940$ \\ \hline
        $0.1 \cdot 2^{-3}$ & $5.9732 \cdot 10^{-3}$ & $2.9091$ & $1.7033 \cdot 10^{-3}$ & $3.8940$ \\ \hline
        $0.1 \cdot 2^{-4}$ & $7.5410 \cdot 10^{-4}$ & $2.9857$ & $1.0731 \cdot 10^{-4}$ & $3.9885$ \\ \hline
        $0.1 \cdot 2^{-5}$ & $9.4509 \cdot 10^{-5}$ & $2.9962$ & $6.7211 \cdot 10^{-6}$ & $3.9969$ \\ \hline
        $0.1 \cdot 2^{-6}$ & $1.1826 \cdot 10^{-5}$ & $2.9985$ & $4.2047 \cdot 10^{-7}$ & $3.9986$ \\ \hline
        $0.1 \cdot 2^{-7}$ & $1.4789 \cdot 10^{-6}$ & $2.9993$ & $2.6294 \cdot 10^{-8}$ & $3.9992$ \\ \hline
        $0.1 \cdot 2^{-8}$ & $1.8491 \cdot 10^{-7}$ & $2.9996$ & $1.6500 \cdot 10^{-9}$ & $3.9942$ \\ \hline
        $0.1 \cdot 2^{-9}$ & $2.3129 \cdot 10^{-8}$ & $2.9991$ & $1.1539 \cdot 10^{-10}$ & $3.8379$ \\ \hline
    \end{tabular}}
    \caption{\rev{The errors and the corresponding orders of method NSSPMS(6,4) using different $\varphi$ functions applied to the SEIR model with $\Pi=0$, plotted in Figure \ref{fig:seir_phiconv}.}}
    \label{tab:seir_phi_errors}
\end{table}

\begin{table}
    \centering
    \ssmall
    \rev{\begin{tabular}{c||c|c||c|c||c|c}
         & \multicolumn{2}{c||}{NSSPMS(4,2)} & \multicolumn{2}{c||}{NSSPMS(4,3) with $\varphi_8$} & \multicolumn{2}{c}{NSSPMS(4,3) with $\varphi_7$}  \\ \hline
        $\Delta t$ & errors & orders & errors & orders & errors & orders  \\ \hline
        $0.05$ & $2.4440 \cdot 10^{-3}$ & - & $2.0610 \cdot 10^{-2}$ & - & $3.4765 \cdot 10^{-2}$  &   -\\ \hline
        $0.05 \cdot 2^{-1}$ & $2.6209 \cdot 10^{-4}$ & $3.2210$ & $1.5549 \cdot 10^{-3}$ & $3.7284$ &  $5.4335 \cdot 10^{-3}$ & $2.6776$ \\ \hline
        $0.05 \cdot 2^{-2}$ & $3.6849 \cdot 10^{-5}$ & $2.8303$ & $9.9708 \cdot 10^{-5}$ & $3.9630$ &  $7.1272 \cdot 10^{-4}$ & $2.9305$ \\ \hline
        $0.05 \cdot 2^{-3}$ & $6.9473 \cdot 10^{-6}$ & $2.4071$ & $6.2038 \cdot 10^{-6}$ & $4.0065$ &  $9.0404 \cdot 10^{-5}$ & $2.9789$ \\ \hline
        $0.05 \cdot 2^{-4}$ & $1.7148 \cdot 10^{-6}$ & $2.0184$ & $3.7664 \cdot 10^{-7}$ & $4.0419$ &  $1.1370 \cdot 10^{-5}$ & $2.9912$ \\ \hline
        $0.05 \cdot 2^{-5}$ & $4.2660 \cdot 10^{-7}$ & $2.0071$ & $2.1926 \cdot 10^{-8}$ & $4.1025$ &  $1.4253 \cdot 10^{-6}$ & $2.9958$ \\ \hline
        $0.05 \cdot 2^{-6}$ & $1.0643 \cdot 10^{-7}$ & $2.0030$ & $1.1617 \cdot 10^{-9}$ & $4.2384$ &  $1.7842 \cdot 10^{-7}$ & $2.9979$ \\ \hline
        $0.05 \cdot 2^{-7}$ & $2.6581 \cdot 10^{-8}$ & $2.0014$ & $5.0217 \cdot 10^{-11}$ & $4.5318$ & $2.2318 \cdot 10^{-8}$   & $2.9990$ \\ \hline
        $0.05 \cdot 2^{-8}$ & $6.6424 \cdot 10^{-9}$ & $2.0006$ & $2.9863 \cdot 10^{-12}$ & $4.0717$ & $2.7908 \cdot 10^{-9}$   & $2.9995$ \\ \hline\hline
    \end{tabular}}
    \rev{\begin{tabular}{c||c|c||c|c||c|c}
         & \multicolumn{2}{c||}{NSSPRK(2,2)} & \multicolumn{2}{c||}{NSSPRK(3,3) with $\varphi_8$} & \multicolumn{2}{c}{NSSPRK(3,3) with $\varphi_7$}  \\ \hline
        $\Delta t$ & errors & orders & errors & orders & errors & orders  \\ \hline
        $0.05$ & $4.2992 \cdot 10^{-4}$ & - & $3.1523 \cdot 10^{-4}$ & -&  $1.7140 \cdot 10^{-3}$ & -  \\ \hline
        $0.05 \cdot 2^{-1}$ & $4.8459 \cdot 10^{-5}$ & $3.1492$ & $1.9170 \cdot 10^{-5}$ & $4.0395$ &  $2.1592 \cdot 10^{-4}$ & $2.9888$ \\ \hline
        $0.05 \cdot 2^{-2}$ & $1.1345 \cdot 10^{-5}$ & $2.0947$ & $1.1269 \cdot 10^{-6}$ & $4.0885$ &  $2.7019 \cdot 10^{-5}$ & $2.9984$ \\ \hline
        $0.05 \cdot 2^{-3}$ & $2.7744 \cdot 10^{-6}$ & $2.0318$ & $6.1496 \cdot 10^{-8}$ & $4.1957$ &  $3.3780 \cdot 10^{-6}$ & $2.9997$ \\ \hline
        $0.05 \cdot 2^{-4}$ & $6.8803 \cdot 10^{-6}$ & $2.0116$ & $2.7269 \cdot 10^{-9}$ & $4.4952$ &  $4.2228 \cdot 10^{-7}$ & $2.9999$ \\ \hline
        $0.05 \cdot 2^{-5}$ & $1.7145 \cdot 10^{-7}$ & $2.0047$ & $1.9683 \cdot 10^{-10}$ & $3.7923$ &  $5.2786 \cdot 10^{-8}$ & $3.0000$ \\ \hline
        $0.05 \cdot 2^{-6}$ & $4.2800 \cdot 10^{-8}$ & $2.0021$ & $2.7238 \cdot 10^{-11}$ & $2.8532$ &  $6.5982 \cdot 10^{-9}$ & $3.0000$ \\ \hline
        $0.05 \cdot 2^{-7}$ & $1.0693 \cdot 10^{-8}$ & $2.0010$ & $3.4914 \cdot 10^{-12}$ & $2.9637$ &  $8.2477 \cdot 10^{-10}$ & $3.0000$ \\ \hline
        $0.05 \cdot 2^{-8}$ & $2.6723 \cdot 10^{-9}$ & $2.0005$ & $5.0260 \cdot 10^{-13}$ & $2.7963$ &  $1.0307 \cdot 10^{-10}$ & $3.0000$ \\ \hline\hline
    \end{tabular}}
    \rev{\begin{tabular}{c||c|c||c|c}
         & \multicolumn{2}{c||}{NSSPRK(6,4)} & \multicolumn{2}{c}{NSSPRK(10,4)}  \\ \hline
        $\Delta t$ & errors & orders & errors & orders \\ \hline
        $0.05$ & $1.1739 \cdot 10^{-1}$ & - & $2.4392 \cdot 10^{-7}$ & - \\ \hline
        $0.05 \cdot 2^{-1}$ & $2.1800 \cdot 10^{-2}$ & $2.4289$ & $1.5246 \cdot 10^{-8}$ & $3.9999$  \\ \hline
        $0.05 \cdot 2^{-2}$ & $1.6444 \cdot 10^{-3}$ & $3.7287$ & $9.5298 \cdot 10^{-10}$ & $3.9999$ \\ \hline
        $0.05 \cdot 2^{-3}$ & $1.0584 \cdot 10^{-4}$ & $3.9576$ & $5.9574 \cdot 10^{-11}$ & $3.9997$ \\ \hline
        $0.05 \cdot 2^{-4}$ & $6.6819 \cdot 10^{-6}$ & $3.9855$ & $3.7333 \cdot 10^{-12}$ & $3.9962$ \\ \hline
        $0.05 \cdot 2^{-5}$ & $4.1960 \cdot 10^{-7}$ & $3.9932$ & $2.4303 \cdot 10^{-13}$ & $3.9413$  \\ \hline
        $0.05 \cdot 2^{-6}$ & $2.6286 \cdot 10^{-8}$ & $3.9966$ & $5.3291 \cdot 10^{-14}$ & $2.1892$  \\ \hline
        $0.05 \cdot 2^{-7}$ & $1.6446 \cdot 10^{-9}$ & $3.9985$ & $6.7391 \cdot 10^{-14}$ & - \\ \hline
        $0.05 \cdot 2^{-8}$ & $1.0248 \cdot 10^{-10}$ & $4.0043$ & $1.0170 \cdot 10^{-13}$ & - \\ \hline
    \end{tabular}}
    \caption{\rev{The errors and the corresponding orders of the different nonstandard methods applied to the SEIR model with $\Pi=0$, plotted in Figure \ref{fig:seir_meth}.}}
    \label{tab:seir_method_errors}
\end{table}








\end{document}